\documentclass[11pt,reqno]{amsart}

\usepackage[T1]{fontenc}
\usepackage[utf8]{inputenc}

\usepackage[top=1in,bottom=1in,left=1in,right=1in]{geometry}
\usepackage{amsmath,amssymb,amsthm,amsfonts}
\usepackage{mathtools}
\usepackage{mathrsfs}
\usepackage[mathscr]{euscript}
\usepackage{fixmath}
\usepackage{bm}
\usepackage{upgreek}
\usepackage{array}
\usepackage{tipa}
\usepackage{enumitem}
\usepackage{graphicx}
\usepackage{epstopdf}
\usepackage{svg}

\usepackage{xcolor}
\usepackage[textsize=small]{todonotes}
\usepackage{comment}
\usepackage{titlecaps}
\usepackage{lipsum}

\usepackage{hyperref}
\usepackage{cleveref}

\newcommand*{\wt}[1]{\widetilde{#1}}
\newcommand{\Teich}{\mathcal{T}}

\newcommand{\R}{\mathbb R}

\newcommand{\PC}{\mathbb{P}\mathcal{C}}
\newcommand{\scc}{\mathrm{scc}}
\newcommand{\cc}{\mathrm{cc}}

\newcommand{\Dil}{\operatorname{Dil}}
\newcommand{\Syst}{\operatorname{Syst}}

\newcommand{\bt}[1]{\boldsymbol{#1}}

\newcommand{\rr}{\mathbb{R}}

\newcommand{\lieg}{\mathfrak{g}}
\newcommand{\liea}{\mathfrak{a}}
\newcommand{\liep}{\mathfrak{p}}
\newcommand{\m}{\mathfrak{m}}
\newcommand{\Int}{\mathrm{int}}

\newcommand{\asf}{\mathsf{A}}
\newcommand{\g}{\mathsf{G}}      
\newcommand{\ko}{\mathsf{K}}     
\newcommand{\w}{\mathsf{W}}      
\newcommand{\p}{\mathsf{P}}      
\newcommand{\f}{\mathscr{F}}    
\newcommand{\nsf}{\mathsf{N}} 

\newcommand{\bg}{\partial\Gamma}
\newcommand{\bgs}{\partial^2\Gamma}
\newcommand{\gh}{\Gamma'}

\newcommand{\ha}{\mathscr{X}_\Theta(\Gamma, \g)}   

\newcommand{\zcone}{\mathscr{L}_r^{\Theta}}
\newcommand{\dzcone}{(\mathscr{L}_r^{\Theta})^*}

\makeatletter
\newtheorem*{rep@theorem}{\rep@title}
\newcommand{\newreptheorem}[2]{%
  \newenvironment{rep#1}[1]{%
    \def\rep@title{#2 \ref{##1}}%
    \begin{rep@theorem}}%
  {\end{rep@theorem}}}
\makeatother

\makeatletter
\newtheorem*{rep@corollary}{\rep@title}
\newcommand{\newrepcorollary}[2]{%
  \newenvironment{rep#1}[1]{%
    \def\rep@title{#2 \ref{##1}}%
    \begin{rep@corollary}}%
  {\end{rep@corollary}}}
\makeatother

\newtheorem{theorem}{Theorem}[section]
\newtheorem*{theorem*}{Theorem*}
\newtheorem{corollary}[theorem]{Corollary}
\newtheorem{lemma}[theorem]{Lemma}
\newtheorem{example}[theorem]{Example}
\newtheorem{proposition}[theorem]{Proposition}

\newtheorem{remark}[theorem]{Remark}

\newtheorem{definition}[theorem]{Definition}
\newtheorem{question}[theorem]{Question}
\newreptheorem{theorem}{Theorem}
\newrepcorollary{corollary}{Corollary}
\newtheorem*{keyproperties}{Key Properties}

\newcommand{\bgamma}{\boldsymbol{\gamma}}
\newcommand{\ba}{\boldsymbol{a}}
\newcommand{\bb}{\boldsymbol{b}}
\newcommand{\bc}{\boldsymbol{c}}
\newcommand{\bxi}{\boldsymbol{\xi}}
\newcommand{\boldeta}{\boldsymbol{\eta}}
\newcommand{\bx}{\boldsymbol{x}}
\newcommand{\by}{\boldsymbol{y}}
\newcommand{\bz}{\boldsymbol{z}}
\newcommand{\bmu}{\boldsymbol{\mu}}
\newcommand{\bnu}{\boldsymbol{\nu}}
\newcommand{\blambda}{\boldsymbol{\lambda}}

\crefname{theorem}{Theorem}{Theorems}
\crefname{lemma}{Lemma}{Lemmas}
\crefname{proposition}{Proposition}{Propositions}
\crefname{corollary}{Corollary}{Corollarys}
\crefname{section}{Section}{Sections}
\crefname{definition}{Definition}{Definitions}
\crefname{question}{question}{Questions}
\crefname{keyproperties}{key property}{key properties}
\Crefname{keyproperties}{Key Property}{Key Properties}

\title{Horoboundary and rigidity of filling geodesic currents}
\author[Meenakshy Jyothis]{Meenakshy Jyothis}
\address{Department of Mathematics\\
         University of Oklahoma\\
         601 Elm Ave Room 423, Norman, OK 73019\\
         USA}
\email{mjyothis@ou.edu}

\author[Martínez-Granado]{Dídac Martínez-Granado}
\address{Department of Mathematics\\University of Luxembourg\\Av. de la Fonte 6, Esch-sur-Alzette, L-4364, Luxembourg}       
\email{didac.martinezgranado@uni.lu}

\date{}

\begin{document}

\begin{abstract}
We endow the space of projective filling geodesic currents on a closed hyperbolic surface with a natural asymmetric metric extending Thurston’s asymmetric metric on Teichmüller space, as well as analogous metrics arising from Hitchin representations. More generally, we show that this metric extends beyond surface groups and geodesic currents, and encompasses metrics associated with Anosov representations of Gromov hyperbolic groups. We identify the horofunction compactification of the space of projective filling currents equipped with this metric with the space of projective geodesic currents. As a consequence, we obtain a rigidity result: the metric spaces of projective filling geodesic currents associated with closed surfaces of distinct genera are not isometric.
\end{abstract}

\maketitle
\section{Introduction}  

Let $S$ be a closed orientable surface of genus at least two. The Teichmüller space $\Teich(S)$ is the space of isotopy classes of marked hyperbolic structures on $S$. Thurston's asymmetric distance $d_{Th}$~\cite{Th1998} is one of the standard metric structures on $\mathcal{T}(S)$, and measures the maximal curve length ratio between two hyperbolic surfaces. We extend the Thurston metric beyond $\Teich(S)$, which is finite dimensional, to two larger infinite dimensional spaces: the \emph{space of projective filling geodesic currents} $\PC_{\mathrm{fill}}(S)$ on $S$, and the \emph{space of metric potentials} $\mathcal{H}_\Gamma^{++}$ defined more generally for any non-elementary Gromov hyperbolic group $\Gamma$~\cite{CT25:Manhattan,CRS24:Joint}. These spaces encompass many classes of geometric structures on $S$ beyond hyperbolic metrics. When $\Gamma = \pi_{1}(S)$, we have the following chain of inclusions:
\[
\Teich(S) \subset \PC_{\mathrm{fill}}(S) \subset \mathcal{H}_\Gamma^{++}.
\]
The Thurston metric $d_{Th}$ on $\Teich(S)$ is obtained as the restriction of an extended metric $d$ on $\mathbb{P}\mathcal{C}_{\mathrm{fill}}$, which itself arises as the restriction of a metric defined on $\mathcal{H}_\Gamma^{++}$.
Since $\mathcal{H}_\Gamma^{++}$ is defined for a general non-elementary Gromov hyperbolic group $\Gamma$, the extended metric $d$ induces metrics on spaces of geometric structures not inherently associated with surfaces, such as \emph{cubulations} and \emph{Anosov representations} of hyperbolic groups. However, a significant part of this paper focuses on the case where $\Gamma = \pi_{1}(S)$, and we study the extended Thurston metric defined on $\mathbb{P}\mathcal{C}_{\mathrm{fill}}(S)$.

The space of \emph{geodesic currents} $\mathcal{C}(S)$ on a surface $S$ was first defined by Bonahon as an infinite-dimensional cone of positive measures on the space of geodesics of $S$  \cite{Bon86}. The space of geodesic currents provides a unified framework for studying various objects of geometric interest associated with $S$. For instance, it contains free homotopy classes of closed curves on $S$~\cite[Proposition~4.4]{Bon86}, the Teichmüller space of $S$ \cite[Theorem~12]{Bon88}, the moduli space of marked flat structures on $S$ \cite[Theorem~1]{DLR2010}, and the negatively curved Riemannian metrics on $S$~\cite[Th\'eor\`eme~1] {Otal90:SpectreMarqueNegative}. Furthermore, the geometric intersection number between closed curves extends to an intersection form $i \colon \mathcal{C} \times \mathcal{C} \to \R$ on the space of currents.

The interior of the cone $\mathcal{C}(S)$ is the space of \emph{filling geodesic currents} $\mathcal{C}_{\mathrm{fill}}(S)$. It is known that $\Teich(S)$ embeds into $\mathcal{C}_{\mathrm{fill}}(S)$ as well as into its projectivization $\mathbb{P}\mathcal{C}_{\mathrm{fill}}(S)$ \cite{Bon88}. For any $x$, $y \in \mathbb{P}\mathcal{C}_{\mathrm{fill}}(S)$ the extended Thurston metric $d$ on $\PC_{\mathrm{fill}}$ is defined as follows:
\begin{equation}
    d(x,y) \coloneqq \log \left( \sup_{\bc \in \cc} 
    \frac{i(\by,\bc)}{i(\bx,\bc)} \cdot \frac{h(\by)}{h(\bx)}\right).
\end{equation}
Here, $\bx$ and $\by$ are representatives of $x$ and $y$ in the space of currents. The function $h \colon \mathcal{C}_{\mathrm{fill}} (S)\rightarrow \R$ is the entropy of currents, measuring the exponential growth of closed curves as a function of their intersection with a given filling current. A more detailed exposition of this metric is provided in Section~\ref{subsec:intro_metric}.

We prove that for two distinct closed surfaces $S$ and $S'$, the metric spaces $(\PC_{\mathrm{fill}}(S), d)$ and ($\PC_{\mathrm{fill}}(S'),d)$ cannot be isometric. To prove this rigidity result, we study the horofunction compactification of the metric space
$(\PC_{\mathrm{fill}}(S), d)$.
Introduced by Gromov, this compactification encodes asymptotic information about a metric space
in a way that reflects its underlying metric structure~\cite{Gro}.
We show that, in our setting, it is homeomorphic to the space of projective geodesic currents
$\PC(S)$, the compact space obtained by considering $\mathcal{C}(S)$ up to scaling by positive
real numbers. 

In what follows we provide a more detailed overview of each of these results.

\subsection{The geometry of projective filling currents}
\label{subsec:intro_projective_filling_currents} 
Over the past few decades, geodesic currents have become an important tool in low-dimensional geometry and dynamics. They have played a crucial role in curve counting results and in higher Teichmüller theory—the study of higher-rank analogues of convex co-compact representations in rank one (see e.g.~\cite{RS19, ES22:GeodesicCount, BridgemanCanaryLabourieSambarino18:SimpleRoots, MZ19:PositivelyRatioed, BIPP21, BIPP24}).  They were also used to construct topological compactification of several marked moduli structures associated to a surface. This includes Thurston’s compactification of Teichmüller space~\cite[Theorem~18]{Bon88}, the Duchin–Leininger–Rafi compactification of spaces of singular Euclidean structures~\cite[Theorem~5]{DLR2010}, and, more recently, compactifications of higher-rank analogues of Teichmüller theory~\cite{MOT24:Ball, OT24:Blaschke, OT23:SO23}. All of these compactifications were constructed by embedding the relevant spaces into the space of projective currents $\mathbb{P}\mathcal{C}(S)$.



The space of projective filling currents $\mathbb{P}\mathcal{C}_{\mathrm{fill}}(S)$
is an open, dense, and locally compact subspace of $\mathbb{P}\mathcal{C}(S)$. It is precisely the projectivization of the interior of the cone $\mathcal{C}(S)$ in the sense that a geodesic current is filling if and only if it has strictly positive intersection with every nonzero geodesic current ~\cite[Theorem~1.3]{BIPP21}. We study $\mathbb{P}\mathcal{C}_{\mathrm{fill}}(S)$ when equipped with the extended Thurston metric $d$ defined in \cref{eq:d}.

\subsubsection{The Horofunction compactification and the horofunction boundary}
Our first result characterizes the horofunction compactification of projective filling geodesic currents. The horofunction compactification is a topological compactification of a metric space that behaves well with respect to the underlying metric. For instance, isometries of a metric space extend to homeomorphisms of its horofunction compactification. It has found applications in the study of asymptotic properties of random walks on weakly hyperbolic groups~\cite{MT18:RandomWalks} and in determining the isometry group of Hilbert geometries~\cite{LW11:Hilbert}.

We show that the horofunction compactification of $\PC_{\mathrm{fill}}(S)$ endowed with the metric $d$ is homeomorphic to $\PC$.

\begin{theorem}\label{horoboundary}
The horofunction compactification of $\mathbb{P}\mathcal{C}_{\mathrm{fill}}(S)$ equipped with the extended Thurston metric $d$ is homeomorphic to the space of projective currents $\mathbb{P}\mathcal{C}(S)$. Moreover, the horofunction boundary corresponds precisely to the space of projective non-filling currents.
\end{theorem}

The horofunction compactification of Teichmüller space has been studied for several of its natural distances. For example, for the Thurston metric it was shown by Walsh~\cite{CW2015} that the horofunction boundary of $\mathcal{T}(S)$ coincides with the space of projective measured laminations $\mathbb{P}\mathcal{ML}(S)$. The horofunction compactification of $(\mathbb{P}\mathcal{C}_{\mathrm{fill}}(S), d)$ constructed in Theorem~\ref{horoboundary} extends Walsh’s horofunction compactification of $(\mathcal{T}(S), d_{Th})$.
The horofunction compactification with respect to the Teichmüller metric on $\mathcal{T}(S)$ is also well understood; see ~\cite{LS14:Horoboundary_Teich,A24:HoroboundaryQualitative}. Other relevant metric spaces whose horofunction compactification has been studied include non-compact symmetric spaces equipped with polyhedral Finsler metrics~\cite{KL18:Finsler}, Hilbert geometry in both finite and infinite dimensions~\cite{Co14:FiniteHilbert,Cor18:HilbertInfinite}, Hadamard manifolds and their quotients ~\cite{DPS12:HoroboundaryNegative}, and the Heisenberg group with the Carnot--Carathéodory metric~\cite{KN10:HeisenbergGroup}.

The construction of the horofunction compactification for $(\mathbb{P}\mathcal{C}_{\mathrm{fill}}(S),d)$ builds on Walsh's framework for the Thurston metric $d_{Th}$ on $\Teich(S)$ ~\cite{CW2015}. However, the passage to the setting of geodesic currents brings genuinely new challenges. In particular, analyzing the behavior of the extension of an isometry 
$\phi \colon \mathbb{P}\mathcal{C}_{\mathrm{fill}}(S) \to \mathbb{P}\mathcal{C}_{\mathrm{fill}}(S)$ to a homeomorphism of the horofunction boundary $\mathbb{P}\mathcal{C}_{\mathrm{fill}}(\infty)$ is considerably more delicate, since boundary points of the space of geodesic currents exhibit a far greater variety of behaviors than in the case of Teichmüller space.

\subsubsection{Isometric Rigidity} As a consequence of \cref{horoboundary} we prove the following rigidity result.

\begin{theorem}\label{isometricrigidity}
Let $S_{g}$ and $S_{g'}$ be closed orientable surfaces of genus at least two, with $g \neq g'$. Then the spaces $(\mathbb{P}\mathcal{C}_{\mathrm{fill}}(S_{g}),d)$ and $(\mathbb{P}\mathcal{C}_{\mathrm{fill}}(S_{g'}),d')$, endowed with their respective asymmetric Thurston metrics $d$ and $d'$, are not isometric.
\end{theorem}

Unlike Teichmüller space, $\mathbb{P}\mathcal{C}_{\mathrm{fill}}(S)$ is infinite dimensional, making its isometry type difficult to distinguish. In related settings, such as infinite-dimensional cones with Hilbert metrics~\cite{Cor18:HilbertInfinite}, deciding the isometry types is a delicate problem.

Furthermore, it follows from classical results in Choquet theory that the spaces of geodesic currents associated to distinct closed hyperbolic surfaces are all homeomorphic~(see Proposition~\ref{prop:geodesiccurrents_homeo}). Hence, to distinguish space of currents arising from distinct closed surfaces we need a finer structure. Theorem~\ref{isometricrigidity} shows that the isometry type of the metric space $(\PC_{\mathrm{fill}}, d)$ is uniquely determined by the genus of the surface. It is not yet known whether the subspaces of filling currents of surfaces of distinct genera can be homeomorphic (non-equivariantly with respect to the action of the mapping class group).


The asymmetric metric $d$ naturally gives analogues of asymmetric metrics on several geometric structures on $S$ which can be realized as geodesic currents (Section~\ref{sec:examples}), 
including normalized negatively curved Riemannian metrics, non-positively curved metrics, and Hitchin representations of surface groups. It would be interesting to understand the geometry of these spaces with respect to this metric. We now discuss the extended Thurston metric in the context of \emph{metric potentials} $\mathcal{H}^{++}_\Gamma$ on a hyperbolic group.

\subsection{Extending the metric beyond geodesic currents: Anosov representations}\label{subsec:further_extensions} 

Beyond geodesic currents, the metric $d$ extends to more general settings, such as (normalized) \emph{cubulations} (see~Section~\ref{ex:cubulation}) and $\Theta$-\emph{Anosov representations} of a Gromov hyperbolic group $\Gamma$ (see Section~\ref{sec:examples}). In fact, these metrics arise from an extension of the asymmetric metric $d$ to the broader space of \emph{metric potentials} $\mathcal{H}^{++}_\Gamma$, in the sense of~\cite{CT25:Manhattan,CRS24:Joint} (see Appendix~\ref{sec:asymmetric} for precise definitions), which we now briefly discuss. 

For a Gromov hyperbolic group let $\Gamma$, $D_\Gamma$ denote \emph{the space of pseudo-metrics} on $\Gamma$. It consists of left-invariant, Gromov hyperbolic  \emph{pseudo-metrics} which are quasi-isometric to a word-metric on $\Gamma$ ~\cite{Fur02:Coarse,Eduardo}. A pseudo-metric is a non-negative function $\Gamma \times \Gamma \to \mathbb{R}$ satisfying symmetry and triangular inequality.  Examples include any word-metric on $\Gamma$ and distances arising from geometric actions of $\Gamma$ on a Gromov hyperbolic geodesic space. We say two pseudo-metrics are \emph{roughly similar} if they differ by a bounded amount after rescaling. The \emph{space of metric structures} $\mathcal{D}_{\Gamma}$ is the quotient of $D_{\Gamma}$ by the equivalence relation of rough similarity.

When $\Gamma=\pi_1(S)$, the space $\mathcal{D}_\Gamma$ contains $\PC_{\mathrm{fill}}(S)$~\cite{StephenEduardo,LucDid}, which already covers many interesting examples discussed above. However, even in the surface setting, $\mathcal{D}_\Gamma$ contains pseudo-metrics which do not arise as geodesic currents, such as quasi-fuchsian representations or cubulations~\cite{FF22:QF,BR24:Approximating}.
The \emph{set of proper metric potentials} of $\Gamma$, denoted by $H^{++}_\Gamma$ is a further generalization of $D_\Gamma$ (introduced in~\cite{CRS24:Joint} building on~\cite{CT25:Manhattan}), which lifts the hypotheses of symmetry and non-negativity of the pseudo-metrics. A proper metric potential has a well-defined notion of \emph{entropy} $h(\psi)$, measuring the exponential growth rate of conjugacy classes of $[\Gamma]$, as well as of \emph{stable length} function, $\ell_\psi \colon [\Gamma] \to \mathbb{R}$, assigning to each $[g]$ the invariant $\lim_n \psi(o,g^n)/n$. For proper metric potentials, entropy and stable length are strictly positive. The \emph{space of proper metric potential}, $\mathcal{H}^{++}_{\Gamma}$ is obtained by considering $H^{++}_{\Gamma}$ up to the equivalence relation of rough similarity.

  \begin{figure}[ht!]
    \centering
  \fontsize{9pt}{9pt}\selectfont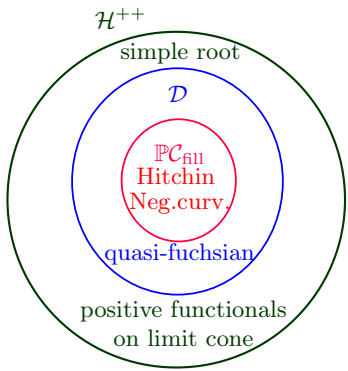
    \caption{Diagram represents the space of proper hyperbolic potentials: When $\Gamma=\pi_1(S)$, we have inclusions $\mathcal{C}_\mathrm{fill}(S) \subset D_\Gamma \subset H^{++}_\Gamma$, which induce inclusions $\PC_\mathrm{fill}(S) \subset \mathcal{D}_\Gamma \subset \mathcal{H}^{++}_\Gamma$ on the respective quotient spaces. Hitchin representations and negatively curved Riemannian metrics are contained in $\mathbb{P}\mathcal{C}_\mathrm{fill}(S)$. Quasi-fuchsian representations are in $\mathcal{D}_{\Gamma}$, but not in $\mathbb{P}\mathcal{C}_\mathrm{fill}(S)$. Furthermore, potentials coming from simple roots or, more generally, positive functionals on the Benoist limit cone, are in $\mathcal{H}^{++}_{\Gamma}$ but not in $\mathcal{D}_\Gamma$. The last statement is part of the content of Theorem~\ref{extensions}.}
    \label{fig:venn}
  \end{figure}
Natural examples of proper hyperbolic functionals not necessarily lying in $\mathcal{D}_\Gamma$ arise from Anosov representation theory.
Let $\mathsf{G}$ be a connected semisimple Lie group and $\Gamma$ a Gromov hyperbolic group. Anosov representations $r\colon \Gamma \to \mathsf{G}$ are indexed by the subsets of simple roots $\Theta$ of a Cartan subspace $\mathfrak{a}$ of $\lieg$, and generalize convex cocompact representations in rank one.
Each such representation $r$ determines a closed subset in the interior of the Cartan subalgebra $\mathfrak{a}^+$ called the Benoist limit cone $\zcone\subset \mathfrak{a}^+$. This cone encodes the asymptotic Cartan projections of $r(\Gamma)$ and records the growth directions of the representation.
Any linear functional $\varphi\colon \mathfrak{a}^+ \to \mathbb{R}$ that is positive on $\zcone\backslash \{0\}$, i.e, $\varphi \in \text{int} (\dzcone)$, yields a conjugacy-invariant function $\ell_r^\varphi \colon \Gamma \to \mathbb{R}$. This function plays the role of $\varphi$-length of $r$. 
Let $\ha$ denote the set of $\Theta$-Anosov representations from $\Gamma$ to $\g$, and fix a functional $\varphi \in \bigcap_{{r}\in\ha} \text{int} (\dzcone).$
Let $\Gamma' \subset \Gamma$ denote the subset of infinite order elements of $\Gamma$, and $[\Gamma']$ its set of conjugacy classes.
In~\cite[Theorem~1.1]{CDPW24:Asymmetric}, for every $\varphi$, the following asymmetric metric is constructed on $\ha$:
\[
d^{\varphi}(r_1,r_2) \coloneqq \sup_{[g] \in [\Gamma']} \frac{h_{r_2}^\varphi \ell_{r_2}^\varphi(g)}{h_{r_1}^\varphi \ell_{r_1}^\varphi(g)}
\]
where $h_r^{\varphi}$ denotes the exponential growth of conjugacy classes $[\Gamma]$ measured with respect to $\ell_r^\varphi$.

We show that $d$ restricts to this metric.

\begin{theorem}\label{extensions}
Given a non-elementary Gromov hyperbolic group $\Gamma$, there exists an extension of the asymmetric metric $d$ to the space of metric potentials $\mathcal{H}_\Gamma^{++}$, given by
\[
d(\psi_1,\psi_2) \coloneqq \log \left( \sup_{[g] \in [\Gamma']}\frac{h(\psi_2) \ell_{\psi_2}([g])}{h(\psi_1) \ell_{\psi_1}([g])} \right).
\]
Moreover, for every $\Theta$-Anosov representation, and every $\varphi \in \text{int} (\dzcone)$ there exists $\psi_{\varphi,r} \in H_\Gamma^{++}$
so that
\[
\ell_{\psi_{r,\varphi}}(g)\coloneqq \lim_n \frac{\psi(o,g^n)}{n} = \ell_r^\varphi(g)
\]
for every $g \in \Gamma$.
In particular, if $\varphi \in \bigcap_{{r}\in\ha} \text{int} (\dzcone),$ then
\[
d(\psi_{r_1,\varphi},\psi_{r_2,\varphi})=d^\varphi(r_1,r_2),
\]
i.e., $d$ extends the metric $d^\varphi$.
\end{theorem}

The content of the above theorem is two-fold:
\begin{enumerate} \item The observation that $d$ can be extended to $\mathcal{H}_\Gamma^{++}$. This follows from an adaptation of arguments in~\cite{CT25:Manhattan}), see Proposition~\ref{prop:asymmetric_defined}.
\item The fact that $\Theta$-Anosov representations yield examples of potentials for any functional in the interior of the Benoist limit cone (see Section~\ref{sec:examples} for definitions), generalizing~\cite[Lemma~7.1]{CT24:Invariant}, and subsuming~\cite[Theorem~1.1]{CDPW24:Asymmetric}.
This result is Proposition~\ref{prop:alpha_metric}. We emphasize that, unlike the construction in~\cite{CDPW24:Asymmetric}, the definition of the extended metric $d$ does \emph{not} rely on the existence of H\"older reparameterizations of the Mineyev flow~\cite{Mineyev}, which only exist in specific settings~\cite{S24:Dichotomy}. In forthcoming work~\cite{CMGR26:GreenMetrics}, the second author and collaborators show that these also induce points in $H_\Gamma^{++}$.
\end{enumerate}
See Appendix~\ref{sec:asymmetric} for a detailed discussion on $\mathcal{D}_\Gamma$ and $\mathcal{H}^{++}_\Gamma$, and how the metric $d$ is defined on these spaces.
See Section~\ref{sec:examples} for the examples discussed in the statement.

Coupled with marked length spectrum rigidity results in the literature, the metric $d$ provides an extension of the classical asymmetric Thurston metric to many other geometric contexts: area-one metrics, negatively curved and non-positively curved metrics or cubulations of entropy 1  (see Section~\ref{sec:examples} and Appendix~\ref{ex:cubulation}).  

\subsection{The extended Thurston metric on projective filling currents}
\label{subsec:intro_metric}
For any two marked hyperbolic structures $x,y \in \mathcal{T}(S)$, Thurston's original metric $d_{Th}(x,y)$ is defined as the logarithm of the minimal Lipschitz constant of a map $f:(S,x) \to (S,y)$ homotopic to the change of marking~\cite{Th1998}. 
For a recent account of this metric, see~\cite{PanSu2024GeometryThurstonMetric}. 
Thurston also showed that $d_{Th}$ admits the following equivalent characterization: for $x \in \mathcal{T}(S)$ and a closed curve $c$ on $S$, let $\ell_{x}(c)$ denote the length of the geodesic representative of $c$ with respect to $x$. Then for $x, y \in \mathcal{T}(S)$, 
\begin{equation}
d_{Th}(x,y) = \log \sup_{c \in \scc} \frac{\ell_y(c)}{\ell_x(c)}, 
\label{eq:dTh}
\end{equation}
where the supremum is taken over all simple closed curves $c \in \scc$.

The metric $d_{Th}$ is asymmetric. Both the space of free homotopy classes of closed curves on $S$ and $\mathcal{T}(S)$ embed into $\mathcal{C}(S)$, and into its projectivization $\PC(S)$. Moreover, the geometric intersection number between closed curves extends to a continuous bilinear form 
\[
i: \mathcal{C}(S) \times \mathcal{C}(S) \to \mathbb{R}_{+}
\]
on geodesic currents~\cite{Bon86}. In particular, for any $x \in \mathcal{T}(S)$ and closed curve $c$ on $S$, one has $i(x,c) = \ell_x(c)$ (see Section~\ref{subsec:currents} for more details).

A natural extension of \cref{eq:dTh} is obtained by replacing geodesic length with the intersection form and taking the supremum over all closed curves. To ensure that the resulting ratio is well defined, we normalize by the topological entropy of currents. Precisely, for $x,y \in \mathbb{P}\mathcal{C}_{\mathrm{fill}}(S)$, choose representatives $\bx,\by$ in the space of filling currents and set
\begin{equation}
    d(x,y) \coloneqq \log \left( \sup_{\bc \in \cc} 
    \frac{i(\by,\bc)}{i(\bx,\bc)} \cdot \frac{h(\by)}{h(\bx)}\right),
\label{eq:d}
\end{equation}
where the supremum is taken over all closed curves $\bc \in \cc$, and the entropy $h:\mathcal{C}(S)\to\mathbb{R}_{+}$ is defined by
\[
h(\bx) \coloneqq \limsup_{n\to\infty} \frac{1}{n}
\log \#\{\bc \in \cc : i(\bx,\bc) < n \}.
\]

In Proposition~\ref{prop:asymmetric_geodesics}, we show that $d$ is a well-defined metric on $\mathbb{P}\mathcal{C}_{\mathrm{fill}}(S)$. 
Sapir recently extended a \emph{symmetrized} version of the Thurston metric to the space of projective filling currents~\cite{Jenya}. 
Reyes~\cite[Remark~3.8]{Eduardo} observed that one can extend $d_{Th}$ to $\mathcal{D}_\Gamma$. In Appendix~\ref{sec:asymmetric}, we observe that the asymmetric distance can be further extended to the space of \emph{metric potentials} $\mathcal{H}^{++}_\Gamma$.
This builds on~\cite{CT25:Manhattan, StephenEduardo, LucDid, CRS24:Joint}.
As mentioned above, this metric generalizes asymmetric metrics previously introduced in the setting of $\Theta$-Anosov representations~\cite{CDPW24:Asymmetric}, which in turn generalize analogous asymmetric metrics arising in the context of convex cocompact representations of hyperbolic manifold groups into $\operatorname{PO}(n,1)$~\cite{Bur93:Manhattan}.

Several related metrics have appeared in the literature that lie outside the scope of the present extension. These include the asymmetric metrics on Culler--Vogtmann Outer space~\cite{FrancavigliaMartino2011}, which use a different normalization (by volume);  metrics on the space of flat structures~\cite{Shi25:Thurston}, where normalization is by area rather than entropy; and extensions of asymmetric metrics to geometrically finite (but not convex cocompact) hyperbolic manifolds~\cite{GK11:GeometricallyFinite}. More recently, asymmetric metrics have been extended to representations of hyperbolic groups into non-semisimple Lie groups, notably in the setting of Margulis spacetimes~\cite{GM25:MargulisSpacetimes}. This construction builds on the asymmetric metric introduced in~\cite{CDPW24:Asymmetric}, arising from the theory of H\"older reparameterizations of the Mineyev flow developed in~\cite{Sam15:Orbital,Quantitative,S24:Dichotomy,BridgemanCanaryLabourieSambarino18:SimpleRoots}.

In forthcoming work, the second author and collaborators~\cite{CMGR26:GreenMetrics} show that these flow reparameterizations induce points in $\mathcal{H}_\Gamma^{++}$, and hence are also subsumed by the extended metric considered here.

\subsection{Notational Convention} From here on, for simplicity, we use the following notation:
\[
\PC_{\mathrm{fill}}(S) = \PC_{\mathrm{fill}} \quad \PC(S) = \mathbb{P}\mathcal{C}, \text{ and } \mathcal{T}(S) = \mathcal{T}
\]
When we discuss spaces associated to more than one surface, as in \cref{isometricrigidity}, we will include the surface in the notation.

As described in \cref{eq:d} in order to define $d(x,y)$ for any $x, y \in \mathbb{P}\mathcal{C}_{\mathrm{fill}}$ we choose representatives $\bx$ of $x$ and $\by$ of $y$ in the space of geodesic currents. This is because both intersection form and entropy of geodesic currents are not projectively invariant.


We see similar definitions later. For example in \cref{sec:compactification}, we define the function $\mathcal{L}_{x}$ for any $x \in \mathbb{P}\mathcal{C}$. But the definition of $\mathcal{L}_{x}$ uses function $Q$ that is defined on the space of currents. 

To be consistent and precise we will use bold letters $\textbf{x}$, $\textbf{y}$, $\textbf{z}$ for currents and $x$, $y$, $z$ for elements in the space of projective currents. In particular, for any $x\in \mathbb{P}\mathcal{C}$ we will use $\textbf{x}$ for a representative of $x$ in the space of currents.

We occasionally write $[\bx]$ to denote the projective class of a geodesic current $\bx$, although this notation is less common. This occurs, for instance in \cref{decompsotion}.

In equations where suprema over closed curves/currents are taken, we will  write $\scc$ for simple closed curves or $\cc$ for closed curves.

\begin{remark} 
\label{rmk:boundaries_punctures} We remark that the hypothesis that $S$ is closed plays an important role throughout the paper. If $S$ had boundary components, it is known that there is no analog of the Liouville current $L_x$ associated to $x \in \Teich(S)$ ~\cite[Proposition~2.2]{T22:ThurstonCompactification}. Moreover, Equation~\ref{eq:dTh} does not define a metric anymore~\cite[Remark~1.1]{LiuSuZhong2015}.
If $S$ has punctures, the Liouville current exists~\cite[Page 56, Remark]{ES22:GeodesicCount}, and the Thurston metric is well-defined~\cite[Section~2.1]{PapadopoulosSu2016}, but then one loses continuity of intersection numbers $i(\mu, \cdot) \colon \mathcal{C}(S)\to \mathbb{R}$ (by~\cite[Section 6]{Sas22:CuspedCurrents}), a fact that is used throughout the paper. Moreover, the results on the decomposition of currents discussed in Section~\ref{decompsotion} are less straightforward. For instance, there exist currents with positive systole whose support is a lamination~\cite[Example~1.1]{BIPP21}
\end{remark}

\subsection{Outline of the paper and key ideas.}
The paper is organized as follows:

\begin{enumerate}

  \item[\S2\ ] This section provides necessary background and is divided into three parts. 
  \begin{enumerate}[label=(\arabic*)]
      \item In the first part, we consider an arbitrary metric space $(X,d)$. We discuss the horofunction boundary $X(\infty)$, the detour cost, Busemann points and certain \hyperref[prop:key]{\emph{key properties}} which, when satisfied by $(X,d)$, ensure that the compactification is well behaved. For instance, if $(X,d)$ satisfies all the key properties, then the natural map $X \to C(X)$ is an embedding, and the symmetrized detour cost $\delta$ defines a (possibly infinite-valued) metric on the horofunction boundary of $X$.

      \item   The second part provides background on $\mathcal{C}$, $\PC$, $\PC_{\mathrm{fill}}$ and Thurston’s metric $d_{Th}$. We discuss the extended Thurston metric $d$ on  $\mathbb{P}\mathcal{C}_{\mathrm{fill}}$ in more detail.

      \item In the third part, we provide a decomposition of a non-filling geodesic current and its projectivization based on the decomposition of geodesic currents in \cite{BIPP21}. 
  \end{enumerate}
  
  \item[\S3\ ] We prove that $(\mathbb{P}\mathcal{C}_{\mathrm{fill}}, d)$ satisfies all the \hyperref[prop:key]{\emph{key properties}}.
  \item[\S4\ ] Using the key properties proved in \S3\, we prove \cref{horoboundary}, the horoboundary result.
  \item[\S5\ ] We prove all non-filling projective currents are Busemann points. For every non-filling projective current we construct an unparametrized linear geodesic in $\PC_{\mathrm{fill}}$ that converges to it.  As a consequence, the symmetrized detour cost $\delta$ serves as a (possibly infinite-valued) metric on projective  non-filling currents.

  \item[\S6\ ] We express the detour cost in terms of the intersection form, and—via the decomposition discussed in \S2\ .3 derive necessary conditions for two non-filling currents $\eta$ and $\xi$ to satisfy $\delta(\eta, \xi) < \infty$.

  \item[\S7\ ] We characterize minimal measured laminations using the detour cost.

  \item[\S8\ ] We prove the isometric rigidity of $(\mathbb{P}\mathcal{C}_{\mathrm{fill}}, d)$. For two distinct surfaces $S_{g}$ and $S_{g'}$, an isometry $f: (\PC_{\mathrm{filll}}(S_{g}),d) \rightarrow (\PC_{\mathrm{filll}}(S_{g'}),d)$ extends to a homeomorphism $f^{*}:\mathbb{P}\mathcal{C}(S_{g}) \rightarrow \mathbb{P}\mathcal{C}(S_{g'})$ on their horofunction compactifications. The map $f^{*}$ takes boundary to boundary and preserves the symmetrized detour cost. Using the characterization from \S6 and \S7 we show that $f^{*}$ maps simple closed curves on $S_{g}$ to minimal measured lamination of $S_{g'}$ and preserves disjointness. We finally deduce rigidity observing that this is only possible when $g = g'$.

  \item[\S9\ ] We discuss examples of spaces of geometric structures on which the extended Thurston metric on geodesic currents $d$ restricts to give an asymmetric metric and discuss further examples. In particular, we give the necessary background on Anosov representations and discuss how they induce proper hyperbolic potentials in the sense of~\cite{CRS24:Joint}, which are defined in the Appendix.

  \item[\S10\ ] We discuss some open questions in relation to the asymetric metric $d$ on $\PC_{\mathrm{fill}}$ and its symmetric counterpart.
  
      \item[Appendix~\ref{sec:asymmetric}] We collect facts about the space of metric structures $\mathcal{D}_\Gamma$ of a non-elementary hyperbolic group $\Gamma$ and the more general space of metric potentials $\mathcal{H}^{++}_\Gamma$. We then discuss how to define an extension of $d$ on $\mathcal{H}^{++}_\Gamma$ and how this extension can be used to recover the asymmetric metric $d$ on $\PC_{fill}$.
\end{enumerate}

\subsection{Acknowledgments} The first author is grateful to her advisor, Jenya Sapir, for insightful comments and mentorship throughout the development of this research. The first author also thanks Lorenzo Ruffoni and Jing Tao for helpful discussions. The second author thanks Eduardo Reyes for helpful conversations. The authors thank Giuseppe Martone for his careful reading of the early versions of the manuscript. The second author was supported by Marie Skłodowska-Curie Action CurrGeo grant (101154865). The results on isometric rigidity and horofunction compactification form part of the first author’s Ph.D. thesis.
\section{Background}\label{sec:background}
The background is organized into three parts. The first part introduces the concepts of horoboundary, detour cost, and Busemann points, defined for any metric space $(X,d)$, without assuming symmetry of the metric.

The second part focuses on the specific metric space relevant to our work. We review geodesic currents, projective currents, Thurston's metric on Teichmüller space, and its extension to the space of projective filling currents.

In third part we discuss the decomposition of a non-filling projective current. This decomposition is used extensively in \cref{sec:detourcost} and  \cref{sccml} characterizing when two projective non-filling currents are finite detour cost apart.

\subsection{Notions for metric spaces.} 

Let $(X,d)$ be a metric space where $d$ is not necessarily symmetric; that is, $d$ satisfies all the properties of a metric except that we might not have $d(x,y) = d(y,x)$ for all $x,y \in X$. Most of the background developed in this paper applies to asymmetric metrics. However, in certain cases, it is convenient to consider the symmetrized metric $d_{\mathrm{sym}}$ defined by $d_{\mathrm{sym}}(x,y) \coloneqq d(x,y) + d(y,x)$.

We refer the reader to \cite{CW2015} for a more detailed discussion of what follows.

\subsubsection{\textbf{Horoboundary}} To define the horofunction compactification of a metric space $X$, we consider $C(X)$, the space of continuous real-valued functions on $X$. Fix a base point $b \in X$, and define the map $\Psi: X \hookrightarrow C(X)$ by
\[
\Psi(z) = \Psi_z,
\]
where $\Psi_z: X \to \mathbb{R}$ is given by
\begin{equation}
\Psi_z(x) = d(x, z) - d(b, z) \quad \text{for all } x, z \in X.
\end{equation}

It is known that the map $\Psi$ is injective and continuous \cite{Nonpositive, CW2015}. The horofunction boundary $X(\infty)$ of $X$ consists of all elements in the closure of $\Psi(X)$ that are not in $\Psi(X)$ itself. 
\[
X(\infty) := \overline{\{ \Psi_z \mid z \in X \}} \setminus \{ \Psi_z \mid z \in X \}.
\]

The closure $\overline{\{ \Psi_z \mid z \in X \}}$ is called the horofunction compactification of $X$. Changing the basepoint $b$ to another $b'$ gives us a homeomorphic horofunction compactification and boundary. We can say more about the horofunction boundary when the metric $d$ satisfies certain nice properties: 

\begin{keyproperties}\label{prop:key}\textcolor{white}{.}
 \begin{itemize}
    \item \textbf{Property A.} The metric $d_{\mathrm{sym}}$ is proper.
    \item \textbf{Property B.} Between any pair of points in $X$, there exists a geodesic with respect to $d$.
    \item \textbf{Property C.} For any point $x$ and sequence $x_n$ in $X$, we have $d(x_n, x) \to 0$ if and only if $d(x, x_n) \to 0$.
    \item \textbf{Property D} For every sequence $x_n$ in X, if $d_{\mathrm{sym}}(b, x_n)$ converges to infinity then so does $d(b, x_n)$.
\end{itemize}     
\end{keyproperties}

For example, if $d$ satisfies properties (A), (B), and (C) listed below, then the map $\Psi: X \hookrightarrow C(X)$ is actually an embedding \cite{CW2015}. It is known that the asymmetric Thurston metric on Teichmüller space satisfies properties (A) through (D) \cite{GuPa}. In \cref{sec:topology} we show that the extended Thurston metric on projective filling currents satisfies properties (A) through (D) as well.

\subsubsection{\textbf{The induced homeomorphism on the horoboundary}}\label{f*} 
Let $(X,d)$ and $(X',d')$ be two metric spaces, and let $f:(X,d) \to (X',d')$ be an isometry. The map $f$ induces a map $f^*: \Psi(X) \to \Psi(X')$ defined by $f^*(\Psi_z) = \Psi_{f(z)}$. It follows from the continuity of $f$ and properties of $d$ that $f^*$ is continuous and extends continuously to the horoboundary. This extension, also denoted $f^{*}: X(\infty) \rightarrow X'(\infty)$ is a homeomorphism between the horoboundaries \cite[Proposition 2.4]{CW2015}. Similarly, the inverse isometry $f^{-1}$ induces $(f^{-1})^*: X'(\infty) \to X(\infty)$. From the definitions of $f^*$ and $(f^{-1})^*$, it follows that these maps are inverses of each other, satisfying $(f^{-1})^* \circ f^* = Id_X$ and $f^* \circ (f^{-1})^* = Id_{X'}$, where $Id_X$ and $Id_{X'}$ denote the identity maps on $X$ and $X'$ respectively.

\subsubsection{\textbf{Geodesics and Busemann Point}}
A path $\gamma:[0, \infty) \rightarrow X$ is an unparametrized geodesic if it satisfies the additivity condition 
\begin{equation}\label{geoeq}
   d(\gamma(0), \gamma(s)) + d(\gamma(s), \gamma(t)) = d(\gamma(0), \gamma(t)) \text{ for all } s,t \in \mathbb{R}_{+} \text{ with } s \leq t.
\end{equation}

We say a parametrized curve $\gamma: [0, \infty) \rightarrow X$ is an \textit{almost-geodesic} if for each $\epsilon >0$,
\[
|d(\gamma(0), \gamma(s)) + d(\gamma(s), \gamma(t)) - t| < \epsilon
\]
for all $s$ and $t$ large enough, with $s \leq t$. All geodesics are almost geodesics. We say a point $z$ on the boundary of $X$ is a \textit{Busemann point} if there exist an almost geodesic $\gamma(t)$ such that $\gamma(t)$ approaches $z$ as $t \to \infty$.

\subsubsection{\textbf{Detour cost}} Let $\xi$ and $\eta$ be two points in the horoboundary $X(\infty)$ of $(X,d)$. The detour cost between $\xi$ and $\eta$ is defined to be 
\begin{equation}\label{detourcosteq}
    H(\xi,\eta) = \displaystyle{\inf _{\gamma} \liminf_{t \to \infty}\Big(d(b, \gamma(t)) + \Psi_\eta(\gamma(t)) \Big)}
\end{equation}
where the infimum is taken over all paths $\gamma(t): \mathbb{R}_{+} \to X$ converging to $\xi$.

\begin{figure}[h]
    \centering
    \includegraphics[width=0.4\linewidth]{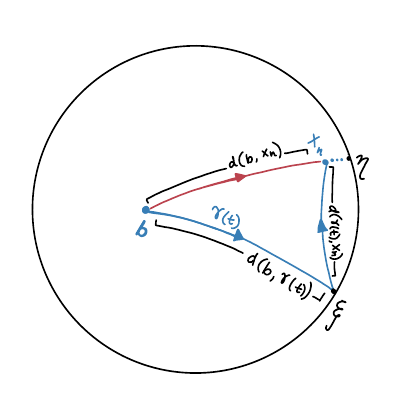}
    \caption{The disc represents the metric space $(X,d)$, with the boundary circle corresponding to the horoboundary $X(\infty)$. The detour cost $H(\xi, \eta)$ compares $d(b, \gamma(t)) + d(\gamma(t), x_n)$ with $d(b, x_n)$ as $n$ and $t$ tend to infinity. Here, $(x_n)$ is a sequence in $X$ converging to $\eta$. And $\gamma(t)$ is a path in $X$ converging to $\xi$ that realizes the infimum in the definition of $H(\xi, \eta)$.}
\label{detourcost}
\end{figure}

The detour cost is not necessarily symmetric; that is, for any two points $\xi$ and $\eta$ on the horoboundary, $H(\xi, \eta)$ is not necessarily equal to $H(\eta, \xi)$. Moreover, the detour cost between two points can be infinite. These differences aside, the detour cost behaves similarly to a metric. More precisely, if $\xi$, $\eta$, and $\upsilon$ are any three points in $X(\infty)$ it is known that:
\begin{enumerate}
    \item $H(\xi, \eta) \geq 0;$
    \item $H(\xi, \eta) \leq H(\xi, \upsilon) + H(\upsilon, \eta).$
\end{enumerate}
Also, if $(X,d)$ satisfies properties (A) through (D) then \(H(\xi, \xi) = 0\) if and only if $\xi$ is a Busemann point \cite[Proposition 5.4]{CW2015}. Therefore,  if $X$ satisfies properties (A) through (D) then the symmetrized detour cost $\delta$ defined below induces a (possibly infinite-valued) metric on the collection of Busemann points of $X$
\[
\delta(\eta, \xi) = H(\eta, \xi) + H(\xi, \eta).
\]
By \cite[Proposition 5.7]{CW2015}, the map $f^{*}: X(\infty) \rightarrow X'(\infty)$ discussed in \cref{f*} is an isometry with respect to $\delta$ in the sense that
\[
\delta(f^{*}(\eta), f^{*}(\xi)) = \delta(\eta, \xi).
\]

\subsection{Teichmüller space, geodesic currents and the extended Thurston metric}
\label{subsec:currents}
Geodesic currents are measures on the space of geodesics of a hyperbolic surface, furnishing a natural completion of the space of weighted closed curves.

In order to get a better behaved topological space, instead of the space of geodesics on $S$ one considers the space of geodesics of $\wt{S}$. Let $G(S)$ denote the space of bi-infinite unoriented unparametrized geodesics of $\wt{S}$. Since $\wt{S}$ is isometric to $\mathbb{H}^2$ (by choice of an arbitrary auxiliary hyperbolic metric) each such geodesic can be parameterized by a pair $(a,b) \in S^1 \times S^1 - \Delta$, where $\Delta=\{ (a,a) \in S^1 \times S^1 \}$. The coordinates $(a,b)$ and the flip $(b,a)$ induce the same unoriented geodesic. Hence, we have $G(S)=(S^1 \times S^1 - \Delta) /\mathbb{Z}_2$, where the quotient is by the action of the orientation flip. 

A \emph{geodesic current} on $S$ is defined as a $\pi_1(S)$-invariant locally finite, positive, Borel measure on $G(S)$. Let $\mathcal{C}(S)$ denote the space of geodesic currents of $S$. 

As a space of measures, $\mathcal{C}(S)$ is equipped with the \emph{weak$^*$}-topology, i.e., the weakest topology so that for every compactly supported continuous function $f \colon G(S) \to \mathbb{R}$, the map
\begin{align*}
F_f \colon \mathcal{C}(S)&\to \mathbb{R}\\
\mu &\mapsto \int f d\mu
\end{align*}
is continuous.
This topology is metrizable, second countable and sequential, and can be equivalently defined in terms of convergence of intersection numbers~\cite{Otal90:SpectreMarqueNegative,DLR2010}: $\mu_n \to \mu$ if and only if $\lim_n i(\mu_n, c)=i(\mu, c)$ for every closed curve $c$.
With respect to this topology, $\mathcal{C}(S)$ is locally compact and $\PC(S)$ is compact~\cite{Bon88}.

A natural example of geodesic current is given by every closed geodesic $c$ on $S$. Indeed, the full lift of $c$ to $\wt{S}$ is given by $\wt{C} \coloneqq \{ g \cdot \wt{c} \}_{g \in\pi_1(S)/\operatorname{stab}(\wt{c})}$, which is a subspace of $G(S)$. Let $\delta_c$ denote the Dirac measure supported on $\wt{C}$. Since $\wt{C}$ is a countable and discrete subspace of $G(S)$, $\delta_c$ is a locally finite Borel measure. Finally, it is $\pi_1(S)$-invariant as its support satisfies  $\pi_1(S)(\wt{C})=\wt{C}$. We will abuse notation and write $\bc$ instead of $\delta_c$ when referring to the geodesic current associated to the closed geodesic $c$. Note that, in a hyperbolic metric, every free homotopy class of closed curves contains a unique closed geodesic. Accordingly, $\bc \in \mathcal{C}(S)$ also denotes the geodesic current associated to the free homotopy class of any closed curve $c$.

The \emph{Teichm\"uller space} $\Teich(S)$ is the space of marked hyperbolic structures on $S$, i.e., a pair $(X,\phi)$ where $X$ is a hyperbolic metric on $S$ and $\phi \colon S \to X$ a homeomorphism, and two $(X,\phi), (X',\phi')$ are considered equivalent if there exists an isometry $I \colon X \to X'$ so that $I \circ \phi$ is homotopic to $\phi'$.

By Bonahon~\cite{Bon86,Bon88}, to every marked metric $\bx = (X, \phi) \in \Teich(S)$ one can associate a \emph{hyperbolic Liouville current} $L_{\bx}$, given as follows. Lift the marking $\phi$ to a homeomorphism of universal covers $\wt{\phi} \colon \wt{S} \to \wt{X}$, and consider the $X$-hyperbolic cross-ratio $[\cdot, \cdot, \cdot, \cdot]_X$ on the boundary of $\wt{X}$-- using the fact that $\wt{X}$ is isometric to $\mathbb{H}$--to asign to every subset of geodesics $[a,b] \times [c,d] \subset G(S)$, where $a,b,c,d$ are counter-clockwise oriented, the number
\[
L_{\bx}([a,b] \times [c,d]) = \log \left[\wt{\phi}(a),\wt{\phi}(b),\wt{\phi}(c),\wt{\phi}(d)\right]_X.
\]
This assignment is finitely additive and $\pi_1(S)$-invariant by the properties of the cross-ratio. Since the cross-ratio is always positive, the finite additivity can be promoted to countable additivity, furnishing a geodesic current.  In fact, this is an example of a \emph{filling geodesic current}, i.e., a geodesic current whose support intersects every geodesic on the surface. 
We use $\mathcal{C}_{\mathrm{fill}}$ to denote the space of \emph{filling geodesic currents}.

By Bonahon~\cite[Proposition~4.5]{Bon86}, the standard geometric intersection number between closed curves extends as a bilinear continuous function  $i(\cdot, \cdot)$ to the space of currents.

By~\cite[Theorem~1.5]{BIPP21}, $\bmu \in \mathcal{C}(S)$ is filling iff $i(\bmu, \bnu)>0$ for every non-zero $\bnu \in \mathcal{C}(S)$ iff the \emph{systole} of $\mu$ is positive, where systole of a geodesic current is defined as follows:
\[ \operatorname{Syst}(\mu) \coloneqq
 \inf \{ i(\mu, \mathbf{c}) \mid c \text{  is a closed geodesic on }S \}.
 \]

The map $\bx \mapsto L_{\bx}$ gives a proper embedding of $\Teich(S)$ into the space of projective filling currents $\PC_\mathrm{fill}(S)$. The space of projective filling currents is a quotient of $\mathcal{C}_{\mathrm{fill}}(S)$ obtained by identifying a filling geodesic current with its scalar multiples. Furthermore, the hyperbolic length of a curve $c$ with respect to $\bx \in \Teich(S)$, denoted by $\ell_{\bx}(c)$, can be recovered via
\begin{equation}
i(L_{\bx}, \bc) = \ell_{\bx}(c)
\label{eq:length_inter}
\end{equation} by~\cite{Bon86, Bon88}. We will eventually abuse notation and refer to $\bx$ instead of $L_{\bx}$ when referring to the geodesic current of a point $\bx \in \Teich(S)$.

Given $\bx \in \PC_\mathrm{fill}$, we can define the \emph{entropy of $\bx$} as
\begin{equation}
h(\bx) \coloneqq \limsup_n \frac{1}{n}\log \# \{c \mbox{ closed geodesics } : i(\bx, \bc) < n \}.
\label{eq:entropy}
\end{equation}

For example, when $\bx \in \Teich(S)$,
\begin{equation}\label{eq:entropy1}
\begin{aligned}
h(L_{\bx}) 
 &= \limsup_n \frac{1}{n}\log \#\{ c \text{ closed geodesics } : i(L_{\bx},\bc) < n \} \\
 &= \limsup_n \frac{1}{n}\log \#\{ c \text{ closed geodesics } : \ell_{\bx}(c) < n \} \\
 &= 1 .
\end{aligned}
\end{equation}
where the first equality is the definition of entropy, the second uses the Equation~\ref{eq:length_inter}, and the last is  by the classical curve counting result of Huber~\cite{entropy1}.

The following properties of entropy are easy to verify.

\begin{lemma}
The entropy satisfies:
\begin{itemize}
\item For every $a >0$ and every $\bx \in \mathcal{C}_\mathrm{fill}$, we have \[h(a \cdot \bx) = \frac{1}{a} \cdot h(\bx).\]
\item For every pair of filling currents $\bx, \by$ satisfying whenever $i(\by, \bc) \leq i(\bx, \bc)$ for every closed curve $\bc$, we have
\[
h(\bx) \leq h(\by).\]
\end{itemize}
\label{lem:entropy_properties}
\end{lemma}

A general fact that we will use repeatedly throughout the paper is that given a pair of filling geodesic currents $\bx, \by$, there exist positive constants $c_1(\bx, \by), c_2(\bx, \by)>0$ so that
\begin{equation}
c_1 < \frac{i(\by, \bc)}{i(\bx, \bc)} < c_2
\end{equation}
for every geodesic current $\bc$.
This follows from the fact that $\frac{i(\by, \cdot)}{i(\bx, \cdot)}$ is a well-defined continuous and positive function on the compact space of projective geodesic currents $\PC$.

From this observation, the following result follows readily.

\begin{lemma}
Given a pair of filling geodesic currents $\bx, \by$, there exist positive constants $c_1, c_2$ so that
\begin{equation}
c_1 < \frac{h(\by)}{h(\bx)} < c_2.
\label{eq:bilip_entropy}
\end{equation}
In particular, $h(\bx)<\infty$ for every $\bx \in \mathcal{C}_\mathrm{fill}$.
\label{lem:finite_entropy}
\end{lemma}
\begin{proof}
The first claim follows since by the paragraph above, we have constants $c_1>0, c_2>0$ so that
\begin{equation}
c_2 i(\bx, \bc) < i(\by, \bc) < c_1 i(\bx, \bc).
\end{equation}
for every closed curve $\bc$.
By Lemma~\ref{lem:entropy_properties}, Equation~\ref{eq:bilip_entropy} follows. From this, by Equation~\ref{eq:entropy1} it follows that $h(\bx)<\infty$ for every $\bx \in \mathcal{C}_\mathrm{fill}$ since $h(L_{\bx})=1$.
\end{proof}

The above discussion motivates the following definition. For $x,y \in \mathbb{P}\mathcal{C}_{\mathrm{fill}}(S)$, choose representatives $\bx,\by$ in the space of filling currents and set
\begin{equation}
    \operatorname{Dil}(x,y) \coloneqq \sup_{\bc \in \cc} 
    \frac{i(\by,\bc)}{i(\bx,\bc)},
\label{eq:dilation}
\end{equation}
which we call the \emph{dilation} of $x$ and $y$.

We can now write the extended asymmetric metric $d$ as
\begin{definition}
For $x$, $y$ in $\PC_{\mathrm{fill}}$, choose representative $\bx$ of $x$ and $\by$ of $y$ in $\mathcal{C}(S)$. The \emph{extended Thurston metric} between $x$ and $y$ is defined as follows:
\[
d(x, y) \coloneqq \log \left(\operatorname{Dil}(x,y) \cdot \frac{h(\by)}{h(\bx)} \right).
\]
\end{definition}

In Proposition~\ref{prop:asymmetric_geodesics} we establish the following result.

\begin{proposition}
The extended Thurston metric $d$ defines an asymmetric metric on $\PC_\mathrm{fill}$, this is for any $x, y, z \in \PC_{\mathrm{fill}}$.  it satisfies:
\begin{itemize}
\item $d(x, y) =0$ if and only if $x = y$.
\item $d(x, z) \leq d(x, y) + d(y, z)$.
\end{itemize}
\end{proposition}

\begin{proposition}
The metric $d$ extends the classical Thurston metric $d_{\mathrm{Th}}$ on $\Teich(S)$, i.e.,
\[
d|_{\Teich(S)} = d_{\mathrm{Th}}.
\]
\label{prop:d_extends_thurston}
\end{proposition}
\begin{proof}
By Equation~\ref{eq:entropy1} it follows that for $\bx, \by \in \Teich(S)$, $h(\bx)=h(\by)=1$.
From Sapir's proof of \cite[Lemma 6.1]{Jenya}, 
it follows that taking the supremum over all closed curves as opposed to only simple closed curves gives the same result in the definition of $d_{\mathrm{Th}}$. Hence, the result follows.
\end{proof}

\subsection{Decomposing a non-filling geodesic current}\label{decompsotion}

In this section we show that any non-filling geodesic current $\boldeta$ can be decomposed into a \emph{filling component} $\boldeta_{S'}$ and a \emph{lamination component} $\blambda_{\boldeta}$,
\[
\boldeta = \boldeta_{S'} + \blambda_{\boldeta}.
\]
This decomposition follows from the decomposition of geodesic current in \cite[Theorem 1.2]{BIPP21}. We will go through some of the details of this construction in order to specify what we classify as $\boldeta_{S}$ and what as $\blambda_{\boldeta}$. Even though the decomposition is defined at the level of geodesic currents, we will see that certain features of the construction remain invariant within a projective class.

For any geodesic current $\boldsymbol{\eta}$, \cite{BIPP21} considers the set $\mathcal{E}_{\boldeta}$ defined below. For a closed curve $c$ on the surface, $\bc$ represents the geodesic current associated with $c$.
\begin{equation}\label{Eset}
\mathcal{E}_{\boldeta} := 
\begin{cases} 
      c \text{ is a closed curve in }S :& i(\boldsymbol{\eta}, \bc) = 0  \text{ and } \\
      
      & i(\boldsymbol{\eta}, \bc') \neq 0 \text{ for every closed curve } \bc' \text{ with } i(\bc, \bc') \neq 0.
\end{cases}
\end{equation}
Note that, $\mathcal{E}_{\boldeta}$ does not vary within the projective class of $\boldeta$.
\[
\mathcal{E}_{\eta} = \mathcal{E}_{\boldeta}, \text{ where } \boldeta \text{ is a representative of } \eta \text{ in the space of currents.}
\]
The set $\mathcal{E}_{\boldeta}$ consists of pairwise disjoint simple closed curves. These simple curves decompose $S$ into a finite number of subsurfaces $\{S_{1}, S_{2}, \dots S_{n}\}$.  We say a geodesic current $\bmu$ \emph{is supported on} $S_{i}$ if the support of $\bmu$ contains geodesic in $\mathbb{H}^2$ that projects to the interior of $S_{i}$. For a geodesic current, \cite{BIPP21} defines systole of $\boldeta$ on $S_{i}$ as follows:

\begin{definition}\label{def:systole} Let $\boldeta$ be a geodesic current and let $S'$ be a subsurface of $S$. The \emph{systole} of $\boldeta$ on $S'$ is defined as 
\[
\operatorname{Syst}_{S'}(\boldsymbol{\eta})
\coloneqq
\inf \left\{
  i(\boldsymbol{\eta}, c)
  \mid
  c \text{ is a closed geodesic in the interior of } S'
\right\}.
\]

We say $\boldsymbol{\eta}$ fills $S'$, if $\Syst_{S'}(\boldsymbol{\eta}) > 0$.
\end{definition}
From in Theorem 1.2 in \cite{BIPP21} we know that a geodesic current $\boldsymbol{\eta}$ can be decomposed into \[
\boldsymbol{\eta} = \displaystyle{ \sum_{i =1}^{n} \boldsymbol{\eta}_{i} + \sum_{c \in \mathcal{E}_{\boldeta}} \kappa_{c}\bc},
\]
where, $\bc$ is the geodesic current corresponding to the simple closed curve $c$, $\kappa_{c}$'s are non negative real numbers and $\boldsymbol{\eta}_{i}$ is a geodesic current supported only on $S_{i}$. The current $\boldsymbol{\eta}_{i}$ fall into exactly on of the following three cases:
\begin{enumerate}
    \item $\boldsymbol{\eta}$ is not supported on $S_{i}$. In this case $\boldsymbol{\eta}_i = 0;$
    \item $\Syst_{S_i}(\boldsymbol{\eta}_{i}) > 0$, that is $\boldsymbol{\eta}_{i}$ fills $S_{i};$
    \item $\Syst_{S_i}(\boldsymbol{\eta}_{i}) = 0$. In this case $\eta_{i}$ is a measured lamination supported only on $S_{i}$ intersecting every curve in the interior $S_{i}$. 
\end{enumerate}

In our decomposition

\begin{equation}\label{decompostioneq}
    \boldsymbol{\eta} = \boldsymbol{\eta}_{S'} + \boldsymbol{\lambda_{\eta}}
\end{equation}
\begin{itemize}
    \item $\boldsymbol{\eta_{S'}}$ consists of currents $\boldsymbol{\eta_{i}}$ that satisfies $\Syst_{S_i}(\boldsymbol{\eta}_{i}) > 0$. For any such $\boldeta_{i}$, if the simple closed curves on the boundary of $S_{i}$ supports $\boldeta$, we take into account those closed curve as well. 
\[
\boldsymbol{\eta}_{S'} = \displaystyle{\sum\limits_{ \Syst_{S_i}(\eta_i) > 0} \boldsymbol{\eta}_{i} + \sum\limits_{\substack{\Syst_{S_i}(\eta_i) > 0\\ c \in \partial S_{i}}} \kappa_{c}c}.
\]

Here $S'$ is the union of all the $S_{i}$'s for which $\Syst_{S_i}(\boldsymbol{\eta}_{i}) > 0$. The subsurface $S'$ does not have to be connected.\\

\item The \textit{measured lamination component} of $\boldsymbol{\eta}$ is $\boldsymbol{\lambda_{\eta}}= \boldsymbol{\eta} - \boldsymbol{\eta}_{S'}$. Because $\boldsymbol{\lambda}_{\boldsymbol{\eta}}$ consists of only measured laminations and disjoint simple closed curves, it is a measured lamination.  
\end{itemize}

\begin{figure}[ht]
    \centering
    \includegraphics[width=\linewidth]{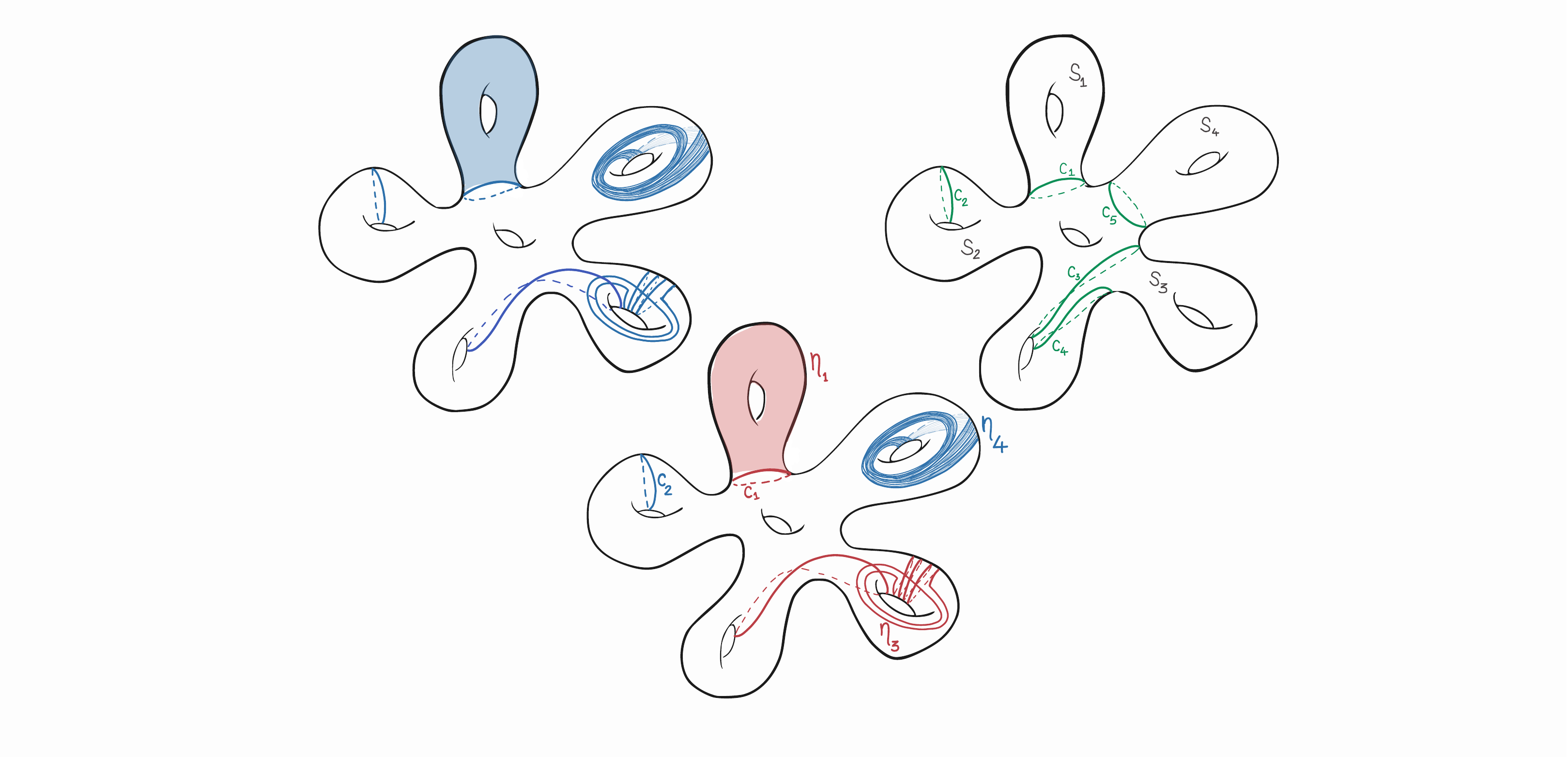}
    \caption{
    Top left: The blue curves represent geodesics in the support of the geodesic current $\boldeta$. If $\boldeta$ is supported on all geodesics in the interior of a subsurface, we indicate this by shading that subsurface. 
    Top right: The green simple closed curves belong to $\mathcal{E}_{\boldeta}$ and decompose the surface into subsurfaces $\{S_1, \dots, S_4\}$. 
    Bottom: The red components contribute to the filling component $\boldeta_{S'}$, while the blue component represents the lamination part $\blambda_{\boldeta}$.
    }
    \label{fig:enter-label}
\end{figure}

Any two geodesic currents $\boldsymbol{\eta}$ and $k \boldsymbol{\eta}$ in the same projective class fill the same subsurface of $S$, and the decomposition of $k\boldsymbol{\eta}$ is $k$ times the right hand side of \cref{decompostioneq}. This allows us to write the decomposition of the projective current $\eta = [\boldeta]$ as $[\boldeta_{S'} + \blambda_{\boldeta}]$

\begin{theorem}\label{thmdecomposition}
    Any non-filling geodesic current $\boldeta \in \mathcal{C}$ can be decomposed into
    \[
    \boldeta = \boldeta_{S'} + \blambda_{\boldeta},
    \]
    where $\boldeta_{S'}$ is a geodesic current that fills a subsurface $S'$ of $S$ and $\blambda_{\boldeta}$ is a measured lamination. We call $\boldeta_{S'}$ the \emph{filling component} of $\boldeta$, $\blambda_{\boldeta}$ the \emph{measured lamination component} of $\boldeta$, and $S'$ the \emph{subsurface filled} by $\boldeta$.
\end{theorem}

.

\section{Topological properties of The extended Thurston metric and the key properties}\label{sec:topology}

In this section, we show that the extended Thurston metric on projective filling currents satisfies the \hyperref[prop:key]{Key Properties}. Jenya Sapir in \cite[Corollary 7.2]{Jenya} establishes that $d_{\mathrm{sym}}$ is proper. Thus, $d$ satisfies Property $A$. It is left to show that it satisfies properties $B$ - $D$. We establish Property C after Property D, as the proof of Property C depends on Property D.

\subsection{Property B: Geodesic metric space} We establish that Property B holds for \( (\mathbb{P}\mathcal{C}_{fill}, d) \), which follows essentially from~\cite[Theorem~1.2]{StephenEduardo}. For the precise definition of $\mathcal{D}_\Gamma$ see Appendix~\ref{sec:asymmetric}.

\begin{proposition}\label{geodesic}[Property B]
    The space of projective filling currents equipped with the extended Thurston metric is geodesic.
\end{proposition}

\begin{proof}
By Proposition~\ref{prop:geod_sym} (Cantrell--Reyes), for any $x, y \in \mathbb{P}\mathcal{C}_{fill}(S)$, a geodesic for $d_{\mathrm{sym}}$ from $x$ to $y$ in the space of metric structures $\mathcal{D}_
{\pi_1(S)}$ is given by  
\[
r(t) = [t \by + \theta_{\bx}^{\by}(t) \bx], \quad 0 \leq t \leq h(\by).
\]
Here, $\theta_{\bx}^{\by}$ is the \textit{Manhattan curve} associated to $x$ and $y$ (see Section~\ref{subsec:manhattan}), and $[\cdot]$ denotes the rough similarity class in $\mathcal{D}_\Gamma$.

Note that $\theta_{\bx}^{\by}(t) \colon \mathbb{R} \to \mathbb{R}$ is a continuously differentiable function 
and it satisfies
\[
\theta_{\bx}^{\by} (0) = h(\bx) \quad \text{and} \quad \theta_{\bx}^{\by}(h(\by)) = 0,
\]
while remaining positive on the interval $[0, h(\by)]$. 

Since positive linear combinations of filling currents are filling, the geodesic $r(t)$ remains within the subspace $\mathbb{P}\mathcal{C}_{fill}\subset \mathcal{D}(\Gamma)$. This proves that $(\mathbb{P}\mathcal{C}_{fill}, d_{\mathrm{sym}})$ is a geodesic metric space.
Finally, by Proposition~\ref{prop:dsym_implies_d_geodesic}, any geodesic wrt $d_{\mathrm{sym}}$ is a geodesic for $d$, so $(\mathbb{P}\mathcal{C}_{fill}, d)$ is geodesic.
\end{proof}

\subsection{Property D: Behavior at infinity} We show that a sequence $\{x_{n}\} \subset \mathbb{P}\mathcal{C}_{fill}$ escapes every compact subset of $\mathbb{P}\mathcal{C}_{fill}$ if and only if if $d(b, x_n) \to \infty$ for any choice of $b \in \mathbb{P}\mathcal{C}_{fill}$. This is proven in the case of $\mathcal{T}$ by \cite{GuPa}. They also show that for a sequence ${x_n}$ in Teichmüller space and some basepoint $b \in \mathcal{T}$, we have $d_{Th}(b, x_n) \to \infty$ if and only if $d_{Th}(x_n, b) \to \infty$. 
We comment on why this equivalence is harder to prove for the space of currents.

\begin{definition}
   We say a  sequence $\{x_n\} \subset \mathbb{P}\mathcal{C}_{fill}$ tends to infinity, and we write $x_n \to \infty$, if no compact subset $K$ contained in $ \mathbb{P}\mathcal{C}_{fill}$ contains infinitely many elements in $\{x_{n}\}$.
\end{definition}


\begin{lemma}\label{unboundedd}
Let $b \in \mathbb{P}\mathcal{C}_{fill}$, and suppose that $\{x_n\}$ is a sequence of filling currents converging to some $\eta \in \mathbb{P}\mathcal{C}$. Then the distance $d(b, x_n)$ tends to infinity if and only if $\eta$ is non-filling.
\end{lemma}

\begin{proof}
If $\{x_{n}\}$ converges to a filling current $\eta$, then by continuity of the metric $d$, it follows that $d(b, x_{n})$ converges to $d(b, \eta)$, which is finite. Hence, $d(b, x_{n})$ does not go to infinity. It remains to show that if $\{x_n\}$ converges to a non-filling current $\eta$, then $d(b, x_n)$ tends to infinity. We will establish the result for $b \in \mathcal{T}$; the general case for $b \in \mathbb{P}\mathcal{C}_{fill}$ will then follow by the triangle inequality.

For each $x_{n}$ in our sequence, we pick a representative $\bx_{n}$ in the space of currents such that $\{\bx_{n}\}$ converges to a representative $\boldeta$ of $\eta$. Note that, for any $b \in \mathcal{T}$, its representative $\bb \in \mathcal{T} \subset \mathcal{C}$ satisfy $h(\bb) = 1$ by Equation~\ref{eq:entropy1}. This gives us
\[
d(b, x_{n}) =  \log \Bigg( \sup_{c \in \mathbb{P}\mathcal{C}} \frac{i(\bx_{n}, \bc)}{i(\bb, \bc)}h(\bx_{n}) \Bigg).
\]  
It follows Proposition~\ref{prop:entropy_proper} that if a  sequence $\{\bx_{n}\}$ in currents converges to a non-filling current, then $\displaystyle{\lim_{n \to \infty}h(\bx_{n}) = \infty}$. Hence, to show $\displaystyle{\lim_{n \to \infty} d(b, x_{n})}$ it suffices to prove that 
\[
\displaystyle{ \lim_{n \to \infty}\sup_{c \in \mathbb{P}\mathcal{C}} \frac{i(\bx_{n}, \bc)}{i(\bb, \bc)}}
\]
is bounded away from zero. Note that,
\[
  \frac{i(\bx_{n}, \bc)}{i(\bb, \bc)} \leq \sup_{c \in \mathbb{P}\mathcal{C}} \frac{i(\bx_{n}, \bc)}{i(\bb, \bc)} \text{ for all } n \in \mathbb{N} \text{ and for all } \bc \in \mathcal{C}
\]
Therefore, 
\begin{align*}
   \frac{i(\boldeta, \bc)}{i(\bb, \bc)} &\leq \lim_{n \to \infty}\sup_{c \in \mathbb{P}\mathcal{C}} \frac{i(\bx_{n}, \bc)}{i(\bb, \bc)} \text{ for all } \bc \in \mathcal{C}\\
\end{align*}
We can always find some geodesic current $\blambda$ such that $i(\boldeta, \blambda) = k_{1} > 0$. Since $\bb$ is filling, we get $i(\bb,\blambda) = k_{2} > 0$. Hence
\[
0 < \frac{k_1}{k_{2}} < \lim_{n \to \infty}\sup_{c \in \mathbb{P}\mathcal{C}} \frac{i(\bx_{n}, \bc)}{i(\bb, \bc)} 
\]
This proves the desired result. Consequently, we get 
\[
\displaystyle{\lim_{n \to \infty}\sup_{c \in \mathbb{P}\mathcal{C}} \frac{i(\bx_{n}, \bc)}{i(\bb, \bc)}h(\bx_{n}) = \infty}
\]
and hence $\displaystyle{\lim_{n \to \infty} d(b, x_{n}) = \infty}$.
\end{proof}

\begin{proposition}\label{inftythm}
    Let $\{x_n\}$ be a sequence in $\mathbb{P}\mathcal{C}_{fill}$, the following three properties are equivalent:\begin{enumerate}
        \item $x_n \to \infty$.
        \item Any convergent subsequences of $\{x_n\}$ has to converge to a non-filling projective current.
        \item For all $b \in \mathbb{P}\mathcal{C}_{fill}$, $d(b, x_n) \to \infty$.
    \end{enumerate}
\end{proposition}

\begin{proof} 
$(1) \implies (2)$: Since $\mathbb{P}\mathcal{C}$ is compact, the sequence $\{x_{n}\}$ has a convergent subsequence $\{x_{n'} \}$. Assume $x_{n'}$ converges to a projective current $\mu$. We claim that $\mu$ cannot be filling.

Assume, for contradiction, that $\mu$ is filling. Then we can find a real number $N$ and a compact neighborhood $K_{\mu}$ of $\mu$ contained in $\mathbb{P}\mathcal{C}_{fill}$ such that $x_{n}' \in K_{\mu}$ for all $n > N$. Since $K_{\mu}$ contains infinitely many elements of $\{x_{n}'\} \subset \{x_{n}\}$, the sequence $\{x_{n}\}$ cannot escape to infinity.

$(2) \implies (1)$ If a compact set $K \subset \mathbb{P}\mathcal{C}_{fill}$ contains infinitely many elements of the sequence $\{x_n\}$, then $K \cap \{x_n\}$ admits a convergent subsequence $\{x_n'\}$ whose limit is a projective filling current contained in $K$, contradicting (2).

$(2) \iff (3)$ Assume there exists some $b \in \mathbb{P}\mathcal{C}_{fill}$ such that $d(b,x_n)$ does not diverge to infinity. Then there is a subsequence $\{x_n'\} \subset \{x_n\}$ for which $d(b,x_n')$ is bounded by some constant $L$. Since $\mathbb{P}\mathcal{C}$ is compact, $\{x_n'\}$ admits a convergent subsequence, which we continue to denote by $\{x_n'\}$. Let $\mu = \lim_{n \to \infty} x_n'$. Then
\[
d(b,\mu) = \lim_{n \to \infty} d(b,x_n') \leq L.
\]
By \cref{unboundedd}, $\mu$ must be filling, contradicting (2). The converse follows by reversing the argument.
\end{proof}

\begin{corollary}[Property D] Let $b \in \mathbb{P}\mathcal{C}_{fill}$ then for every sequence $x_n$ in $\mathbb{P}\mathcal{C}_{fill}$, if $d_{\mathrm{sym}}(b, x_n)$ converges to infinity, then so does $d(b, x_n)$.
\end{corollary}

\begin{proof} Let us assume $d_{\mathrm{sym}}(b, x_n)$ converges to infinity. By \cite[Corollary 7.2]{Jenya} infinitely many elements in $\{x_{n}\}$ cannot be contained in a compact set in $\mathbb{P}\mathcal{C}_{fill}$. Therefore, $x_n \to \infty$. But then by \cref{inftythm} we have $d(b, x_{n}) \to \infty$.
\end{proof}

\begin{remark}\label{remark:boundeddistance}
It is conceivable that $d(x_n, b)$ remains bounded along certain sequences  $\mathbb{P}\mathcal{C}_{fill}$ that escape to infinity. Indeed, choose a basepoint $b \in \mathcal{T}(S)$. Let $\bb \in \mathcal{T} \subset \mathcal{C}$ denote its representative. By \cref{inftythm}, we can find a subsequence of $\{x_n\}$ that converges to a projective non-filling current $\eta$. We continue to denote this subsequence by $\{x_n\}$. Let $\bx_n$ be a representative of $x_n$ such that the sequence $\{\bx_n\}$ converges to a representative $\boldeta$ of $\eta$. Then
\[
\lim_{n \to \infty} d(x_n, b) = \lim_{n \to \infty} \log\left( \sup_{c \in \mathbb{P}\mathcal{C}} \frac{i(\bb, \bc)}{i(\bx_n, \bc)} \frac{1}{h(\bx_n)} \right).
\]
For a sequence $\{\bx_n\}$ converging to a non-filling current $\eta$, the entropy $h(\bx_n)$ goes to infinity. Therefore, by Proposition~\ref{prop:entropy_proper}, we have
\[
\lim_{n \to \infty} \frac{1}{h(\bx_n)} = 0.
\]
On the other hand, since $\boldeta$ is non-filling, there exists a sequence of closed curves $\blambda_n$ such that $i(\boldeta, \blambda_n) \to 0$ \cite[Proposition~5.1]{BIPP21}. Consequently, $\lim_{n \to \infty} i(\bx_n, \blambda_n) = 0$. Since the numerator in the ratio considered below is bounded away from zero for all $\blambda_n$s, it follows that
\[
\lim_{n \to \infty} \sup_{c \in \mathbb{P}\mathcal{C}} \frac{i(\bb, \bc)}{i(\bx_n, \bc)} = \infty.
\]
This means that $\lim_{n \to \infty} d(x_n, b)$ compares the rate at which $h(\bx_n)$ is increasing and the rate at which the supremum value of the ratio is increasing. We thus expect this limit to be finite for some escaping sequences.
\end{remark}

\subsection{Property C: Distance approaching zero} Let $\{x_n\}$ be a sequence in $\mathbb{P}\mathcal{C}_{fill}$ and let $x \in \mathbb{P}\mathcal{C}_{fill}$. In this section, we show that $d(x_n, x) \to 0$ if and only if $d(x, x_n) \to 0$. We also show that this is equivalent to the sequence $\{x_n\}$ converging to $x$.

We will first prove the following lemma.

\begin{lemma} \label{zerolemma}
Let $b \in \mathbb{P}\mathcal{C}_{fill}$ and let $\{x_{n}\}$ be a sequence in $\mathbb{P}\mathcal{C}_{fill}$ such that $d(x_n, b) \to 0$, then any limit point $y$ of $x_{n}$ in $\mathbb{P}\mathcal{C}$ is a filling projective current.
\end{lemma}

\begin{proof}
Suppose $\{x_{n}\}$ has a subsequence $\{ x_{n}' \}$ that converges to a non-filling projective current $\eta$. By \cref{inftythm}, $x_{n}' \to \infty$. We will prove that for any $b \in \mathbb{P}\mathcal{C}_{fill}$, $d(x_{n}',b)$ cannot converge to 0, and consequently
\[
\lim_{n \to \infty} d(x_{n}, b) \neq 0.
\]
Let $K$ be a compact subset of $\mathbb{P}\mathcal{C}_{fill}$ such that $b$ in contained in the interior of $K$. Since $x_{n}' \to \infty$, we can find an $N$ such that $x_{n}' \notin K $ for all $n > N$. Since we choose $K$ to be a proper subset of $\mathbb{P}\mathcal{C}_{fill}$ with non-empty interior, the boundary of $\partial K$ is non-empty.
We define
\[
D = \displaystyle{\inf_{ y \in \partial K} d(y,b)}.
\]
As $K$ is compact and $\mathbb{P}\mathcal{C}$ is Hausdorff, $\partial K$ is closed and therefore compact. Hence, the infimum of $d(y,b)$ is attained at some $z \in \partial K$. Note, that $z \neq b$, since $b$ is in the interior of $K$. Hence, we have 
\[
D = d(z, b) \neq 0
\]
But now, for all ${x_n'}$ not in $K$, we have $d(x_n', b) > D$. This is because $(\mathbb{P}\mathcal{C}_{fill}, d)$ is a geodesic metric space and for each $x_{n}'$ not in $K$, we can find a geodesic $\gamma_{n'}: [0,1] \rightarrow \mathbb{P}\mathbb{C}_{fill}$ such that $\gamma_{n'}(0) = x_{n}'$ and $\gamma_{n'}(1) = b$. Since $\gamma_{n'}$ is continuous, the preimages of the interiors of $K$ and of $\mathbb{P}\mathcal{C} \setminus K$ under $\gamma_{n'}$ are disjoint, non-empty open subsets of $[0,1]$. Hence, there exists some $t \in [0,1]$ such that $\gamma_{n'}(t)$ lies neither in the interior of $K$ nor in the interior of $\mathbb{P}\mathcal{C} \setminus K$.
This means $\gamma_{n'}(t) \in \partial K$. Therefore, we get
\begin{align*}
    d(x_{n}', b) &= d(x_{n'}, \gamma_{n'}(t)) + d(\gamma_{n'}(t), b)\\
                 &\geq D + d(\gamma_{n'}(t), b)
\end{align*}
Since $d(\gamma_{n'}(t), b)$ is non-zero, $d(x_{n}', d) > D$. But then
\[
\lim_{n \to \infty} d(x_{n}', b) \neq 0
\]
completing our proof.

\end{proof}

We now show that the entropy is a continuous function on filling currents.

\begin{lemma}
Let $(\bx_n)$ be a sequence of filling geodesic currents converging to $\bx$. Then $h(\bx_n) \to h(\bx)$.
\label{lem:entropy_continuous}
\end{lemma}
\begin{proof}
Suppose $\bx_n \to \bx$.
Since \[\Dil(\bx, \by) \coloneqq \sup_{\gamma \in \cc}\frac{i(\by, \gamma)}{i(\bx, \gamma)} = \max_{\nu \in \mathcal{P}\mathcal{C}} \frac{i(\by, \nu)}{i(\bx, \nu)},\]
the function $F(\cdot) =\frac{i(\by, \cdot)}{i(\bx, \cdot)}$ is a continuous function on projective geodesic currents, and the space of projective geodesic currents is compact, by Berge’s
 maximum theorem~\cite[Theorem~17.31]{AB06:InfiniteDimensional}, $\Dil(\bx, \cdot)$ is continuous. 
We have
\[
\Dil(\bx, \bx_n)^{-1} i(\bx_n, \gamma) \leq i(\bx, \gamma) \leq \Dil(\bx, \bx_n) i(\bx_n, \gamma) 
\]
for every $\gamma$.
So, for every $\epsilon>0$, there exists $n_0$ so that for every $n \geq n_0$,
we have
\[
1-\epsilon   \leq \frac{i(\bx, \gamma)}{i(\bx_n, \gamma)}\leq 1 + \epsilon. 
\]
Hence, if $G(L, \bx) \coloneqq \{ \gamma \in \cc : i(\bx, \gamma) < L \}$, we have
\[
G(L/(1+\epsilon), \bx_n) \subset G(L, \bx) \subset G(L/(1-\epsilon), \bx_n).
\]
Hence 
\[
\frac{1}{1+\epsilon} h(\mu_i) \leq h(\mu) \leq \frac{1}{1-\epsilon} h(\mu_i)
\]
thus, taking logarithms, for $\epsilon$ small enough, we have
\[
- \epsilon \leq \log( h(\bx_n) ) - \log( h(\bx) )\leq \epsilon.
\]
and hence $h(\bx_n) \to h(\bx)$.
\end{proof}

Next, we prove $h$ is a proper function on filling currents.

\begin{proposition}
Let $\bx_n$ be a sequence of filling geodesic currents converging to a geodesic current $\bx$. Then $\bx$ is non-filling if and only if $\lim h(\bx_n)\to +\infty$.
\label{prop:entropy_proper}
\end{proposition}
\begin{proof}
Suppose first $\bx$ is filling. Since the space of filling geodesic currents is locally compact, there exists a compact neighborhood $K$ of $\bx$.
Hence, there exists $n_0$ so that for $n \geq n_0$, $\bx_n \in K$, and since $h$ is continuous by Lemma~\ref{lem:entropy_continuous}, there exists a positive constant $C>0$ so that $h(\bx_n) \leq C$ for all $n \geq n_0$.
Suppose now $\bx$ is non-filling, and, for the sake of contradiction, suppose that there exist a constant $C>0$ so that for all $n$, $h(\bx_n)<C$.
Then, for each $n$, there let $h_n \coloneqq h(\bx_n)>0$.
We have that $h(h_n \bx_n)=1$. 
Let $\bx'$ be a limit point for a subsequence $h_{n_k} \bx_{n_k}$.
Then, by \cite[Lemma~3.7]{StephenEduardo}, it  follows that $\bx'$ is filling.
Since $h_{n_k}$ is bounded, up to further subsequence, $h_{n_k} \to h' \geq 0$.
Since $\bx_{n_k} \to \bx$, we have $\bx' = h' \bx$. If $h'=0$, then $\bx'$ is the zero current, which is a contradiction. If $h'>0$, then $\bx'$ is non-filling, which is a contradiction. This finishes the proof.
\end{proof}
\begin{proposition}[Property C]\label{propertyC}
     For any sequence $\{x_n\}$ in $\mathbb{P}\mathcal{C}_{fill}$ and for any $x \in \mathbb{P}\mathcal{C}_{fill}$ the following three properties are equivalent:
    \begin{enumerate}
        \item $x_n$ converges $x$ in $\mathbb{P}\mathcal{C}_{fill}$
        \item $d(x, x_{n}) \to 0$ 
        \item $d(x_n,x) \to 0$ 
    \end{enumerate}
\end{proposition}

\begin{proof}
$(1) \implies (2)$ Let $\bx_{n}$ be a representative of $x_n$ in the space of currents so that the sequence $\{ \bx_{n} \}$ converges to a representative $\bx \in \mathcal{C}$ of $x\in \mathbb{P}\mathcal{C}_{fill}$. Intersection form is continuous by Bonahon~\cite[Proposition~4.5]{Bon86}, and so is entropy by Lemma~\ref{lem:entropy_continuous}. Hence, we get: 
\[
\lim_{n \to \infty} i(\bx_{n}, \bc) = i(\bx, \bc)
\]
for all geodesic currents $\bc \in \mathcal{C}$ and
\[
 \lim_{n \to \infty} h(\bx_{n}) = h(\bx) 
\]
As a consequence 
\[
\displaystyle{\lim_{n \to \infty} d(x, x_n) = \lim_{n \to \infty} \log \Bigg( \sup_{c \in \mathbb{P}\mathcal{C}} \frac{i(\bx_{n},\bc)}{i(\bx, \bc)} \frac{h(\bx_n)}{h(\bx)} \Bigg) = 0}
\]
A similar argument shows $(1) \implies (3)$.

$(2) \implies (1)$ Papadopoulos and Th\'eret proved this statement for the Thurston metric on Teichmüller space in \cite[Theorem 2]{GuPa}. The proof here follows a similar strategy.

Assume, for contradiction, that $d(x, x_n) \to 0$ but the sequence $\{x_n\}$ does not converge to $x$. Then there exists a neighborhood $N_x$ of $x$ such that, for every $k \in \mathbb{N}$, there exists $n_k > k$ with $x_{n_k} \notin N_x$. We have that $d(x, x_{n_k}) \to 0$.
By compactness of $\mathbb{P}\mathcal{C}$, the sequence $\{x_{n_k}\}$ admits a convergent subsequence $\{x_{n_k'}\}$ with limit $y \in \mathbb{P}\mathcal{C}$. Since $\{x_{n_k'}\}$ is a subsequence of $\{x_{n_k}\}$ we still have $
d(x, x_{n_k'}) \to 0$.
Since, $d(x, x_{n_k'})$ does not go to infinity, by \cref{unboundedd}, the limit $y$ of the subsequence $\{x_{n_k'}\}$ must be a filling current. Applying the triangle inequality, we obtain
\[
d(x,y) \leq d(x, x_{n_k'}) + d(x_{n_k'}, y) \text{ for all } n_k'.
\]

We know that $d(x, x_{n_k'}) \to 0$. Since the subsequence $\{x_{n_k'}\}$ converges to $y$, it follows that $d(x_{n_k'}, y) \to 0$ as well. Hence, $d(x, y) = 0$ and  $x = y$.

This implies that $x$ is a limit point of the sequence $\{x_{n_k}\}$. In particular, every neighborhood of $x$, including $N_x$, must contain infinitely many elements of the subsequence $\{x_{n_k'}\}$. However, this is a contradiction by the construction of $\{x_{n_k'}\}$.

The proof of $(3)\implies(1)$ is similar, except that we use \cref{zerolemma} instead of \cref{unboundedd} to argue that the limit point $y$ of $\{x_{n_k'}\}$ is a filling current.
\end{proof}

\section[Horofunction compactification of PC\textsubscript{fill} is PC.]{Horofunction compactification of $\mathbb{P}\mathcal{C}_{fill}$ is $\mathbb{P}\mathcal{C}$.}
\label{sec:compactification}

In this section we prove that the horofunction compactification of $(\mathbb{P}\mathcal{C}_{fill}, d)$ is homeomorphic to the space of projective currents $\mathbb{P}\mathcal{C}$, and the horofunction boundary of $(\mathbb{P}\mathcal{C}_{fill},d)$ consists of projective non-filling currents. This result extends a theorem of Walsh, which states that the horofunction compactification of $(\mathcal{T}, d_{Th})$ is homeomorphic to Thurston's compactification of Teichmüller space \cite{CW2015}.

To establish the horofunction compactification of $\mathbb{P}\mathcal{C}_{fill}$ is $\mathbb{P}\mathcal{C}$, we define a map $\Psi: \mathbb{P}\mathcal{C} \to C(\mathbb{P}\mathcal{C}_{fill})$ such that for any $z$ in $\mathbb{P}\mathcal{C}_{fill}$, the image of $z$ is $\Psi_{z}: \mathbb{P}\mathcal{C}_{fill} \rightarrow \mathbb{R}$ is 
\[
\Psi_{z}(x) = d(x,z) - d(b,z) \text{ for all } x \in \mathbb{P}\mathcal{C}_{fill}.
\]
We assume that our basepoint $b \in \mathcal{T}$. The asymmetric metric $d$ is only defined on projective filling currents. Consequently, when $z$ is non-filling we cannot work with the above definition of $\Psi_{z}$. In what follows we will extend the definition of the map $\Psi$ to \emph{all} projective currents. In \cref{sec2.5} and \cref{sec2.6} we will prove that the map $\Psi$ thus obtained is an embedding and hence is homeomorphic to its image.

\subsection{Defining \texorpdfstring{$\Psi$}{Psi} on projective currents}
In order to define the map $\Psi_z$ for non filling currents, we will first have to define the following component maps. These component maps are the same as those defined by Walsh in \cite{CW2015}, but adapted to the setting of $(\mathbb{P}\mathcal{C}_{fill}, d)$.
\begin{itemize}
    \item For any $\bz$ in $\mathcal{C}$ we define 
\[
Q(\mathbf{z}) 
= \sup_{\bc \in \mathcal{C}} \frac{i(\mathbf{z},\bc)}{i(\mathbf{b},\bc)}\,\frac{1}{h(\mathbf{b})} 
= \sup_{\bc \in \mathbb{P}\mathcal{C}} \frac{i(\mathbf{z},\bc)}{i(\mathbf{b},\bc)}\,\frac{1}{h(\mathbf{b})}
\]
The entropy of any current representing a point in Teichmüller space is equal to 1 by Equation~\ref{eq:entropy1}. This allows us to rewrite $Q(\textbf{z})$ as 
\[
\displaystyle{Q(\bz) =
\sup_{c \in \mathbb{P}\mathcal{C}}\frac{i(\bz,\bc)}{i(\bb,\bc)}}.
\]
Since the supremum is taken over a compact set and the denominator is never zero the supremum is always achieved and is positive and finite. 

\item For any $z \in \mathbb{P}\mathcal{C}$, we define 
\begin{center}
   $\mathcal{L}_{z}: \mathbb{P}\mathcal{C}: \rightarrow \mathbb{R}_{+}$\\
   $\displaystyle{\mathcal{L}_{z}(\mu) = \frac{i(\textbf{z}, \mu)}{Q(\textbf{z})}}$.
\end{center}
To define $\mathcal{L}_{z}$, we chose a representative $\bz$ of $z$ in the space of currents. But, since $\displaystyle{Q(k \textbf{z}) = k Q(\textbf{z})}$ for any positive real number $k$, $\mathcal{L}_{z}$ does not depend on the choice of the representative.
\end{itemize}
Now, for any $z \in \mathbb{P}\mathcal{C}$, we define $\Psi_{z}: \mathbb{P}\mathcal{C}_{fill} \rightarrow \mathbb{R}$ such that for all $x \in \mathbb{P}\mathcal{C}_{fill}$
\begin{equation}\label{eq3}
    \displaystyle{\Psi_{z}(x) = \log \Bigg( \sup_{c \in \mathbb{P}\mathcal{C}} \frac{\mathcal{L}_{z}(\bc)}{i(\bx, \bc)}\frac{1}{h(\bx)} \Bigg)}.
\end{equation}
The above function is well-defined on the projective class of $x$ since, for any positive constant $k$, we have, by Lemma~2.1~\ref{lem:entropy_properties},
\[
i(k\bx, \bc) \, h(k\bx) = k \, i(\bx, \bc) \cdot \frac{1}{k} h(\bx) = i(\bx, \bc) \, h(\bx).
\]

\begin{lemma}
    For any $z \in \mathbb{P}\mathcal{C}_{fill}$ the map $\Psi_{z}$ defined in \cref{eq3} satisfies $\Psi_z(x) = d(x,z) - d(b,z)$ for all $x \in \mathbb{P}\mathcal{C}_{fill}$.
\end{lemma}

\begin{proof}
Consider any $z \in \mathbb{P}\mathcal{C}_{fill}$. Then we have 
   \begin{align*}
    \displaystyle
    \Psi_{z}(x) &= \log \Bigg( \sup_{c \in \mathbb{P}\mathcal{C}} \frac{\mathcal{L}_{z}(\bc)}{i(\bx, \bc)}\frac{1}{h(\bx)} \Bigg).
\end{align*}
Since $z$ is filling, by Lemma~\ref{lem:finite_entropy}, any representative $\bz$ of $z$, $h(\bz)$ is finite. Therefore
substituting $\displaystyle{\mathcal{L}_z(c)} = \frac{i(\textbf{z}, \bc)}{Q(\bz)}$ and $Q(\textbf{z})$ as $\displaystyle{\sup_{c \in \mathbb{P}\mathcal{C}} \frac{i(\textbf{z}, \bc)}{i(\textbf{b}, \bc)}}$ in the above equation and simplifying gives us 
\begin{align*}
& \log \Bigg( \sup_{c \in \mathbb{P}\mathcal{C}} \frac{i(\textbf{z}, \bc) h(\bz)}{i(\textbf{x}, \bc) h(\bx)} \Bigg) -\log \Bigg( \sup_{c \in \mathbb{P}\mathcal{C}}\frac{i(\textbf{z},\bc)}{i(\textbf{b},\bc)}\frac{h(\bz)}{h(\bb)} \Bigg)\\
= & d(x,z) - d(b,z).
\end{align*} 
\end{proof}

\subsection{The map \texorpdfstring{$\Psi$}{Psi} is injective}\label{sec2.5}

We prove $\Psi$ is injective by considering two separate cases. To establish the two cases we make the following definitions:

\begin{definition}
    We say geodesic currents $\textbf{x}$ and $\textbf{y}$ have the same null support if 
    \begin{equation}
        i(\textbf{x}, \bc) = 0 \iff i(\textbf{y}, \bc) = 0 \text{ for all geodesic currents }\bc
    \end{equation}
    and we will denote it by $|\textbf{x}|_{0} = |\textbf{y}|_{0}$.
\end{definition}
Note that the following property is invariant for the projective classes of $\textbf{x}$ and $\textbf{y}$. This allows us to define projective currents with the same null support. 
\begin{definition}
 We say projective currents $x$ and $y$ have the same null support if
    \begin{equation}
        i(\textbf{x}, c) = 0 \iff i(\textbf{y}, c) = 0 
    \end{equation}
For any representative $\bx$ of $x$ and $\by$ of $y$ and for any geodesic currents $\bc \in \mathcal{C}$    
We will denote it by $|x|_{0} = |y|_{0}$.
\end{definition}

\subsubsection{\underline{Case 1: Projective currents with distinct null support}}
\begin{lemma}\label{lemma2.5}
    Let $x$, $y$ be two elements in $\mathbb{P}\mathcal{C}$ such that $|x|_{0} \neq |y|_{0}$. Then $\Psi_{x} \neq \Psi_{y}$.
\end{lemma}

\begin{proof} It is known that the map $\Psi:\mathbb{P}\mathcal{C} \rightarrow C(\mathbb{P}\mathcal{C}_{fill})$ is injective on $\mathbb{P}\mathcal{C}_{fill}$ \cite[Proposition 2.1]{CW2015}.

We are interested in the case when at least one of $x$ or $y$ in non-filling. Assume without loss of generality that $x$ is non-filling. Let $\bx$ and $\by$ be representatives of $x$ and $y$ in the space of currents. By our initial assumption, we can find a geodesic current $\boldsymbol{\mu}$ such that $i(\textbf{x}, \boldsymbol{\mu}) = 0$ and $i(\textbf{y}, \boldsymbol{\mu}) \neq 0$. To prove that $\Psi_{x} \neq \Psi_{y}$ we need to find an element $p \in \mathbb{P}\mathcal{C}_{fill}$ where the maps $\Psi_{x}$ and $\Psi_{y}$ disagree. 

Let  $\bnu$ be a filling multicurve. Consider the sequence $\textbf{p}_{n} = \textbf{x} + \frac{1}{n} \bnu$, which consists of filling currents converging to $\textbf{x}$. For each $n$, the following supremum is finite and is attained, as it is taken over a compact set and the denominator is never zero:
\[
\displaystyle 
    \sup_{c \in \mathbb{P}\mathcal{C}}\frac{i(\bx, \bc)}{Q(\bx)i(\textbf{p}_n, \bc)}.
\]
Using the properties of intersection form and supremum we get the following computation:
\begin{align}
    \displaystyle 
    \sup_{c \in \mathbb{P}\mathcal{C}}\frac{i(\bx, \bc)}{Q(\bx)i(\textbf{p}_n, \bc)} 
    &= \frac{1}{Q(\textbf{x})}\sup_{c \in \mathbb{P}\mathcal{C}}\frac{i(\bx, \bc)}{i(\textbf{p}_n, \bc)} \\ 
    &= \frac{1}{Q(\textbf{x})} \sup_{c \in \mathbb{P}\mathcal{C}}\frac{i(\bx, \bc)}{i(\bx, \bc) + \frac{1}{n} i(\bnu, \bc)} \notag \notag \\ 
    &= \frac{1}{Q(\textbf{x})} \Bigg( \frac{1}{1 + \frac{1}{n}\inf_{c \in \mathbb{P}\mathcal{C}}\frac{i(\bnu, \bc)}{i(\bx, \bc)}} \Bigg)\notag
\end{align} 
Since the infimum $\displaystyle{\inf_{c \in \mathbb{P}\mathcal{C}} \frac{i(\bnu, \bc)}{i(\bx, \bc)}}$ is taken over a compact set and $\bnu$ has non-zero intersection with all currents we have 
\[
\displaystyle{\inf_{c \in \mathbb{P}\mathcal{C}} \frac{i(\bnu, \bc)}{i(\bx, \bc)}} >0.
\]
Thus, we have
\[
\lim_{n \to \infty} \sup_{c \in \mathbb{P}\mathcal{C}}\frac{i(\bx, \bc)}{Q(\bx)i(\textbf{p}_n, \bc)} = \frac{1}{Q(\bx)}.
\]

On the other hand, since there exist a $\bm{\mu}$ such that $i(\bt{x},\bm{\mu}) = 0$ and $i(\bt{y}, \bm{\mu}) \neq 0$, we have $i(\bt{p}_n, \bm{\mu})$ approaches 0 as $n$ approaches infinity. This tell us  
\[
\displaystyle{\lim_{n \to \infty}\frac{1}{Q(\bt{y})}\sup_{c \in \mathbb{P}\mathcal{C}}\frac{i(\bt{y}, \bc)}{i(\bt{p}_n, \bc)} = \infty}
\]

We can pick a $N$ large enough so that ,
\[
 \displaystyle{
    \sup_{\lambda \in\mathbb{P} \mathcal{C}}\frac{i(\bt{y}, \bc)}{Q(\bt{y})i(\bt{p_N}, \bc)} \neq \sup_{c \in \mathbb{P}\mathcal{C}}\frac{i(\bt{x}, \bc)}{Q(\bt{x})i(\bt{p_N}, \bc)} .
    }
 \]
 Subsequently $\Psi_{y}(p_N) \neq \Psi_{x}(p_N)$.
\end{proof}

\subsubsection{\underline{Case 2: Projective currents with the same null support}}

In \cref{lemma2.8} we prove that for $x$, $y$ in $\mathbb{P}\mathcal{C}_{fill}$ satisfies $|x|_{0} = |y|_{0}$ then $\Psi_{x} \neq \Psi_{y}$. In order to prove that we introduce the following definition and \cref{Lemma2.7}.

\begin{definition}\label{def2.6}
    For a geodesic current $\textbf{x}$ we will use $\mathcal{C}_{x}$ to denote the collection of all geodesic currents that has non zero intersection with $\textbf{x}$. So,
    \[
    \mathcal{C}_{x} \coloneqq \{ \bm{\lambda} \in \mathcal{C} \mid i(\textbf{x} , \bm{\lambda}) \neq 0 \}.
    \]
    For every element in the projective class of $\textbf{x}$, the set $\mathcal{C}_{x}$ is the same.
\end{definition}

\begin{lemma}\label{Lemma2.7}
Let $\textbf{x}$, $\textbf{y}$ be two distinct geodesic currents such that $|\textbf{x}|_{0} = |\textbf{y} |_{0}$. Then  
\[
\text{either }
\displaystyle{\sup_{\lambda \in \mathcal{C}_{y}}\frac{i(\textbf{x}, \lambda)}{i(\textbf{y}, \lambda)} > 1}
\text{ or }
\displaystyle{\sup_{\lambda \in \mathcal{C}_{x}}\frac{i(\textbf{y}, \lambda)}{i(\textbf{x}, \lambda)} > 1}.
\]
\end{lemma}

\begin{proof} Assume there exists distinct geodesic currents $\textbf{x}$ and $\textbf{y}$ that satifies the following inequalities 
\[
\displaystyle{\sup_{\blambda \in \mathcal{C}_{y}}\frac{i(\textbf{x}, \blambda)}{i(\textbf{y}, \blambda)} \leq 1} \text{ and }    \displaystyle{\sup_{\blambda \in \mathcal{C}_{x}}\frac{i(\textbf{y}, \blambda)}{i(\textbf{x}, \blambda)} \leq 1}
\]
Since $|\textbf{x}|_{0} = |\textbf{y}|_{0}$ we have $\mathcal{C}_{x} = \mathcal{C}_{y}$.

The first inequality implies that $\displaystyle{i(\textbf{x}, \blambda) \leq i(\textbf{y}, \blambda)}$ for all $\blambda \in \mathcal{C}_{y}$. Similarly, 
 from the second inequality we get  $\displaystyle{i(\textbf{y}, \blambda) \leq i(\textbf{x}, \blambda)}$ for all $\blambda \in \mathcal{C}_{x}$. This implies 
 \begin{center}
 $\displaystyle{i(\textbf{y}, \blambda) = i(\textbf{x}, \blambda)}$ for all $\blambda \in \mathcal{C}_{x}$.
 \end{center}
 But, we also have 
 \begin{center}
 $\displaystyle{i(\textbf{y}, \blambda) = 0 = i(\textbf{x}, \blambda)}$ for all $\blambda \in \mathcal{C} \setminus \mathcal{C}_{x}$.
 \end{center}
Thus, the currents $\bx$ and $\by$ have the same intersection number with all closed curves. However, by a result of Otal \cite[Th\'eor\`eme~2]{Otal90:SpectreMarqueNegative}, this implies that the geodesic currents $\bx$ and $\by$ are equal, which yields a contradiction.
\end{proof}

\begin{lemma}\label{lemma2.8}
 Let $x$, $y$ be two elements in $\mathbb{P}\mathcal{C}$ such that $|x|_{0} = |y|_{0}$. Then $\Psi_{x} \neq \Psi_{y}$.
\end{lemma}

\begin{proof} As discussed earlier, the map $\Psi: \mathbb{P}\mathcal{C} \rightarrow C(\mathbb{P}\mathcal{C}_{fill})$ is known to be injective on $\mathbb{P}\mathcal{C}_{fill}$. Suppose at least one of $x$ or $y$ is non-filling. Since $|x|_{0} = |y|_{0}$, it follows that both $x$ and $y$ must be non-filling.

Pick a representative $\textbf{x}$ of $x$ and $\textbf{y}$ of $y$ in the space of currents. Let $\textbf{x}' = Q(\textbf{y})\textbf{x}$ and let $\textbf{y}' = Q(\textbf{x})\textbf{y}$. By \cref{Lemma2.7} we have $\displaystyle{\sup_{\lambda \in \mathcal{C}_{x}}\frac{i(\textbf{y}', \lambda)}{i(\textbf{x}', \lambda)} > 1}$ or $\displaystyle{\sup_{\lambda \in \mathcal{C}_{y}}\frac{i(\textbf{x}', \lambda)}{i(\textbf{y}', \lambda)} > 1}$. Without loss of generality, we assume the former.

Let $\bnu$ be a filling multicurve. Consider the sequence $\textbf{p}_{n} = \textbf{x} + \frac{1}{n} \bnu$ in $\mathbb{P}\mathcal{C}_{fill}$. We claim that for $n$ big enough $\Psi_{y}(p_{n}) \neq \Psi_{x}(p_{n})$.

To see this, take 
\begin{align*}
 \Psi_{y}(p_{n}) - \Psi_{x}(p_{n})  & = \displaystyle{\log \Bigg( \sup_{c \in \mathbb{P}\mathcal{C}} \frac{i(\textbf{y},\bc)}{Q(\textbf{y})i(\textbf{p}_{n}, \bc)}\frac{1}{h(\textbf{p}_{n})} \Bigg) - \log \Bigg( \sup_{c \in \mathbb{P}\mathcal{C}} \frac{i(\textbf{x},\bc)}{Q(\textbf{x})i(\textbf{p}_{n}, \bc)}\frac{1}{h(\textbf{p}_{n})} \Bigg)}.\\
\end{align*}
Since $Q(\textbf{x})$, $Q(\textbf{y})$ and $h(\textbf{p}_{n})$ are constants for all $n$ it does not affect the supremum. Therefore the above expression equals  
\begin{align}\label{eq8}
\log\displaystyle{\Bigg( \frac{\frac{Q(\textbf{x})}{Q(\textbf{y})}\sup_{c \in \mathbb{P}\mathcal{C}} \frac{ i(\textbf{y},\bc)}{i(\textbf{p}_{n}, \bc)}}{\sup_{c \in \mathbb{P}\mathcal{C}} \frac{i(\textbf{x},\bc)}{i(\textbf{p}_{n}, \bc)}}\Bigg)}.
\end{align}
The numerator $\displaystyle{\frac{Q(\textbf{x})}{Q(\textbf{y})} \sup_{c \in \mathbb{P}\mathcal{C}}\frac{i(\textbf{y}, \bc)}{i(\textbf{p}_n, \bc)}} =  \displaystyle{ \sup_{c \in \mathbb{P}\mathcal{C}}\frac{i(\textbf{y'}, \bc)}{Q(\textbf{y})i(\textbf{p}_n, \bc)}}$ converges to $\displaystyle{\sup_{c\in \mathbb{P}\mathcal{C}}\frac{i(\textbf{y'}, \bc)}{i(\textbf{x'}, \bc)}}$.  

We know that 
\[
\displaystyle{\sup_{c\in \mathbb{P}\mathcal{C}}\frac{i(\textbf{y'}, \bc)}{i(\textbf{x'}, \bc)}} \geq \displaystyle{\sup_{\bc \in \mathcal{C}_{x}}\frac{i(\textbf{y}', \lambda)}{i(\textbf{x}', \lambda)} > 1}
\]

We can pick $N$ large enough, so that numerator of the quantity in \cref{eq8} is greater than 1.

Now, we will run a similar argument as in the proof of the \cref{lemma2.5} to see that the denominator in \cref{eq8} is less than 1.
\begin{align*}
    \displaystyle{\sup_{c \in \mathbb{P}\mathcal{C}} \frac{i(\textbf{x},\bc)}{i(\textbf{p}_{n}, \bc)}}
    &= \displaystyle{\sup_{c \in \mathbb{P}\mathcal{C}} \frac{i(\textbf{x},\bc)}{i(\textbf{x}, \bc) + \frac{1}{n}i(\bnu, \bc)}}
    = \displaystyle{ \frac{1}{1 + \frac{1}{n} \inf_{c \in \mathbb{P}\mathcal{C}}\frac{i(\bnu, \bc)}{i(\textbf{x}, \bc)}}}
\end{align*}
Since the infimum $\displaystyle{\inf_{c \in \mathbb{P}\mathcal{C}} \frac{i(\bnu, \bc)}{i(\bx, \bc)}}$ is taken over a compact set $\mathbb{P}\mathcal{C}$ and $\bnu$ has non-zero intersection with all currents we get 
\[
\inf_{c \in \mathbb{P}\mathcal{C}} \frac{i(\bnu, \bc)}{i(\bx, \bc)} > 0 
\] 

This implies that the denominator of the quantity in \cref{eq8}  approaches 1 from the left hand side. This tell us that the for all $n$,
$ \displaystyle{\sup_{c \in \mathbb{P}\mathcal{C}} \frac{i(\bx,\bc)}{i(\textbf{p}_{N}, \bc)}
< 1}$.

To conclude, for $N$ large enough we get
\[
\Psi_{y}(p_{N}) - \Psi_{x}(p_{N}) = \log\displaystyle{\Bigg( \frac{\frac{Q(\textbf{x})}{Q(\textbf{y})}\sup_{c \in \mathbb{P}\mathcal{C}} \frac{ i(\textbf{y},\bc)}{i(\textbf{p}_{N}, \bc)}}{\sup_{c \in \mathbb{P}\mathcal{C}} \frac{i(\textbf{x},\bc)}{i(\textbf{p}_{N}, \bc)}}\Bigg)} > 0
\]
\[
\Psi_{y}(p_{N}) \neq \Psi_{x}(p_{N})
\]
\end{proof}

\begin{proposition}\label{injectivity} The map $\Psi: \mathbb{P}\mathcal{C} \rightarrow C(\mathbb{P}\mathcal{C}_{fill})$ is injective.
\end{proposition}

\begin{proof}
   The proposition follows directly from \cref{lemma2.5} and \cref{lemma2.8}
\end{proof}

\subsection{The map \texorpdfstring{$\Psi$}{Psi} is an embedding}\label{sec2.6}
The following two lemmas follow from the proofs of Proposition 3.2 and Lemma 3.5 in \cite{CW2015}. Walsh establishes these results in the setting of Thurston's compactification of Teichmüller space, $\mathcal{T}^{\mathcal{T}}$, rather than in $\mathbb{P}\mathcal{C}$. However, the proofs do not rely on any properties specific to Teichmüller space and remain valid when $\mathcal{T}^{\mathcal{T}}$ and $\mathbb{P}\mathcal{ML}$ are replaced with $\mathbb{P}\mathcal{C}$ in the argument.

\begin{lemma}
    If a sequence $x_{n}$ in $\mathbb{P}\mathcal{C}$ converges to a point $x \in \mathbb{P}\mathcal{C}$, then $\mathcal{L}_{x_n}$ converges to $\mathcal{L}_{x}$ uniformly on compact set of $\mathbb{P}\mathcal{C}$.
\end{lemma}

\begin{lemma}\label{lemma2.10} The map $\Psi: \mathbb{P}\mathcal{C} \rightarrow C(\mathbb{P}\mathcal{C}_{fill})$ that takes $z$ to $\Psi_{z}$ is continuous.
\end{lemma}

\begin{reptheorem}{horoboundary}
The horofunction compactification of $\mathbb{P}\mathcal{C}_{fill}$ equipped with the extended Thurston metric is homeomorphic to the space of projective currents $\mathbb{P}\mathcal{C}$. Moreover, the horofunction boundary corresponds precisely to the space of projective non-filling currents. \end{reptheorem}

\begin{proof}
By \cref{injectivity} and \cref{lemma2.10}, the map $\Psi: \mathbb{P}\mathcal{C} \to C(\mathbb{P}\mathcal{C}_{fill})$ is a continuous injection from a compact space to a Hausdorff space, and hence a homeomorphism onto its image.
\[
\mathbb{P}\mathcal{C} \cong \Psi(\mathbb{P}\mathcal{C}).
\]

The horofunction compactification of $\mathbb{P}\mathcal{C}_{fill}$ is the closure of $\Psi(\mathbb{P}\mathcal{C}_{fill})$ in $C(\mathbb{P}\mathcal{C}_{fill})$. Since $\Psi(\mathbb{P}\mathcal{C})$ is a compact subset of a Hausdorff space it is closed. It follows that
\[
\overline{\Psi(\mathbb{P}\mathcal{C}_{fill})} \subseteq \Psi(\mathbb{P}\mathcal{C})
\]

On the other hand, by continuity, $\Psi^{-1}(\overline{\Psi(\mathbb{P}\mathcal{C}_{fill})})$ is closed and contains $\mathbb{P}\mathcal{C}_{fill}$, so
\[
\Psi(\mathbb{P}\mathcal{C}) \subseteq \overline{\Psi(\mathbb{P}\mathcal{C}_{fill})}   
\]

Thus, $\Psi(\mathbb{P}\mathcal{C}) = \overline{\Psi(\mathbb{P}\mathcal{C}_{fill})}$, and the horofunction compactification of $\mathbb{P}\mathcal{C}_{fill}$ is homeomorphic to $\mathbb{P}\mathcal{C}$.
\end{proof}

\begin{remark}\label{extendedcompactification}
    Walsh proved that the horofunction boundary of Teichmüller space $\mathcal{T}$ with respect to the Thurston metric $d_{Th}$ is homeomorphic to the space of projective measured laminations $\mathbb{P}\mathcal{ML}$ \cite{CW2015}. Since projective measured laminations are contained in the space of projective non-filling currents, it follows that the horofunction boundary of $(\mathcal{T}, d_{Th})$ is contained in the horofunction boundary of $(\mathbb{P}\mathcal{C}_{fill}, d)$.

Recall that by Proposition~\ref{prop:d_extends_thurston}, for any $x, y \in \mathcal{T}$, we have $d_{Th}(x, y) = d(x, y)$. As a consequence, for each $y \in \mathcal{T}$, the horofunction map $\Psi_y: \mathcal{T} \to C(\mathcal{T})$ associated to the Thurston metric is the restriction of the horofunction map $\Psi_y: \mathbb{P}\mathcal{C}_{fill} \to C(\mathbb{P}\mathcal{C}_{fill})$ associated to the metric $d$.
\end{remark}




\section{All projective non-filling currents are Busemann points}
\label{sec:busemann}
From the previous section, we know that the horofunction boundary of $\mathbb{P}\mathcal{C}_{fill}$ consists of projective non-filling currents. We now show that every non-filling projective current can be approached along a geodesic path in $\mathbb{P}\mathcal{C}_{fill}$. Given a base point $b \in \mathcal{T}$ and a projective non-filling current $z$, in \cref{busemann2} we exhibit an example of an unparametrized geodesic ray $\gamma(t): [0, \infty) \to \mathbb{P}\mathcal{C}_{fill}$ that approaches $z$ as $t$ goes to infinity. The fact that every point on the horoboundary of $\mathbb{P}\mathcal{C}_{fill}$ is a Busemann point is a direct consequence of this lemma. 

First, we extend the definition of dilation $\Dil(\bx, \by)$ to allow one of the currents to be non-filling.

\begin{definition}
    Let $\bx \in \mathcal{C}_\mathrm{fill}$ and $\by \in \mathcal{C}$. We define $\Dil(\bx, \by)$ as 
\[
\Dil(\bx, \by) =  \sup_{c\in \mathbb{P}\mathcal{C}} \frac{i(\by, \bc)}{i(\bx, \bc)}.
\]
Note that since $\bx$ is filling the above supremum is always attained and is finite as the supremum is taken over the compact space $\PC$ and the denominator is never zero.
\end{definition}

The following useful fact is used in proving \cref{busemann2}.

\begin{lemma}\label{busemann1}Let $\bz$ be a non-filling geodesic current and let $\bb \in \mathcal{T}$. Consider the path $\bgamma(t): \mathbb{R}_{+} \to \mathcal{C}_{fill}$ given by $\bgamma(t) = \bb + t\bz$. For any $t \in \mathbb{R}_{+}$, $\Dil(\bb,\bgamma(t)) = 1 + \Dil(\bb, \bz)$.
\end{lemma}

\begin{proof} We compute:
    \[
\displaystyle{\Dil(\bb, \bgamma(t)) = \sup_{c \in \mathbb{P}\mathcal{C}}\frac{i(\bgamma(t), \bc)}{i(\bb, \bc)} = \sup_{c \in \mathbb{P}\mathcal{C}}\frac{i(\bb, \bc) + t. i(\bz ,\bc)}{i(\bb, \bc)} = \sup_{c \in \mathbb{P}\mathcal{C}}\Bigg( 1 + t\frac{i(\bz, \bc)}{i(\bb, \bc)}\Bigg) = 1 + t \Dil(\bb, \bz)}.
\]
This proves the lemma.
\end{proof}

Now we construct unparameterized geodesic rays.
Note that Cantrell--Reyes already construct bi-infinite geodesics within the broader space of pseudo-metrics in~\cite{StephenEduardo}, but their corresponding rays are not naturally contained in the space of geodesic currents.
Hence, the following construction is necessary.

\begin{lemma}\label{busemann2}
Let $z \in \mathbb{P}\mathcal{C}$ be a projective non-filling current, and let $\bb \in \mathcal{T}$. Fix a representative $\bz$ of $z$, and define a curve $\bgamma: \mathbb{R}_{+} \to \mathcal{C}_{fill}$ by
\[
\bgamma(t) = \bb + t\bz.
\]
Then the projective class $\gamma(t) = [\bb + t\bz]$ defines a geodesic $\gamma: \mathbb{R}_{+} \to \mathbb{P}\mathcal{C}_{fill}$ that converges to $z$ as $t \to \infty$.
\end{lemma}

\begin{proof}
Since the following limit holds in the space of currents
\[
\displaystyle{\lim_{t \to \infty} \frac{1}{t}( \textbf{b} + t \textbf{z}) = \lim_{t \to \infty} \frac{\textbf{b}}{t} + \textbf{z} = \textbf{z}}
\]
we have $\displaystyle{\lim_{t \to \infty} \gamma(t)} = z$ in projective currents. It remains to show that $\gamma(t)$ is a geodesic. For that we want to show that $\gamma(t)$ satisfies the following relation
\begin{equation}
    d(b, \gamma(s)) + d(\gamma(s), \gamma(t)) = d(b, \gamma(t)), \text{ where } s, t \in \mathbb{R}_{+} \text{ with } s \leq t.
\end{equation}
When we use the definition of Thurston's extended metric and simplify using the properties of logarithm, all the terms involving entropy cancel out and we get the following equivalent equation:
\begin{equation}\label{Dileq1}
    \Dil(\bb, \bgamma(s)). \Dil(\bgamma(s), \bgamma(t)) = \Dil(\bb, \bgamma(t))
\end{equation}

In the rest of the proof we show \cref{Dileq1} holds.

A simplification of the triangle inequality  $d(b, \gamma(s)) + d(\gamma(s), \gamma(t)) \geq d(b, \gamma(t))$ tells us that
\begin{equation}\label{Dileq2}
    \Dil(\bb, \bgamma(s)). \Dil(\bgamma(s), \bgamma(t)) \geq \Dil(\bb, \bgamma(t)) 
\end{equation}
We must show that the above equation is an equality. As a consequence of \cref{Dileq2} and \cref{busemann1} we get 
\[
\Dil(\bgamma(s), \bgamma(t)) \geq \frac{\Dil(\bb, \bgamma(t))}{ \Dil(\bb, \bgamma(s))} = \frac{1 + t \Dil(\bb,\bz)}{1 + s \Dil(\bb,\bz)}.
\]
Assume \cref{Dileq2} is a strict inequality, that is
\begin{equation}\label{inequality}
 \Dil(\bgamma(s), \bgamma(t)) > \frac{1 + t \Dil(\bb,\bz)}{1 + s \Dil(\bb,\bz)}.  
\end{equation}

Rewriting $\Dil(\bgamma(s), \bgamma(t))$ in terms of $s$ and $t$, we obtain:
\[
    \Dil(\bgamma(s), \bgamma(t)) = \sup_{c\in \mathbb{P}\mathcal{C}} \frac{i(\bgamma(t), \bc)}{i(\bgamma(s), \bc)} =  \sup_{c \in \mathbb{P}\mathcal{C}} \frac{1 + t\frac{i(\bz,\bc)}{i(\bb,\bc)}}{1 + s\frac{i(\bz,\bc)}{i(\bb,\bc)}}.
\]
Substituting this into \cref{inequality} we get
\[
\sup_{c \in \mathbb{P}\mathcal{C}} \frac{1 + t\frac{i(\bz,\bc)}{i(\bb,\bc)}}{1 + s\frac{i(\bz,\bc)}{i(\bb,\bc)}} > \frac{1 + t \Dil(\bb,\bz)}{1 + s \Dil(\bb,\bz)}.
\]

Let us say the supremum on the left is achieved by some $\blambda_{s}$. Then,
\begin{equation*}
    \displaystyle\frac{1 + t\frac{i(\bz,\blambda_s)}{i(b,\blambda_s)}}{1 + s\frac{i(\bz,\blambda_s)}{i(\bb,\blambda_s)}} > \frac{1 + t \Dil(\bb,\bz)}{1 + s \Dil(\bb,\bz)}. \nonumber 
\end{equation*}

Distributing and simplifying the above inequality we get
\begin{align*}
\frac{i(\bz,\blambda_s)}{i(\bb,\blambda_s)} &> \Dil(\bb,\bz).
\end{align*}
However, $\Dil(\bb,\bz)$ is defined to be a supremum and hence cannot be strictly less than the left-hand side. This forces \cref{Dileq2} to be an equality, and thus we conclude that \cref{Dileq1} holds, proving that 
$\gamma(t)$ is a geodesic. 
\end{proof}

\begin{remark} The argument above does not use the fact that the pseudo-metrics arise from geodesic currents, and it works more generally in the setting of general pseudo-metrics of a non-elementary hyperbolic group (see Appendix~\ref{sec:asymmetric} for definitions in that framework).
\end{remark}

\begin{theorem}\label{busemann3}
All projective non-filling currents are Busemann points.
\end{theorem}
\begin{proof}
The theorem is direct consequence of \cref{horoboundary} and \cref{busemann1}.
\end{proof}

\section{Characterizing finite detour cost on non-filling projective currents}\label{sec:detourcost}

For any pair of points $\eta$ and $\xi$ on the horofunction boundary, the symmetrized detour cost $\delta$ is defined by:
\[
\delta(\eta, \xi) = H(\eta, \xi) + H(\xi, \eta).
\]

In this section, we study the detour cost $H(\xi, \eta)$ and the symmetrized detour cost $\delta(\xi, \eta)$ for pairs of non-filling projective currents $\xi$ and $\eta$.

\cref{detourcost1} describes $H(\xi, \eta)$ in terms of the intersection form. In \cref{finiteDC}, we provide a necessary condition that $\xi$ and $\eta$ must satisfy for $H(\xi, \eta)$ to be finite. Building on this and the decomposition of geodesic currents from the previous section, \cref{finitedetorcost2} shows that non-filling projective currents $\xi$ and $\eta$ with
\[
\delta(\xi, \eta) < \infty
\]
must fill the same subsurface of $S$ and have measured lamination components with the same support.

\begin{lemma}\label{detourcost1}
    Let $\xi$ and $\eta$ be non-filling projective currents, and let $b$ be a basepoint in $\mathcal{T}$. Then the detour cost is given by
\begin{align}\label{detourcost2}
H(\xi, \eta) 
= \inf_{\bgamma} \liminf_{t \to \infty} \Bigg( 
    & \log\left( \sup_{c \in \mathbb{P}\mathcal{C}} \frac{i(\bgamma(t), \bc)}{i(\bb, \bc)} \right)
    + \log\left( \sup_{c \in \mathbb{P}\mathcal{C}} \frac{i(\boldeta, \bc)}{i(\bgamma(t), \bc)} \right) \nonumber\\
    & \quad - \log\left( \sup_{c \in \mathbb{P}\mathcal{C}} \frac{i(\boldeta, \bc)}{i(\bb, \bc)} \right)
\Bigg),
\end{align}
where $\bxi$, $\boldeta$, and $\bb \in \mathcal{C}$ are representatives of $\xi$, $\eta$, and $b$, respectively. The infimum is taken over all paths $\bgamma: \mathbb{R}_{+} \to \mathcal{C}_{fill}$ such that $\bgamma(t)$ converges to a scalar multiple of $\bxi$ as $t \to \infty$. Note that the expression above is independent of the choice of representatives in the space of currents. 
\end{lemma}

\begin{proof} From the definition of detour cost in \cref{detourcosteq}
we know that 
\begin{equation*}
    H(\xi, \eta) = \displaystyle{\inf _{\gamma} \liminf_{t \to \infty}\Big(d(b, \gamma(t)) + \Psi_{\eta}(\gamma(t)) \Big)}
\end{equation*}
where the infimum is taken over all paths $\gamma(t): \mathbb{R}_{+} \to \mathbb{P}\mathcal{C}_{fill}$ converging to $\xi$. For each $t \in \mathbb{R}_{+}$, we can choose a representative of $\gamma(t)$ in the space of currents $\mathcal{C}$ such that $\gamma(t)$ converges to a representative of $\xi$ in $\mathcal{C}_{fill}$. We define $\bgamma(t): \mathbb{R}_{+} \rightarrow \mathcal{C}_{fill}$ to be this representative path. Substituting the expression for $\Psi_{\eta}(\gamma(t))$ from \cref{eq3} and simplifying, all terms involving entropy cancel out.
\begin{equation*}
     H(\xi, \eta) = \displaystyle{\inf _{\gamma} \liminf_{t \to \infty}\Bigg( \log \Bigg(\sup_{c \in \mathbb{P}\mathcal{C}} \frac{i(\bgamma(t), \bc)}{i(\bb,\bc)} \Bigg) + \log \Bigg(\sup_{c \in \mathbb{P}\mathcal{C}} \frac{i(\boldeta, \bc)}{i(\bgamma(t),\bc)} \Bigg)  - \log \Bigg(\sup_{c \in \mathbb{P}\mathcal{C}} \frac{i(\boldeta, \bc)}{i(\bb,\bc)} \Bigg) \Bigg)}.
\end{equation*}
For every ray $\bgamma: \mathbb{R}_{+} \rightarrow \mathcal{C}_{fill}$ converging to a scalar multiple of $\bxi$ we get a ray $\gamma: \mathbb{R}_{+} \rightarrow \mathbb{P}\mathcal{C}_{fill}$ converging to $\xi$. Thus, the above expression is equivalent to \cref{detourcost2}, completing the proof of the lemma.
\end{proof}

\begin{remark}[Notational Convention]\label{Notations2}
Throughout the rest of the paper, we use $x$ instead of $\bx$ if it appears in a statement that does not depend on the choice of representative of $x$. For example, the expression $\sup_{c \in \mathbb{P}\mathcal{C}} \frac{i(\boldsymbol{x},c)}{i(\boldsymbol{y},c)}$ is independent of the choice of representative for $c$, so we write $c$ rather than $\bc$. 

Another example appears when discussing zero versus non-zero intersection numbers. If $i(\bx, \by) = 0$ for some representatives $\bx$ of $x$ and $\by$ of $y$, then $i(\bx', \by') = 0$ for any other choice of representatives. Thus, we simply write $i(x, y) = 0$.

This convention is particularly relevant when discussing finiteness properties: the finiteness of most the quantities we consider in this paper does not depend on the choice of representatives. For instance, when we write
\[
\sup_{c \in \mathbb{P}\mathcal{C}} \frac{i(x,c)}{i(y,c)} < \infty,
\]
we mean that the supremum is finite regardless of the choice of \(\boldsymbol{x} \in [x]\), \(\boldsymbol{y} \in [y]\), and \(\boldsymbol{c} \in [c]\). Note that, while finiteness is independent of the choice of representatives, the actual numerical value of the supremum may vary depending on the representatives chosen.
\end{remark}

\subsection{Necessary conditions for finite detour cost}

\begin{lemma}\label{finitedetorcost1} Let $\xi$ and $\eta$ be non-filling geodesic currents and let $b \in \mathcal{T}$. The detour cost $H(\xi, \eta)$ is finite if and only if there exist some curve $\bgamma(t): \mathbb{R}_{+} \rightarrow \mathcal{C}_{fill}$ converging to a representative of $\xi$ in $\mathcal{C}$ for which 
\[
\displaystyle{\liminf_{t \to \infty}\sup_{c \in \mathbb{P}\mathcal{C}} \frac{i(\eta, c)}{i(\bgamma(t),c)} < \infty}.
\]
\end{lemma}

\begin{proof}
Note that the quantity $\displaystyle{\sup_{c \in \mathbb{P}\mathcal{C}} \frac{i(\boldeta, c)}{i(\bb,c)}}$ appearing \cref{detourcost2} is finite regardless of the choice of the representative $\bb$ of $b$ and $\boldeta$ of $\eta$. This is because the supremum is taken over a compact set and the denominator is bounded away from zero. 

Similarly, for any $t \in \mathbb{R_{+}}$ and for any $\bgamma(t): \mathbb{R}_{+} \to \mathcal{C}_{fill}$ converging to a representative $\bxi$ of $\xi$, the supremum 
$\displaystyle{\sup_{c \in \mathbb{P}\mathcal{C}} \frac{i(\bgamma(t), c)}{i(b,c)} }$ is always achieved. That is, there exist some $c_{t} \in \mathbb{P}\mathcal{C}$ such that
\[
\displaystyle{\sup_{c \in \mathbb{P}\mathcal{C}} \frac{i(\bgamma(t), c)}{i(\bb,c)} = \frac{i(\bgamma(t), c_t)}{i(\bb,c_t)}}.
\]
We can find a sequence of points $\{c_{n}\}$ contained in $\{c_{t} \mid t \in \mathbb{R_{+}}\}$ that converges to projective current $\kappa$. This tells us that 
\[
\displaystyle{\liminf_{n \to \infty} \sup_{c \in \mathbb{P}\mathcal{C}} \frac{i(\bgamma(n), c)}{i(\bb,c)} \leq \lim_{n \to \infty} \frac{i(\bgamma(n), c_n)}{i(\bb,c_n)} = \frac{i(\bxi, \kappa)}{i(\bb,\kappa)} < \infty }.
\]
Note that, the finiteness of the two supremum values discussed above does not depend of the choice of representative $\bxi$, $\bgamma(t)$, and $\bb$. Therefore by \cref{detourcost2} $H(\xi,\eta)$ is finite if and only if there exist some curve $\bgamma(t): \mathbb{R}_{+} \rightarrow \mathcal{C}_{fill}$ converging to a representative of $\xi$ in $\mathcal{C}$ for which 
\[
\displaystyle{\liminf_{t \to \infty}\sup_{c \in \mathbb{P}\mathcal{C}} \frac{i(\boldeta, c)}{i(\bgamma(t), c)} < \infty}.
\]
The finiteness of this supremum does of depend on the choice of curve $\bgamma(t)$ and $\boldeta$ in their respective projective class. We write this as 
\[
\displaystyle{\liminf_{t \to \infty}\sup_{c \in \mathbb{P}\mathcal{C}} \frac{i(\eta, c)}{i(\bgamma(t), c)} < \infty}.
\]
\end{proof}

\begin{lemma}\label{finitedetourcost1.2}
    Let $\xi$ and $\eta$ be non-filling geodesic currents and let $b \in \mathcal{T}$. The detour cost $H(\xi, \eta)$ is finite if and only if 
\[
\displaystyle{\sup_{\substack{c \in \mathbb{P}\mathcal{C}}}\frac{i(\eta, c)}{i(\xi,c)} \text{ is bounded above.}}  
\]
\end{lemma}

\begin{proof}
Let $\bxi$ and $\boldeta$ be representatives of $\xi$ and $\eta$ in the space of currents. Suppose that $H(\xi, \eta)$ is finite. Then by \cref{finitedetorcost1}, there exists a path $\bgamma: \mathbb{R}_{+} \to \mathcal{C}_{\text{fill}}$ converging to $\bxi$ such that
\[
    \liminf_{t \to \infty} \sup_{c \in \mathbb{P}\mathcal{C}} \frac{i(\boldeta, c)}{i(\bgamma(t), c)} < \infty.
\]

Since $\bgamma(t)$ is filling for each $t \in \mathbb{R}_{+}$, the supremum is achieved for each $t$. Let's say that for any $t \in \mathbb{R}_{+}$ supremum is achieved at $\bc_{t} \in \mathcal{C}$
\begin{equation}
    \frac{i(\boldeta, \bc)}{i(\bgamma(t), \bc)} \leq \frac{i(\boldeta, \bc_t)}{i(\bgamma(t), \bc_t)} \quad \text{for all } \bc \in \mathcal{C}.
\end{equation}

Since the $\liminf$ is finite, we can find a sequence $\{t_{n}\} \to \infty$ such that
\[
\lim_{n \to \infty}\frac{i(\boldeta, \bc_{t_n})}{i(\bgamma(t_n),\bc_{t_n})} \leq L
\]
for some positive constant $L$. As $\mathbb{P}\mathcal{C}$ is compact, we can find a subsequence $\{\bc_{m}\}$ of $\{\bc_{t_n}\}$ that converges projectively to some $\blambda \in \mathcal{C}$. That is, there exists a sequence of positive scalars $k_n$ such that
\[
    \lim_{n \to \infty} k_n \bc_n = \blambda.
\]
By continuity of the intersection form and convergence of $\bgamma(t_n) \to \bxi$, we obtain
\begin{align*}
   \frac{i(\boldeta, \bc)}{i(\bxi, \bc)} &\leq \frac{i(\boldeta, \blambda)}{i(\bxi,\blambda)} \leq L \text{ for all } \bc \in \mathcal{C}.
\end{align*}
Hence, 
\[
\sup_{c \in \mathbb{P}\mathcal{C}} \frac{i(\boldeta, \bc)}{i(\bxi, \bc)} \leq L.
\]
Since the choice of representatives $\boldeta$ and $\bxi$ is arbitrary we have:
\[
\sup_{c \in \mathbb{P}\mathcal{C}} \frac{i(\eta, \bc)}{i(\xi, \bc)} < \infty.
\]

Now, to prove the converse we assume the supremum is bounded above by a constant $L$, and so, for any $\bc \in   \mathcal{C}$ we have
    \begin{equation}\label{supeq2}
        i(\boldeta, \bc) \leq L\, i(\bxi, \bc).
    \end{equation}
Let $\bmu$ be a filling current. We define a path, $\bgamma(t): \mathbb{R}_{+} \rightarrow \mathcal{C}_{fill}$ by 
\[
\bgamma(t) = \bxi + \frac{1}{t} \bmu.
\]
Note that, $\bgamma(t)$ approaches $\bxi$ as $t \to \infty$. Each $\bgamma(t)$ is a filling current, so the supremum
\[
    \sup_{c \in \mathbb{P}\mathcal{C}} \frac{i(\boldeta, c)}{i(\bgamma(t), c)}
\]
is finite and attained. Let $c_t$ be a projective current that realizes the supremum at $t$. By \cref{supeq2} we get 
\begin{align*}
\sup_{ c \in \mathbb{P}\mathcal{C}} \frac{i(\boldeta,c)}{i(\bgamma(t),c)} = \frac{i(\boldeta, c_t)}{i(\bgamma(t), c_t)} &= \frac{i(\boldeta, c_t)}{i(\bxi, c_t) + \frac{1}{t} i(\bmu, c_t)} \leq \frac{L \cdot i(\bxi, c_t)}{i(\bxi, c_t) + \frac{1}{t} i(\bmu, c_t)} \leq \frac{L}{1 + \frac{1}{t} \cdot \frac{i(\bmu, c_t)}{i(\bxi, c_t)}}.
\end{align*}
Since the numerator is bounded above and the denominator exceeds $1$ for all $t$, it follows that
    \[
    \liminf_{t \to \infty} \sup_{c \in \mathbb{P}\mathcal{C}} \frac{i(\boldeta, c)}{i(\bgamma(t), c)} < \infty.
    \]
The proof does not depend on the choice of representatives $\boldeta$ and $\bxi$, and by \cref{finitedetorcost1} we have $H(\xi, \eta) < \infty$.
\end{proof}

\begin{lemma}\label{finiteDC}Let $\xi$ and $\eta$ be two non filling projective currents. If the detour cost $H(\xi, \eta)$ is finite then for all $x \in \mathbb{P}\mathcal{C}$, $i(\xi, x) = 0$ implies $i(\eta, x) = 0$.
\end{lemma}

\begin{proof} 
    Suppose there exists a $x \in \mathbb{P}\mathcal{C}$ such that $i(\xi, x) = 0$ but $i(\eta, x) \neq 0$. Then for any path $\bgamma(t): \mathbb{R}_{+} \to \mathcal{C}_{fill}$ converging to a representative $\bxi$ of $\xi$, we have
    \[
    \lim_{t \to \infty} i(\bgamma(t), x) = 0.
    \]
    Consequently,
    \[
    \lim_{t \to \infty} \frac{i(\eta, x)}{i(\bgamma(t), x)} = \infty.
    \]
    Note that,
    \begin{align*}
        \frac{i(\eta, x)}{i(\bgamma(t), x)} &\leq \sup_{c \in \mathbb{P}\mathcal{C}} \frac{i(\eta, c)}{i(\bgamma(t), c)} \text{ for all t},
     \end{align*}
     and therefore
     \begin{align*}
      \lim_{t \to \infty} \frac{i(\eta, x)}{i(\bgamma(t), x)} &\leq \liminf_{t \to \infty} \sup_{c \in \mathbb{P}\mathcal{C}} \frac{i(\eta, c)}{i(\bgamma(t), c)}.
     \end{align*}
        By \cref{finitedetorcost1}, it follows that $H(\xi, \eta) = \infty$. 
\end{proof}

\begin{proposition}\label{finitedetorcost2}
If the symmetrized detour cost $\delta(\xi, \eta)$ is finite for two non-filling projective currents $\xi$ and $\eta$, then for any representative $\boldeta$ of $\eta$ and $\bxi$ of $\xi$ in the space of currents, the following two conditions are satisfied:
\begin{enumerate}
\item The filling components $\bxi$ and $\boldeta$ fill the same subsurface $S'$ of $S$; and
\item The measured lamination components of $\bxi$ and $\boldeta$ have the same support.
\end{enumerate}
\end{proposition}

\begin{proof} From \cref{finiteDC} we know that if $\delta(\xi, \eta) = H(\xi, \eta) + H(\eta, \xi)$ is finite then following condition has to be satisfied:
\begin{equation}\label{iff}
    i(\eta, c) = 0 \iff i(\xi, c) = 0 \text{ for any geodesic current } c.
\end{equation}
Any two projective currents $\eta$ and $\xi$ satisfying \cref{iff} the set $\mathcal{E}_{\eta}$ and $\mathcal{E}_{\xi}$ defined \cref{Eset} are the same. Therefore, both $\eta$ and $\xi$ decomposes $S$ into the same collection of subsurfaces $\{S_i\}_{i =1}^{n}$. Now, if there is a subsurface $S_{i}$ that is filled by $\eta$ and not by $\xi$, then we can find a geodesic current $c$ intersecting exactly one of $\eta$ and $\xi$ contradicting \cref{iff}. Let $S'$ be the union of all $S_{i}$'s filled by both $\eta$ and $\xi$. The filling components $\xi$ and $\eta$ fill the same subsurface $S'$ of $S$.

To show that the measured lamination components of $\xi$ and $\eta$ have the same support, we pick a representative $\bxi$ of $\xi$ and $\boldeta$ of $\eta$ in the space of currents. Let $\boldeta = \boldeta_{S'} + \blambda_{\boldeta}$ and $\bxi = \bxi_{S'} + \blambda_{\bxi}$ denote the decompositions of these currents (as per Equation~\ref{decompostioneq}). The lamination components $\blambda_{\bxi}$ and $\blambda_{\boldeta}$ must be supported on the same subsurfaces $S_i$ of $S$. Otherwise, we could find a closed geodesic that intersects one lamination but not the other, contradicting \cref{iff}. Now, if $\boldeta$ fills a subsurface $S_i$, then by the decomposition we know that it intersects every curve in the interior of $S_i$. If $\blambda_{\boldeta}$ and $\blambda_{\bxi}$ do not share the same support, we would have $i(\blambda_{\boldeta}, \blambda_{\bxi}) \neq 0$, leading to
\[
i(\boldeta, \blambda_{\boldeta}) = i(\boldeta_{S'}, \blambda_{\boldeta}) + i(\blambda_{\boldeta}, \blambda_{\boldeta}) = 0,
\]
and
\[
i(\bxi, \blambda_{\boldeta}) = i(\bxi_{S'}, \blambda_{\boldeta}) + i(\blambda_{\bxi}, \blambda_{\boldeta}) \neq 0,
\]
which again contradicts \cref{iff}.
\end{proof}

\section{Characterizing projective minimal  measured laminations}\label{sccml}


In this section, we characterize minimal projective  measured laminations using the symmetrized detour cost $\delta$. Specifically, we show that a non-filling projective current $\eta$ is a minimal measured lamination if and only if the only non-filling projective current that is a finite $\delta$-distance from $\eta$ is $\eta$ itself. That is, $\eta$ is a minimal projective measured lamination if and only if any non-filling projective current $\xi$ satisfying
\[
\delta(\eta, \xi) < \infty
\]
must equal $\eta$.

A special case of minimal measured laminations is that of simple closed curves. When $\eta$ is a simple closed curve, it follows from \cref{finitedetorcost2} that any non-filling projective current $\xi$ with $\delta(\xi, \eta) < \infty$ must be projectively equivalent to $\eta$. In \cref{minimal ML}, we extend this result to arbitrary minimal projective measured laminations. In the remainder of the section, we prove that no other non-filling projective currents satisfy this rigidity property.


In \cref{deltaML} we consider a non-minimal projective measured lamination $\eta$ and provide an explicit example of another projective measured lamination finite $\delta$-distance away from $\eta$. \cref{deltafilling} establishes a similar result for non-filling projective currents with filling components. In particular, we show that if a non-filling projective current $\eta$ fills a subsurface $S' \subset S$, and $b$ is a boundary curve of $S'$, then $\delta(\eta, \eta + b) < \infty$. \cref{H_{Th}} compares the detour cost $H(\xi, \eta)$ defined on non-filling projective currents to $H_{Th}(\xi, \eta)$ defined in \cite{CW2015} on only projective measured laminations. We use this in the proof of \cref{deltaML}. \cref{finitedelta} follows directly from these lemmas.


The following definition is from \cite{CW2015}.
\begin{definition}
Let $\bxi$ and $\boldeta$ be measured laminations, and suppose
\[
\bxi = \sum_{i=1}^n \bxi_i
\]
is the decomposition of $\bxi$ into its minimal components $\bxi_i$. We write $\boldeta \ll \bxi$ if $\boldeta$ can be expressed as
\[
\boldeta = \sum_{i=1}^n k_i \bxi_i,
\]
where each $k_i$ is a non-negative real number.
\end{definition}

\begin{lemma}\label{minimal ML} Let $\eta$ be a minimal measured lamination. If a non-filling projective current $\xi$ satisfies $\delta(\xi, \eta) < \infty$, then $\xi = \eta$.
\end{lemma}

\begin{proof}
Assume $\xi$ is a non-filling projective current that satisfies $\delta(\xi, \eta) < \infty$. It follows from \cref{finitedetorcost2} that the decomposition of $\xi$ cannot contain any filling component. Hence, $\xi$ must be a projective measured lamination.

Let $\boldeta$ and $\bxi$ be measured laminations representing the projective classes $\eta$ and $\xi$, respectively. Suppose for contradiction that $\eta \neq \xi$. Since $\boldeta$ is minimal, and $\bxi$ is not projectively equivalent to $\boldeta$, it follows that
\[
\bxi \not\ll \boldeta.
\]
By \cite[Lemma~6.4]{CW2015}, this implies that the following supremum is infinite:
\[
\sup_{\lambda \in \mathbb{P}\mathcal{ML}} \frac{i(\bxi, \lambda)}{i(\boldeta, \lambda)} = \infty.
\]
However, since $\mathbb{P}\mathcal{ML} \subset \mathbb{P}\mathcal{C}$, we also have
\[
\sup_{c \in \mathbb{P}\mathcal{C}} \frac{i(\bxi, c)}{i(\boldeta, c)} \geq \sup_{\lambda \in \mathbb{P}\mathcal{ML}} \frac{i(\bxi, \lambda)}{i(\boldeta, \lambda)}
\]
This implies the supremum of the ratio taken over projective currents is infinite, and consequently by \cref{finitedetourcost1.2} 
\[
H(\xi, \eta) = \infty
\]
But this contradicts our assumption that $\delta(\eta, \xi) < \infty$
\end{proof}

\begin{lemma}\label{H_{Th}}
    For any two projective measured laminations $\xi$, $\eta$, the detour cost $H_{Th}(\xi, \eta)$ with respect to the Thurston metric on the Teichmüller space $(\mathcal{T}(S), d_{Th})$ is  greater than or equal to the detour cost $H(\xi, \eta)$ with respect to $(\mathbb{P}\mathcal{C}_{fill}, d)$, i.e.,
    \[
    H(\xi, \eta) \leq H_{Th}(\xi, \eta)
    \]
\end{lemma}

\begin{proof} 
By Proposition~\ref{prop:d_extends_thurston}, we know that the extended Thurston metric on projective filling currents agrees with the classical Thurston metric on the Teichmüller space. That is, $d_{Th}(x, y) = d(x, y)$ for all $x, y \in \mathcal{T}(S)$. The inequality stated in the lemma then follows from this fact and from the definition of detour cost given in \cref{detourcosteq}. Since $H(\xi, \boldeta)$ is defined as the infimum over a larger collection of paths $\gamma(t): \mathbb{R}_{+} \rightarrow \mathbb{P}\mathcal{C}_{fill}$ converging to $\xi$, $H(\xi, \eta) \leq H_{Th}(\xi, \eta)$.

To make this more precise, consider 
\begin{align*}
    H_{Th}(\xi, \eta) &= \inf_{\gamma \in \mathcal{T}} \liminf_{t \to \infty} \left( d_{Th}(b, \gamma(t)) + \eta(\gamma(t)) \right).
\end{align*}
Here, we have $b \in \mathcal{T}(S)$, and the infimum is taken over all paths $\gamma(t): \mathbb{R}_{+} \rightarrow \mathcal{T}$ converging to $\xi$. We have chosen the basepoint in the definition of the horofunction boundary to lie in $\mathcal{T}(S)$. As a consequence, and since $\gamma(t) \in \mathcal{T}(S)$ for each $t$, the value of $\eta(\gamma(t))$ is the same whether considered as a point in the horofunction boundary of $(\mathcal{T}(S), d_{Th})$ or of $(\mathbb{P}\mathcal{C}_{fill}, d)$. Hence, we have
\begin{align*}
    H_{Th}(\xi, \eta) &= \inf_{\gamma \in \mathcal{T}} \liminf_{t \to \infty} \left( d(b, \gamma(t)) + \eta(\gamma(t)) \right)\\
    &\geq \inf_{\gamma \in \mathbb{P}\mathcal{C}_{fill}} \liminf_{t \to \infty} \left( d(b, \gamma(t)) + \eta(\gamma(t)) \right)\\
    &= H(\xi, \eta)
\end{align*}
Since the infimum in the second equation is taken over a larger collection of paths $\gamma(t) \colon \mathbb{R}_{+} \rightarrow \mathbb{P}\mathcal{C}_{fill}$, it is the smaller.
\end{proof}

\begin{lemma}\label{deltaML}
Let $\boldeta$ be a measured lamination that decomposes as a finite sum of non-zero minimal laminations:
\[
\boldeta = \sum_{i = 1}^{n} \boldeta_{i}.
\]
For any choice of positive real numbers $k_i > 0$, consider the measured lamination
\[
\bxi = \sum_{i =1}^{n} k_{i} \boldeta_{i}.
\]
Then the projective classes $\eta = [\boldeta]$ and $\xi = [\bxi]$ satisfy $\delta(\xi, \eta) < \infty$.
\end{lemma}

\begin{proof} The projective measured laminations $\eta$ and $\xi$ satisfies $\eta \ll \xi$ and $\xi \ll \eta$. By \cite{CW2015} we know that 
\[
H_{Th}(\xi, \eta) < \infty \text{ and } H_{Th}(\eta, \xi) < \infty,
\]
and thus, by \cref{H_{Th}}, we get 
\[
H(\xi, \eta) < \infty \text{ and } H(\eta, \xi) < \infty.
\]
Therefore, $\delta(\xi, \eta) < \infty$ proving our lemma.
\end{proof}

\begin{lemma}\label{lemma:systolineqality} Let $\boldeta_{S'}$ be a geodesic current that fills a connected subsurface $S'$ of $S$, and let $b$ be a simple closed curve forming one of the boundary components of $S'$. Let $\bb$ denote the geodesic current representing the simple closed curve $b$. Then
\[
\sup_{c \in \mathbb{P}\mathcal{C}}\frac{i(\bb,c)}{i(\boldeta_{S'}, c)} \leq \frac{4}{\Syst_{S'}(\boldeta_{S'})}
\]
where $\Syst_{S'}(\boldeta_{S'})$ is defined in \cref{def:systole}.
\end{lemma}

\begin{proof}
Let $a$ be a closed curve on $S$ that has non-zero intersection with $b$, and let $\ba$ be its representative in the space of currents. We have,
\[
i(\bb, \ba) \neq 0
\]
In this proof, we show that any such $\ba$ satisfies the inequality
\[
i(\boldeta_{S'}, \ba) \geq \frac{i(\bb, \ba)}{4} \cdot \operatorname{Syst}_{S'}(\boldeta_{S'}).
\]

Since $b$ is a boundary curve of $S'$, it is contained in $\mathcal{E}_{\eta_{S'}}$. We recall the definition of  $\mathcal{E}_{\eta_{S'}}$ here:
\[
\mathcal{E}_{\boldeta_{S'}} :=
\begin{cases} 
      c \text{ is a closed geodesics in }S :& i(\boldeta_{S'}, \bc) = 0  \text{ and } \\
      
      & i(\boldeta_{S'}, \bc') \neq 0 \text{ for every closed curve } \bc' \text{ with } i(\bc, \bc') \neq 0.
\end{cases}
\]

Therefore, we have $i(\boldeta_{S'}, \ba) \neq 0$. Choose a representative of $a$ in the deck group. This gives an action of $a$ on $\mathbb{H}^2$. Let $x$ and $y$ denote the fixed points of the action of $a$ on $\mathbb{H}^2$, and let $L$ be the geodesic in $\mathbb{H}^2$ that is the axis of this action (See \cref{fig3}). Choose a point $q \in L$ such that $q$ lies in a fundamental domain of the cover corresponding to $S \setminus S'$. That is, $\pi(q) \notin S'$. The element $a$ acts on $L$ by translation along $a$. Let $(q, a \cdot q]$ denote the oriented arc along $L$ from $q$ to $a \cdot q$ including $aq$ excluding $q$. Let $G_{(q, c \cdot q]}$ be the set of all geodesics that intersect the arc $(q, a \cdot q]$. By \cite[Lemma 4.4]{MZ19:PositivelyRatioed} we know that
\begin{equation}\label{MZeq}
i(\boldeta_{S'}, \ba) = \boldeta_{S'}(G_{(q, a\cdot q]}).
\end{equation}

\begin{figure}[ht]
    \centering
    \includegraphics[width=0.9\linewidth]{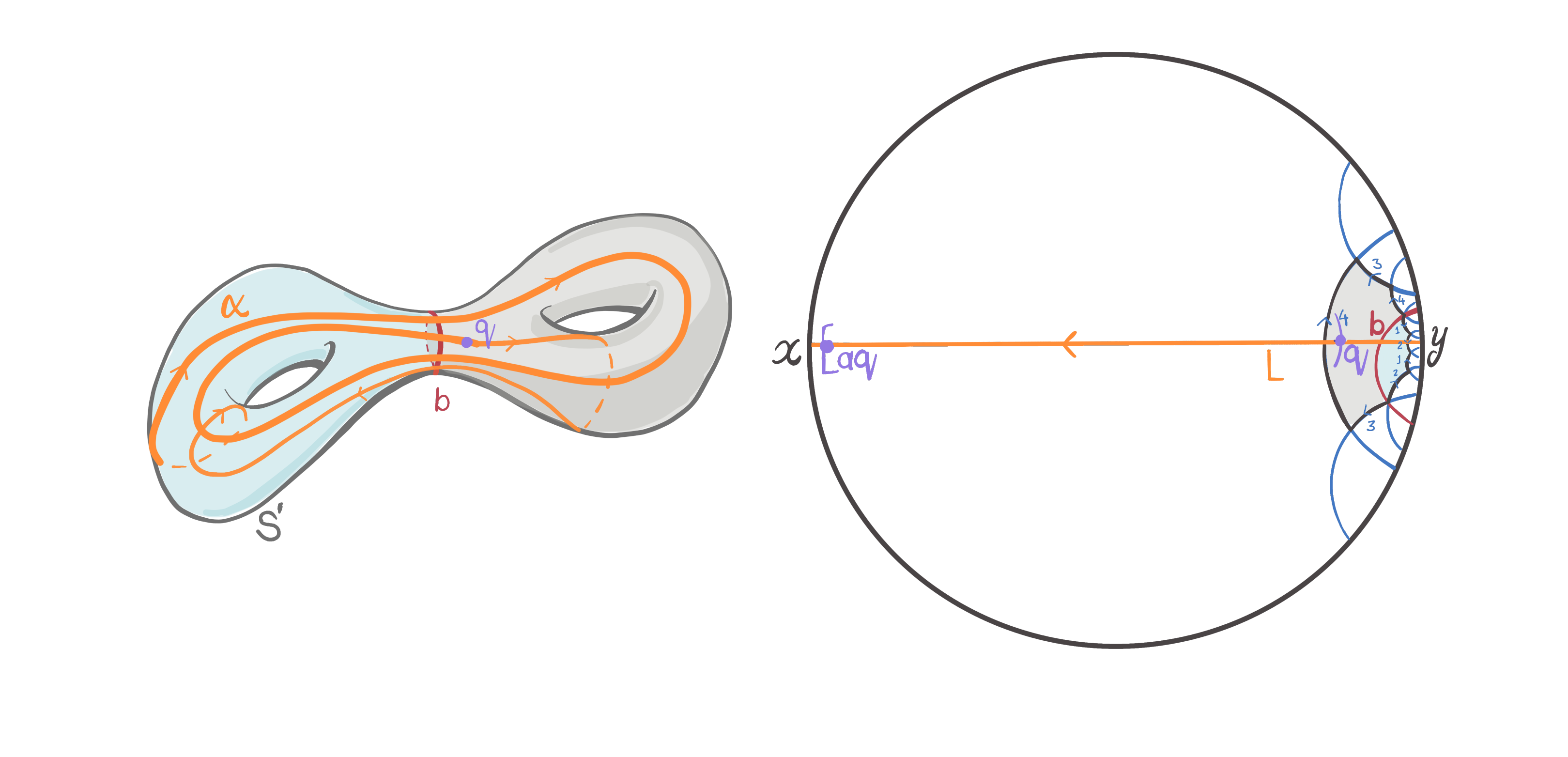}
  \caption{An example where $S$ is a genus two surface and the subsurface $S'$ filled by $\boldeta_{S'}$ is a one-holed torus. The figure shows fundamental domains for $S$ and $S'$, a lift of the boundary curve $b$, the point $q$, the arc $(q, aq]$, and the projected closed curve $\alpha$.}
  \label{fig3}
\end{figure}

Let $\alpha:[0,1] \rightarrow S$ be a closed curve in the homotopy class of $a$, obtained by projecting the arc $[q, a \cdot q]$ onto the surface $S$. We choose $\alpha(0) = \alpha(1)$ to be the image of $q$ and equip $\alpha$ with the orientation induced by the arc. 

Let $\{\alpha(t_1), \alpha(t_2), \dots, \alpha(t_n)\}$ be all the points at which $\alpha$ crosses \emph{any} boundary curve of $S'$, ordered so that $t_1 < t_2 < \dots < t_n$. Since the boundary components of $S'$ separate $S'$ from the rest of the surface, the total number of intersections, $n$, is an even number. Note that
\begin{equation}
  i(\bb,\ba) \leq n.
\end{equation}

This is because some of the intersection points $\alpha(t_i)$ may lie on a boundary component $\bb'$ of $S'$ that is distinct from $\bb$.

Consider an arc $\alpha: (t_1, t_2] \rightarrow S$ that lies entirely within the subsurface $S'$. Let $\alpha_{12}$ be the lift of the arc to $\mathbb{H}^2$ such that $\alpha_{12} \subset (q, a \cdot q]$. By our construction, the next subarc $\alpha: (t_2, t_3] \rightarrow S$ lies outside of $S'$. Since the current $\boldeta_{S'}$ is supported only on $S'$, any geodesic intersecting the subarc $\alpha_{23}$ does not contribute to $\boldeta_{S'}$. Thus,
\begin{equation} \label{0int}
   \boldeta_{S'}(G_{\alpha_{23}}) = 0.
\end{equation}

\begin{figure}[ht]
    \centering
    \includegraphics[width=\linewidth]{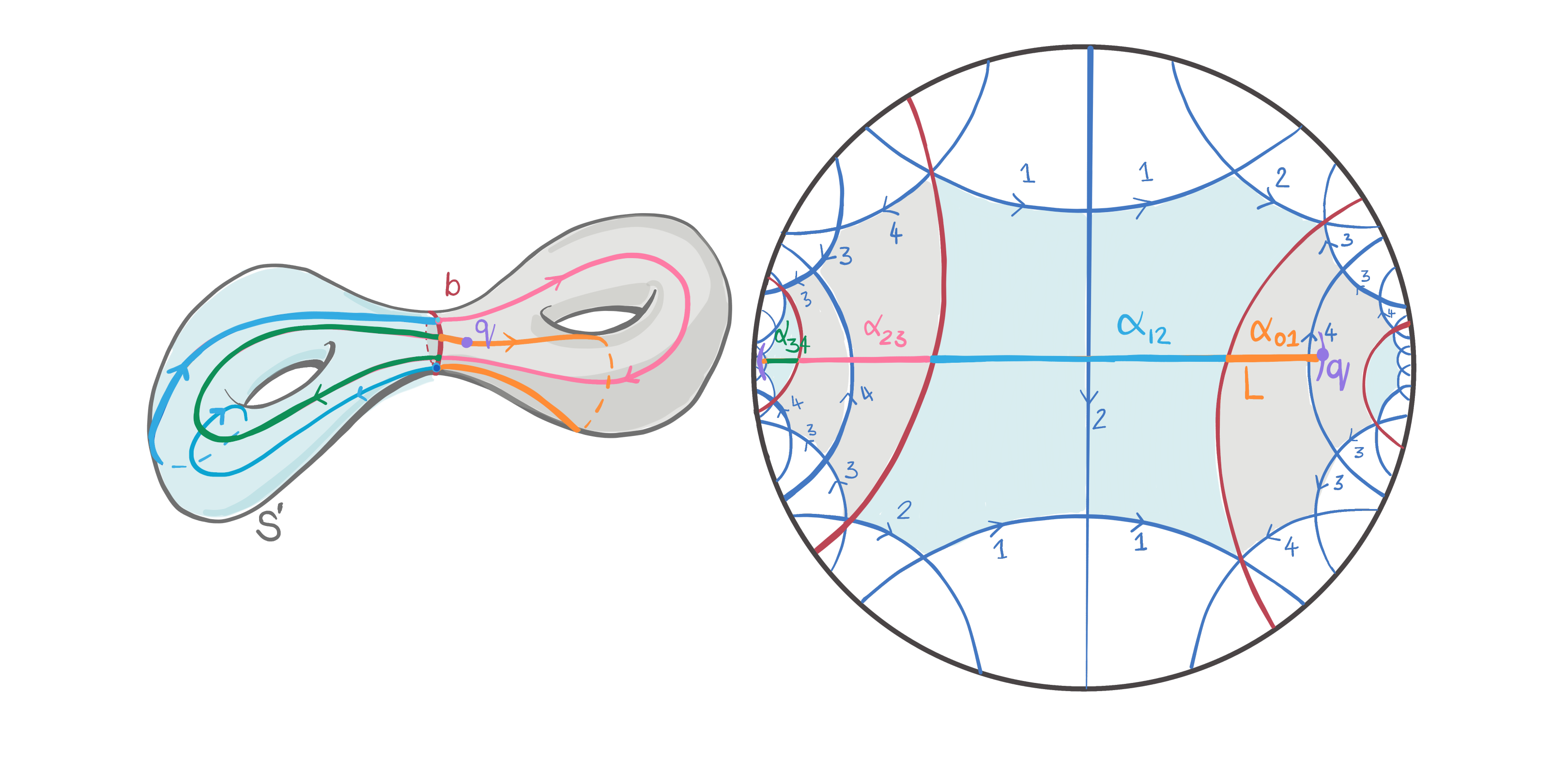}
    \caption{ The subarcs $\alpha_{(j-1)j}$ and its projection onto the surface $S$ considered for the same example as in Figure 3.}
    \label{fig4}
\end{figure}

By following a similar construction, we obtain a collection of subarcs that partition $(q, a \cdot q]$ 
\[
\{ \alpha_{01}, \alpha_{12}, \alpha_{23}, \alpha_{34}, \dots, \alpha_{(n-1)n}, \alpha_{n0}\}.
\]
For any integer $1 \leq i \leq n-1$, the projection of the arc $\alpha_{i(i+1)}$ onto the surface S is the arc $\alpha:(t_{i},t_{i+1}] \to S$. The projection of $\alpha_{01}$ is the arc $\alpha:(0,t_{1}] \to S$ , and the projection of $\alpha_{n0}$ is the arc $\alpha:(t_{n},1] \to S$.

Note that, $\alpha:(t_{i},t_{i+1}] \to S$. lies outside $S'$ when $i$ is even. Therefore,
\begin{equation} \label{0int:1}
   \boldeta_{S'}(G_{\alpha_{i(i+1)}}) = 0 \text{ when } i \text{ is even}.
\end{equation}
We also have,
\[
(q, a \cdot q] = \bigcup_{i=0}^{n-1} \alpha_{i(i+1)}.
\]

This implies,
\begin{equation}
  i(\boldeta_{S'}, a) = \boldeta_{S'}(G_{\alpha_{12}}) + \boldeta_{S'}(G_{\alpha_{34}}) + \dots + \boldeta_{S'}(G_{\alpha_{(n-1)n}}).
\end{equation}
In the latter half of the proof, we will show that each arc $\alpha_{(j-1)j}$ with $j$ even satisfies
\begin{equation}\label{eq:syst}
  \Syst_{S'}(\boldeta_{S'}) \leq 2\boldeta_{S'}(G_{\alpha_{(j-1)j}}), \quad \text{for } 2 \leq j \leq n, \text{ with $j$ even}.  
\end{equation}

Assuming \cref{eq:syst}, the remainder of the proof proceeds as follows:

\[
\frac{n}{2} \Bigg( \frac{ \Syst_{S'}(\boldeta_{S'})}{2} \Bigg) \leq \sum_{\substack{j =2 \\ j \text{ even}}}^{n} \boldeta_{S'}(G_{(\alpha_{(j-1)j}}) = i(\boldeta_{S'}, \ba).
\]
Hence for any curve $\ba$ with $i(\bb,\ba) \neq 0$ we have
\[
\frac{i(\bb, \ba)}{i(\boldeta_{S'}, \ba)} \leq \frac{n}{i(\boldeta_{S'}, \ba)} \leq\frac{4n}{n {\Syst}_{S'}(\boldeta_{S'})} = \frac{4}{{\Syst}_{S'}(\boldeta_{S'})}.
\]
Because closed curves are dense in projective currents, it follows that
\[
\sup_{c \in \mathbb{P}\mathcal{C}}\frac{i(\bb,c)}{i(\boldeta_{S'}, c)} = \sup_{c \in \cc}\frac{i(\bb,c)}{i(\boldeta_{S'}, c)} \leq \frac{4}{{\Syst}_{S'}(\boldeta_{S'}) }.
\]

Finally, we prove that the systole is a lower bound for the intersection of $\eta_{S'}$ with the arcs $\alpha_{(j-1)j}$.
We claim that
\[
\operatorname{Syst}_{S'}(\boldeta_{S'}) \leq 2 \boldeta_{S'}(G_{\alpha_{(j-1)j}}).
\]
For each arc $\alpha(t_{(j-1)}, t_j]$ on $S'$, we construct a closed curve $\gamma_{{(j-1)j}}$ such that $\gamma_{{(j-1)j}}$ lies in the interior of $S'$ and satisfies 
\[
i(\boldeta_{S'}, \bgamma_{(j-i)j}) \leq 2 \boldeta_{S'}(G_{\alpha_{(j-1)j}}).
\]
In order to construct this curve we consider the following two cases:\\

\textbf{Case I.} $\alpha(t_{j-1})$ and $\alpha(t_{j})$ lie on the same boundary curve $\bb'$ of $S'$.

Note that, $\bb' \in \mathcal{E}_{\eta_{S'}}$, and $\bb'$ does not necessarily have to be same as $\bb$. 

The arc $\alpha(t_{j-1}, t_j]$ is \emph{essential} in $S'$; that is, it is not homotopic (relative to its endpoints) to an arc entirely contained in $\partial S'$. This follows from the fact that $\alpha(t)$ is a geodesic representative that realizes the minimal intersection with $\eta_{S'}$. We define $\gamma_{(j-1)j} : [0,1] \to S'$ to be a closed curve obtained by following the arc $\alpha[t_{j-1}, t_j)$ and then concatenating it with one of the arcs $\beta_1$, $\beta_2$ or $\beta_2$ illustrated in \cref{fig5}. Since $\alpha[t_{j-1}, t_j)$ is essential, we can always choose one of $\beta_1$, $\beta_2$ or $\beta_3$ such that the resulting concatenation is not homotopic to a point or to a boundary component of $S'$. We choose that $\beta_i$ to close up the curve $\gamma_{ij}$.

\begin{figure}[ht]
    \centering
    \includegraphics[width=0.8\linewidth]{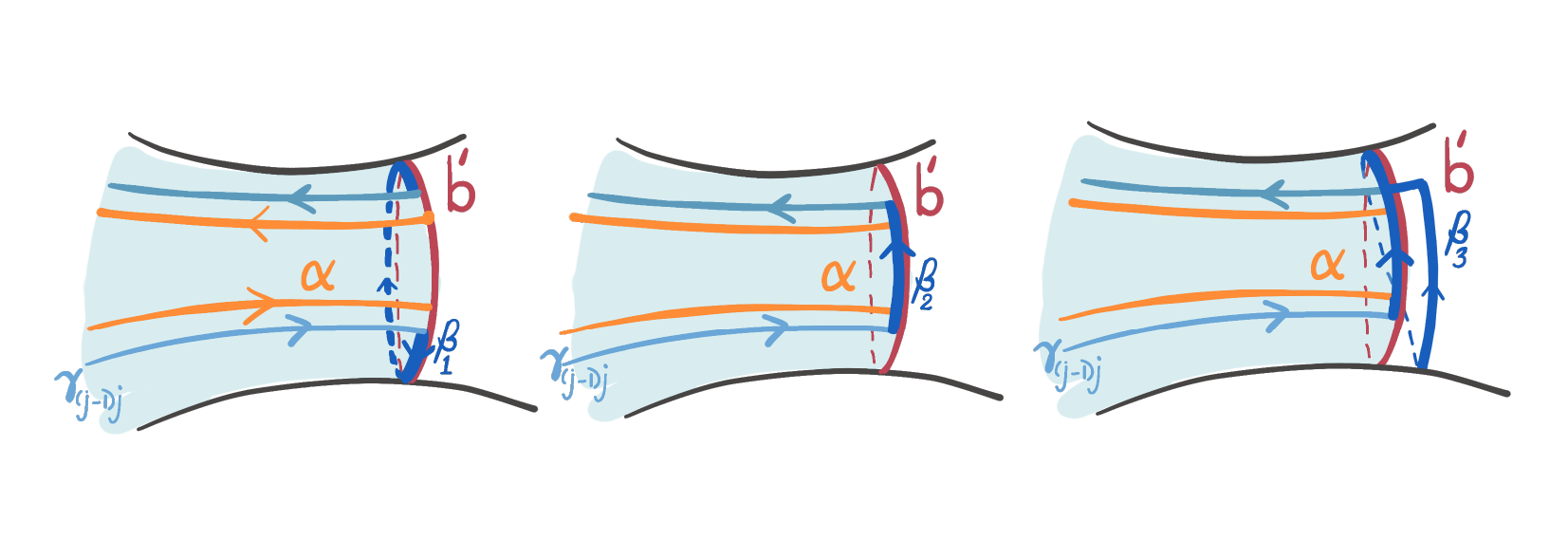}
    \caption{ The subarcs $\beta_1, \beta_2$ and $\beta_3$ that can be used to construct the closed curve $\gamma_{(j-1)j}$ in case I}.
    \label{fig5}
\end{figure}

\textbf{Case II.} $\alpha(t_{j-1})$ and $\alpha(t_{j})$ lie on distinct boundary curves $\bb'$ and $\bb''$ of $S'$.

Suppose that $\alpha(t_{j-1})$ lies on the boundary component $\bb'$ and $\alpha(t_j)$ lies on a different boundary component $\bb''$. In this case, we define $\gamma_{(j-1)j}$ as the closed curve constructed through the following steps:
\begin{enumerate}
    \item  Follow the arc $(\alpha(t_{j-1}), \alpha(t_j)]$ on $S'$. This will give us an arc start at $\bb'$ and ending at $\bb''$.
    \item  Follow the boundary curve $\bb''$ until it meets the arc constructed in the previous step again.
    \item Follow the arc $(\alpha(t_{j-1}), \alpha(t_j)]$ in the reverse. This gives us an arc from $b''$ to $b' $ (See \cref{case 2}).
    \item Concatenate  with one the arcs $\beta_1$, or $\beta_2$ given in \cref{case 2} so that the resulting closed curve is not homotopic to a point or to a boundary component of $S'$. For similar reasons as discussed in case I, we will always we able to find one such $\beta_{i}$.
\end{enumerate}

\begin{figure}[ht]
    \centering
    \includegraphics[width=\linewidth]{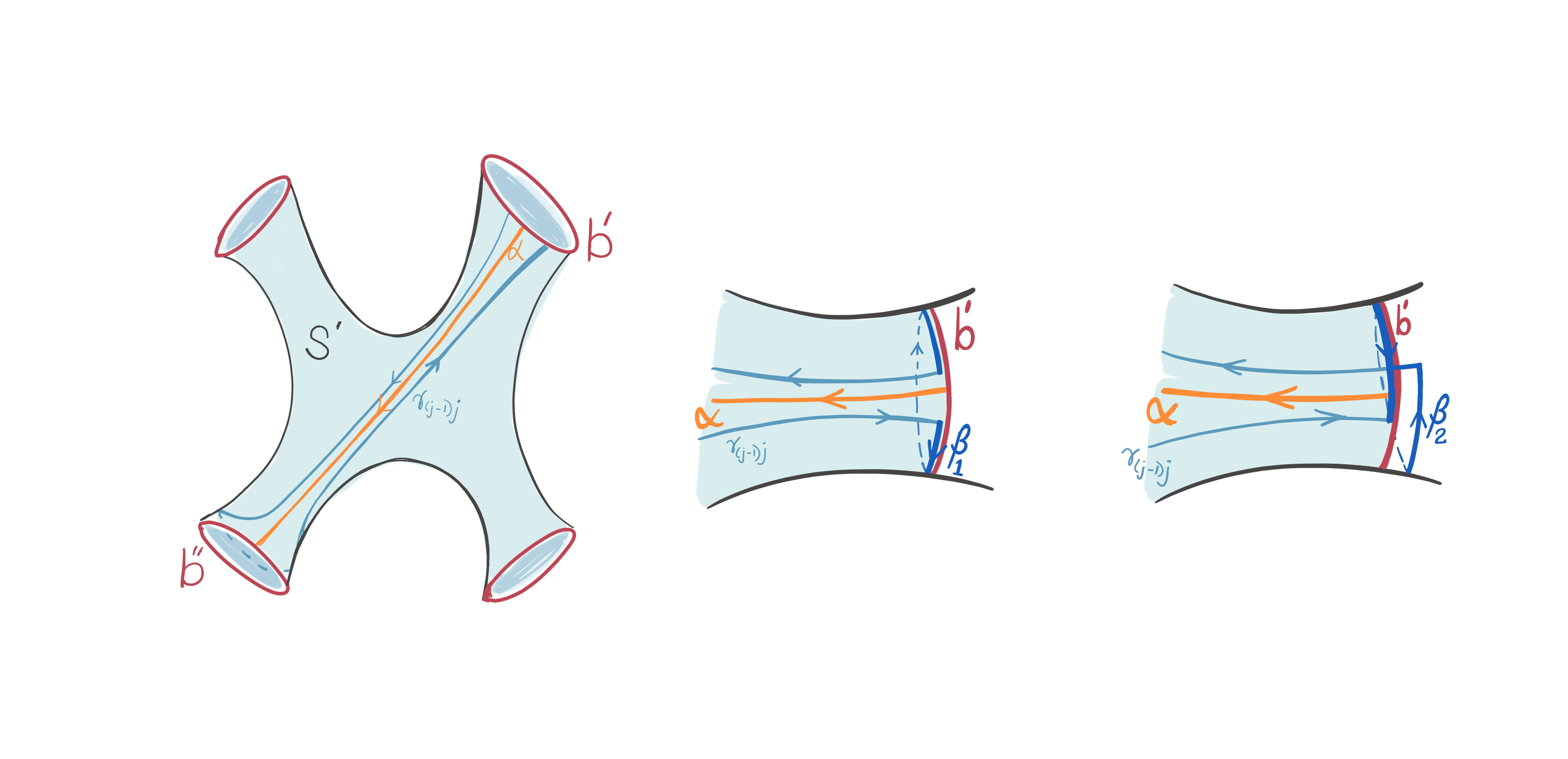}
    \caption{An example where $S'$ is a four-holed sphere. The figure illustrates the construction of the arc $\gamma_{(j-1)j}$ before it is completed by concatenating one of the arcs $\beta_1$, or $\beta_2$ shown on the right. In this specific example, the arc used to complete the construction is $\beta_2$.}
    \label{case 2}
\end{figure}

For the closed curve $\gamma_{(j-1)j}$ constructed in Case I, we consider its lift to $\mathbb{H}^2$ that follows the geodesic $\alpha_{(j-1)j}$. This lifted curve in $\mathbb{H}^2$ first traces along $\alpha_{(j-1)j}$, then continues along the lift of the subarc $\beta_i$ used to complete the construction of $\gamma_{(j-1)j}$. We denote this lift of $\beta_i$ by $\tilde{\beta_{i}}$. The complete lift of $\gamma_{(j-1)j}$ to $\mathbb{H}^2$ travels from one point on the boundary of $\mathbb{H}^2$ to another repeating this pattern along different lifts of $\alpha(t_{j-i}, t_{j}]$ and $\beta_{i}$ \cite[Lemma 4.4]{MZ19:PositivelyRatioed}.

The lift of $\gamma_{(j-1)j}$ we constructed here has the same endpoints at infinity as its geodesic representative. From $\pi_{1}$-invariance of geodesic currents one can conclude that:
\begin{equation}\label{eq1}
    i(\boldeta_{S'}, \bgamma_{{(j-1)}j}) \leq  \boldeta_{S'}(G_{\alpha_{(j-1)j}}) + \boldeta_{S'}(G_{\tilde{\beta_i}}).
\end{equation}

One can run a similar argument for the closed curve $\bgamma_{{(j-1)}j}$ constructed in Case II to get:
\begin{equation}\label{eq2}
    i(\boldeta_{S'}, \bgamma_{{(j-1)}j}) \leq  2\boldeta_{S'} (G_{\alpha_{(j-1)j}}) + i(\boldeta_{S'}, \bb'' ) + \boldeta_{S'}(G_{\tilde{\beta_i}}).
\end{equation}
Note that, the lifts of the arc $\beta_{i}$ lies along the lifts of the boundary curve $\bb'$. Since $i(\boldeta_{S'}, \bb') = 0$ we have 
\begin{equation}\label{eq5}
    \boldeta_{S'}(G_{\tilde{\beta_i}}) = 0.
\end{equation}

From \cref{eq1}, \cref{eq2} and \cref{eq5}, regardless of whether we consider Case I or II, we always obtain
\[
i(\boldeta_{S'}, \bgamma_{{(j-1)}j}) \leq 2\boldeta_{S'}(G_{\alpha_{(j-1)j}})
\]
and therefore, 
\[
\Syst_{S'}(\boldeta_{S'}) \leq 2\boldeta_{S'}(G_{\alpha_{(j-1)j}}).
\]
This completes the proof of the claim and thereby the proof of the Lemma.
\end{proof}

\begin{lemma}\label{boundarycurve}
Let $\eta_{S'}$ be a non-filling projective current that fills a connected subsurface $S' \subset S$.  Let $b$ be a simple closed curve in $\mathcal{E}_{\eta_{S'}}$ that forms one of the boundary components of $S'$.  
Then  
\[
H(\eta_{S'}, b) < \infty.
\]
\end{lemma}

\begin{proof}
By \cref{finitedetourcost1.2}, to prove $H(\eta_{S'}, b) < \infty$ it is enough to show that  
\begin{equation}\label{eq:lemma8.6}
\sup_{c \in \mathbb{P}\mathcal{C}} \frac{i(b,c)}{i(\eta_{S'}, c)} < \infty.    
\end{equation}

Let $\bb$ be the geodesic current corresponding to the closed curve $b$, and let $\boldeta_{S'}$ be a geodesic current in the projective class of $\eta_{S'}$.  
Assign the geodesic current $\bb$ a positive weight $k > 0$.  
By \cref{lemma:systolineqality} we have  
\begin{equation}\label{eq:systinequality}
   \sup_{c \in \mathbb{P}\mathcal{C}} \frac{i(k\bb,c)}{i(\boldeta_{S'}, c)}
   = k \left( \sup_{c \in \mathbb{P}\mathcal{C}} \frac{i(\bb,c)}{i(\boldeta_{S'}, c)} \right) 
   \leq \frac{4k}{\Syst_{S'}(\boldeta_{S'})}.
\end{equation}
Since $\boldeta_{S'}$ fills $S'$, $\Syst_{S'}(\boldeta_{S'}) > 0$.  
Consequently, the quantity in \cref{eq:systinequality} is finite.  
As \cref{eq:systinequality} does not depend on the choice of representative $\boldeta_{S'}$ or the weight $k$, we obtain \cref{eq:lemma8.6}, which proves the lemma.
\end{proof}

\begin{lemma}\label{deltafilling}
    Let $\eta_{S'}$ be a non-filling geodesic current that fills a subsurface $S' \subset S$. Let $b$ be a simple closed curve that is a boundary component of $S'$. Then $\delta(\eta_{S'}, \eta_{S'} + b) < \infty$.
\end{lemma}
\begin{proof}
    First, note that
    \[
    \sup_{c \in \mathbb{P}\mathcal{C}}\frac{i(\eta_{S'}, c)}{i(\eta_{S'} + b, c)} \leq 1.
    \]
So, by \cref{finitedetourcost1.2}, we have $H(\eta_{S'}+b, \eta_{S'}) < \infty$. To show finiteness of $H(\eta_{S'}+b, \eta_{S'})$ we consider 
\begin{align*}
   \sup_{c \in \mathbb{P}\mathcal{C}}\frac{i(\eta_{S'} +b, c)}{i(\eta_{S'}, c)} =   \sup_{c \in \mathbb{P}\mathcal{C}}1 + \frac{i(b,c)}{i(\eta_{S'}, c)}.
\end{align*}
By \cref{boundarycurve} we know that $\displaystyle{\sup_{c\in \mathbb{P}\mathcal{C}}\frac{i(b,c)}{i(\eta_{S'},c)}}$ is finite. This proves that the above supremum is finite and therefore $\delta(\eta_{S'}, \xi) < \infty$. 
\end{proof}

\begin{proposition}\label{finitedelta}
A non-filling projective current $\eta$ is a minimal measured lamination if and only if every non-filling projective current $\xi$ with $\delta(\eta, \xi) < \infty$ satisfies $\xi = \eta$.
\end{proposition}

\begin{proof}
    If $\eta$ is a minimal measured lamination then by \cref{finitedetorcost2} and \cref{minimal ML} we know that the only projective non-filling current finite $\delta$-distance away from $\eta$ is $\eta$ itself. 

    Let's say $\eta$ has a non-trivial filling component $\eta_{S'}$ that fills a subsurface $S'$ of $S$. Note that $S'$ need not be connected. Let $b \in \mathcal{E}_{\eta_{S'}}$. Here $b$ is a boundary curve of one of the connected components of $S'$. We show that $\delta(\eta + b, \eta)$ is finite. Let $\boldeta$ be a representative of $\eta$ in the space of currents. We have 
    \[
    \boldeta = \boldeta_{S'} + \blambda_{\boldeta}.
    \]
     Note, that the following supremum is always less than one as the denominator is always bigger than the numerator
    \[
    \sup_{c \in \mathbb{P}\mathcal{C}}\frac{i(\boldeta, c)}{i(\boldeta + \bb, c)}
    \]
    and by \cref{finitedetourcost1.2} it follows that $H(\eta+ b, \eta) < \infty$. Now,
    \[
    \sup_{c \in \mathbb{P}\mathcal{C}}\frac{i(\boldeta + \bb, c)}{i(\boldeta, c)} =  \sup_{c \in \mathbb{P}\mathcal{C}}\frac{i(\boldeta_{S'} + \bb + \blambda_{\boldeta}, c)}{i(\boldeta_{S'} + \blambda_{\boldeta}, c)} \leq \sup_{c \in \mathbb{P}\mathcal{C}}\frac{i(\boldeta_{S'} + \bb, c)}{i(\boldeta_{S'} + \blambda_{\boldeta}, c)} + \sup_{c \in \mathbb{P}\mathcal{C}}\frac{ i(\blambda_{\boldeta}, c)}{i(\boldeta_{S'} + \blambda_{\boldeta}, c)}.
    \]
    From \cref{boundarycurve} it follows the first summand in the above equation is finite
\[
\sup_{c \in \mathbb{P}\mathcal{C}}\frac{i(\boldeta_{S'} + \bb, c)}{i(\boldeta_{S'} + \blambda_{\boldeta}, c)} \leq \sup_{c \in \mathbb{P}\mathcal{C}}\frac{i(\boldeta_{S'} + \bb, c)}{i(\boldeta_{S'} , c)} < \infty.
\]
And for the second summand we get 
\[
\sup_{c \in \mathbb{P}\mathcal{C}}\frac{ i(\blambda_{\boldeta}, c)}{i(\boldeta_{S'} + \blambda_{\boldeta}, c)} \leq 1.
\]
Again, by \cref{finitedetourcost1.2}, we have $H(\eta, \eta + b) < \infty$. Therefore, no projective current with a filling component can satisfy the condition in the proposition. Now, suppose that $\eta $ is a non-minimal projective measured lamination. That is, we can write a representative $\boldeta$ of $\eta$ as 
\[
\displaystyle{\boldeta = \sum_{i = 1}^{n} k_{i} \boldeta_{i}}
\]
where each $k_{i}$ is a positive real number, $\boldeta_{i}$'s are minimal measured laminations, and $n$ is at least 2. Then, by \cref{deltaML}, the projective class $\xi$ of $\bxi$ defined as 
\[
\bxi = \sum_{i =1}^{n} k_{i}' \boldeta_{i}, \text{ with } k_{i}' >0 \text{ and } k_{i} \neq k_{i}.
\]
satisfies $\delta(\eta, \xi) < \infty$. This completes the proof of the proposition.
\end{proof}
\section{Isometric Rigidity}

Let $S_g$ and $S_{g'}$ be closed, orientable, smooth surfaces of genus $g$ and $g'$, respectively. Without loss of generality, assume $g' \leq g$. Let $f: (\mathbb{P}\mathcal{C}_{fill}(S_g), d) \to (\mathbb{P}\mathcal{C}_{fill}(S_{g'}), d)$ be an isometry, then $f$ induces a homeomorphism $f^*: \mathbb{P}\mathcal{C}_{S_g}(\infty) \to \mathbb{P}\mathcal{C}_{S_{g'}}(\infty)$ between the corresponding horoboundaries, as described in \cref{f*}.

In this section, we show that the induced homeomorphism $f^*$ maps simple closed curves to minimal projective  measured laminations (\cref{f(scc)1}), and also preserves disjointness (\cref{disjoint}). However, this is only possible when $g = g'$ (\cref{isometricrigidity}).

\begin{lemma}\label{f(scc)1}
The homeomorphism $f^*: \mathbb{P}\mathcal{C}_{S_{g}}(\infty) \rightarrow \mathbb{P}\mathcal{C}_{S_{g'}}(\infty)$ maps simple closed curves to  minimal projective  measured laminations.
\end{lemma}
\begin{proof}
Let $\lambda$ be a simple closed curve on $S_{g}$. Consider $f^{*}(\lambda) \in \mathbb{P}\mathcal{C}_{S_{g'}}(\infty)$. Suppose $\eta$ is a non-filling projective current in $\mathbb{P}\mathcal{C}(S_{g'})$ such that $\delta(f^*(\lambda), \eta) < \infty$. Since $f^{*}$ is an isometry with respect to $\delta$, we get $\delta(\lambda, (f^{*})^{-1}(\eta)) < \infty$. By \cref{finitedelta} the only projective current with finite symmetrized detour cost to $\lambda$ is $\lambda$ itself.  Hence, $(f^{*})^{-1}(\eta) = \lambda$, which implies $\eta = f^*(\lambda)$. Thus, any geodesic current $\eta$ satisfying $\delta(f^*(\lambda), \eta) < \infty$ must be $f^*(\lambda)$. By \cref{finitedelta}, $f^*(\lambda)$ is a minimal projective measured lamination.
\end{proof}

\begin{lemma}\label{disjoint} Let $\eta_1$ and $\eta_2$ be disjoint simple closed curves. Then the images of $\eta_1$ and $\eta_2$ under the map $f^*: \mathbb{P}\mathcal{C}_{S_{g}}(\infty) \rightarrow \mathbb{P}\mathcal{C}_{S_{g'}}(\infty)$ are disjoint:
\[
i(f^{*}(\eta_1), f^*(\eta_2)) = 0.
\]
\end{lemma}

\begin{proof} 

Since $\eta_1$ and $\eta_2$ are disjoint simple closed curves, $\eta_1 + \eta_2$ is a measured lamination. Using the decomposition of geodesic currents and \cref{finiteDC}, one can conclude that any non-filling projective current $\xi$ satisfying
\[
H(\eta_1 + \eta_2, \xi) < \infty
\]
must be of the form $\xi = k_1 \eta_1 + k_2 \eta_2$, where $k_1$ and $k_2$ are non-negative real numbers. 

We want to show that $i(f^*(\eta_1), f^*(\eta_2)) = 0$. This is equivalent to proving that $f^*(\eta_1) + f^*(\eta_2)$ is a measured lamination. We prove this in the following steps:\\

\underline{\textbf{Step 1:} We claim that $f^*(\eta_1) +  f^*(\eta_2)$  is a non-filling projective current.}\\
\vspace{-0.2cm}

By \cref{deltaML} we know that
\[
H(\eta_1 + \eta_2, \eta_i) < \infty \text{ for } i = 1 \text{ and } 2.
\]
Since $f^*$ preserves finiteness of detour cost \cite[Proposition 5.3]{CW2015}, we get:
\begin{equation}\label{eqdisjointness}
  H(f^*(\eta_1 + \eta_2), f^*(\eta_i) ) < \infty \text{ for } i = 1 \text{ and } 2. 
\end{equation}
Observe that  $f^*(\eta_1 + \eta_2)$ is a non-filling projective current, as it is the image of a non-filling projective current under $f^*$. This means that we can find $c \in \mathbb{P}\mathcal{C}$ such that 
\[
i(f^*(\eta_1 + \eta_2), c) = 0.
\]
By \cref{finiteDC}, we have $i(f^*(\eta_1), c) = 0$ and $i(f^*(\eta_2), c) = 0$. Therefore, 
\[
i(f^*(\eta_1) + f^*(\eta_2), c) = 0.
\]
Hence, $f^*(\eta_1) + f^*(\eta_2)$ is non-filling.

\underline{\textbf{Step 2:} We claim that $f^*(\eta_1) +  f^*(\eta_2) = f^*(k_1 \eta_1 + k_2 \eta_2)$ for some non-zero constants $\eta_1$ and $\eta_2$.}\\
\vspace{-0.1cm}

By \cref{finitedetourcost1.2}, \cref{eqdisjointness} is equivalent to the finiteness of the following supremum:
\[
\sup_{c \in \mathbb{P}\mathcal{C}} \frac{i(f^*(\eta_i), c)}{i(f^*(\eta_1 + \eta_2), c)} < \infty \quad \text{for } i = 1 \text{ and } 2.
\]
This implies
\[
\sup_{c \in \mathbb{P}\mathcal{C}} \frac{i(f^*(\eta_1) + f^*(\eta_2), c)}{i(f^*(\eta_1 + \eta_2), c)} \leq 
\sup_{c \in \mathbb{P}\mathcal{C}} \frac{i(f^*(\eta_1), c)}{i(f^*(\eta_1 + \eta_2), c)} +
\sup_{c \in \mathbb{P}\mathcal{C}} \frac{i(f^*(\eta_2), c)}{i(f^*(\eta_1 + \eta_2), c)} < \infty.
\]
By \cref{finitedetourcost1.2}, this is equivalent to
\[
H(f^*(\eta_1 + \eta_2), f^*(\eta_1) + f^*(\eta_2)) < \infty.
\]
Since $f$ is an isometry, it is invertible. Moreover, $(f^{*})^{-1} = (f^{-1})^*$ also preserves finiteness of the detour cost. Therefore,
\[
H(\eta_1 + \eta_2, (f^{-1})^* \big(f^*(\eta_1) + f^{*}(\eta_2)\big)) < \infty.
\]
As discussed earlier in the proof, this is only possible if 
\begin{align*}
   (f^{-1})^* \big(f^*(\eta_1) + f^{*}(\eta_2)\big) = k_{1}\eta_{1} + k_{2}\eta_{2}  
\end{align*}
for non negative constants $k_1$ and $k_2$. Hence, 
\[
f^*(\eta_1) + f^{*}(\eta_2) = f^*(k_{1}\eta_{1} + k_2\eta_{2}).
\]

\underline{\textbf{Step 3:} We claim that $f^*(\eta_1) +  f^*(\eta_2)$ does not have a filling component.}\\
\vspace{-0.1cm}

Assume $i(f^*(\eta_1), f^*(\eta_2)) \neq 0$, or equivalently $f^{*}(\eta_1) + f^*(\eta_2)$ is not a measured lamination. This means that any geodesic current representative of $f^{*}(\eta_1) + f^*(\eta_2)$ has a filling component that fills a proper subspace $S'$ of $S$. Let $b \in \mathcal{E}_{f^{*}(\eta_1) + f^*(\eta_2)}$ be a boundary curve of $S'$. Note that, 
\[
b \neq f^*(\eta_i) \text{ for }i = 1 \text{ or } 2.
\]
This is true because by our construction we know that $i(f^*(\eta_1) + f^*(\eta_2) , b) = 0$, and hence both $i(f^*(\eta_1), b) = 0$ and $i(f^*(\eta_2), b) = 0$

We have
\[
H\big( k_1\eta_1 + k_2\eta_2, (f^*)^{-1}(b) \big) = H(f^*(k_1\eta_1 + k_2\eta_2), b) = H(f^*(\eta_1) + f^{*}(\eta_2), b) < \infty.
\]
where the first equality uses $f^*$-invariance of $H$, the second one uses Step 2, and the last inequality uses \cref{boundarycurve}
Here, $k_1\eta_1 + k_2\eta_2$ is a weighted multicurve consisting of disjoint simple closed curves. Using an argument similar to the one at the beginning of this proof, we conclude that the above detour cost can be finite only if
\[
(f^*)^{-1}(b) = \alpha_1 \eta_1 + \alpha_2 \eta_2,
\]
where $\alpha_1$ and $\alpha_2$ are non-negative real numbers. Since $b$ is a simple closed curve $(f^*)^{-1}(b)$ is a minimal measured lamination. Therefore, at least one of $\alpha_1$ or $\alpha_2$ must be zero. Without loss of generality, we assume $\alpha_2 = 0$. Since projective currents are considered up to scaling, we get
\begin{align*}
    (f^*)^{-1}(b) &= \eta_1\\
    f^*(\eta_1) &= b.\\
\end{align*}
However, we showed that $b$ is distinct from both $f^*(\eta_1)$ and $f^*(\eta_2)$, leading to a contradiction. Therefore, $f^*(\eta_1) + f^*(\eta_2)$ must be a measured lamination, completing the proof.
\end{proof}

\begin{reptheorem}{isometricrigidity}
Let $S_{g}$ and $S_{g'}$ be closed orientable surfaces of genus at least two, with $g \neq g'$. Then the spaces of projective filling currents $\mathbb{P}\mathcal{C}_{fill}(S_{g})$ and $\mathbb{P}\mathcal{C}_{fill}(S_{g'})$, endowed with the extended Thurston metric $d$, are not isometric.
\end{reptheorem}

\begin{proof}
Without loss of generality, assume $g' < g$. Suppose, for contradiction, that there exists an isometry 
\[
f: \mathbb{P}\mathcal{C}_{fill}(S_g) \to \mathbb{P}\mathcal{C}_{fill}(S_{g'})
\]
with respect to the extended Thurston metric. By \cref{f(scc)1}, the induced map on the horoboundaries
\[
f^* : \mathbb{P}\mathcal{C}_{S_g}(\infty) \to \mathbb{P}\mathcal{C}_{S_{g'}}(\infty)
\]
sends simple closed curves to minimal projective measured laminations. Let $\{\gamma_1, \gamma_2, \dots, \gamma_{3g - 3}\}$ be a pants decomposition of $S_g$. Then, by \cref{disjoint}, the set $\{f^*(\gamma_1), f^*(\gamma_2), \dots, f^*(\gamma_{3g - 3})\}$ consists of pairwise disjoint projective measured laminations on $S_{g'}$. However, this contradicts the fact that a surface of genus $g'$ admits at most $3g' - 3$ disjoint measured laminations, and by assumption $g' < g$. 
\end{proof}

Although Theorem~\ref{isometricrigidity} may at first seem natural, it is in fact non-trivial, since the space of geodesic currents is infinite-dimensional.
For instance, the topology of $S$ alone does not determine the topological type of $\mathbb{P}\mathcal{C}(S)$: indeed, $\mathbb{P}\mathcal{C}(S_g)$ is homeomorphic to $\mathbb{P}\mathcal{C}(S_{g'})$ for all genera $g,g'$, a fact that follows from results in Choquet theory. We now recall the necessary preliminaries to show this.

For the following preliminaries, we direct the reader to the references~\cite{Dow91:Choquet,phelps2001choquet}.
Let $K$ be a convex subset of a topological vector space. 
A point $x \in K$ is called an \emph{extreme point} of $K$ if whenever 
\[
x = t y + (1-t) z \quad \text{with } y,z \in K,\ 0 < t < 1,
\]
we have $x = y = z$. 
The set of all extreme points of $K$ is denoted $\mathrm{ext}(K)$.

If $K$ is moreover compact, 
by the Krein--Milman theorem, every $x \in K$ can be represented as the barycenter of a probability measure $\mu$ supported on the set of extreme points $\mathrm{ext}(K)$, i.e.
\[
x = \int_{\mathrm{ext}(K)} y \, d\mu(y).
\]
If for each $x \in K$ this representing measure $\mu$ is \emph{unique}, then $K$ is called a \emph{Choquet simplex}.
Examples of Choquet simplices are finite-dimensional simplices in $\mathbb{R}^n$, and infinite dimensional spaces such as probability measures on a compact metric space
 (with weak* topology).

The \emph{Poulsen simplex} is the unique (up to affine homeomorphism) separable Choquet simplex whose set of extreme points is dense in the simplex~\cite{poulsen1961simplex,LOS78:Poulsen}.

\begin{proposition}
All spaces of projective geodesic currents are homeomorphic.
\label{prop:geodesiccurrents_homeo}
\end{proposition}
\begin{proof}
By the above, it suffices to show that for any closed hyperbolic surface $S$, $\PC(S)$ is homeomorphic to a Poulsen simplex.
The space of geodesic currents $\mathcal{C}$ is affinely homeomorphic 
to the space $\mathcal{M}_{\text{flip},\text{inv}}(T^1S)$ of flip invariant and geodesic flow invariant
locally finite Borel measures on $T^1 S$~\cite{Bon86,Bon88, ES22:GeodesicCount}, which is a closed subspace of $\mathcal{M}_{\text{inv}}(T^1S)$, the space of geodesic flow invariant measures on $T^1 S$.
By normalizing by mass (which is a homeomorphism in the weak$^*$ since $T^1 S$ is compact), $\mathbb{P}\mathcal{M}_{\text{inv}}(T^1S)$ is homeomorphic to the space $\mathcal{P}_{\text{inv}}(T^1S)$ of geodesic flow invariant Borel probability measures on $T^1 S$, which is known to be a Poulsen simplex.
Similarly, $\PC(S)$ is homeomorphic to the projectivized space
$\mathbb{P}\mathcal{M}_{\text{flip},\text{inv}}(T^1S)$, which is homeomorphic to $\mathcal{P}_{\text{flip},\text{inv}}(T^1S)$.
Since $\mathcal{P}_{\text{flip},\text{inv}}(T^1S)$ is a closed and convex subspace of $\mathcal{P}_{\text{inv}}(T^1S)$, it is a convex compact space: this is because the flip map $\sigma_* \colon \mathcal{M}(T^1 S) \to \mathcal{M}(T^1S)$ induced by the orientation reversing homeomorphism is an affine homeomorphism, so it preserves convexity. Hence $\mathcal{P}_{\text{flip},\text{inv}}(T^1S)$ is a Choquet simplex. Since a geodesic current supported on a closed curve is extremal, and the space of such measures is dense in $\mathcal{P}_{\text{flip},\text{inv}}(T^1S)$~\cite[Proposition~4.4]{Bon86}, it follows that $\mathrm{ext}(\mathcal{P}_{\text{flip},\text{inv}}(T^1S))$ is dense, hence $\mathcal{P}_{\text{flip},\text{inv}}(T^1S)$ is a Poulsen simplex, and, since $\PC(S)$ is homeomorphic to it, the result follows. 
\end{proof}

\section{Examples of induced asymmetric metrics on geometric structures}
\label{sec:examples}
As mentioned in the introduction, the space of projective filling geodesic currents $\PC_{fill}$ contains many geometric structures of interest. As we sketched in the introduction and is explained in detail in Appendix~\ref{sec:asymmetric}, the space of projective filling geodesic currents on a surface embeds in the broader space of metric structures $\mathcal{D}_{\pi_1(S)}$ 
~\cite{LucDid,StephenEduardo}, which further embeds into the space of proper metric potentials $\mathcal{H}^{++}_{\pi_1(S)}$, in the sense of~\cite{CRS24:Joint}. More generally, these two spaces can be defined for any non-elementary hyperbolic group $\Gamma$. The extended Thurston metric $d$ is then defined on $\mathcal{H}^{++}_\Gamma$, and when $\Gamma=\pi_1(S)$ it restricts to the metric $d$ on $\PC_{fill}$ used in Section~\ref{sec:compactification} (see Proposition~\ref{prop:asymmetric_defined}).

We now describe how to induce an asymmetric metric on collections of geometric structures.

Let $\mathcal{M}$ be a \emph{family of geometric structures} associated to $\Gamma$, i.e., a class where for each $m \in \mathcal{M}$ there is a \emph{length function} $\ell_m \colon [\Gamma] \to \mathbb{R}$, satisfying $\ell_m(g)=\ell_m(g^{-1})$ for every $g \in \Gamma$.
Typically, $\mathcal{M}$ will be 
defined as a set of equivalence classes of geometric objects, where the equivalence will depend on the context, e.g., diffeomorphism isotopic to the identity, conjugacy, etc. 

To apply this framework, two conditions are required:
\begin{itemize}
\item \textbf{Requirement 1: realizable.} Every $m \in \mathcal{M}$ can be \emph{realized} as a proper hyperbolic potential $\Gamma$, i.e., there exists a map assigning to each $m \in \mathcal{M}$ an element $\psi_m \in H_\Gamma^{++}$ so that $\ell_m(g)=\ell_{\psi_m}(g)$ for every $g \in \Gamma$.
\item \textbf{Requirement 2: normalized marked length spectrum.} The assignment $m \mapsto \psi_m$ defines an injective map into $\mathcal{H}^{++}_\Gamma$. Note that $\mathcal{H}^{++}_\Gamma$ is a space of rough similarity classes, where two proper hyperbolic metric potentials are equivalent if they differ by a scaling up to additive error (see Appendix~\ref{sec:asymmetric} for definitions). We say that $\mathcal{M}$ is \emph{closed under scaling}, i.e., there is an action $\mathbb{R}^+$ on $\mathcal{M}$ so that $\lambda \in \mathbb{R}^+$ and $\lambda \cdot m \in \mathcal{M}$ whenever $\lambda \in \mathbb{R}^+$, $m \in \mathcal{M}$. In this case, if $\mathcal{M}$ satisfies Requirement 1  and the stable length scales homogeneously $\ell_{\psi_{\lambda \cdot m}} = \lambda \ell_{\psi_m}$, then, in order for $\mathcal{M}$ to satisfy Requirement 2, one needs to restrict to a class of metrics where a certain homogeneous invariant is normalized (such as entropy 1, area 1, etc.).
\end{itemize}

\subsection*{Examples realizable as geodesic currents}

In many interesting cases, we will be able to satisfy Requirement 1 with $\psi_m \in \mathcal{D}_\Gamma$. 
In fact, if $\Gamma=\pi_1(S)$, sometimes Requirement 1 can be realized by embedding $m$ into $\PC_{fill}$ via a geodesic current $\mu_m$. The assignment $m \mapsto \mu_m$ gives then an embedding into projective filling currents, i.e., if $[\mu_m] = [\mu_{m'}]$ then $m=m'$. From here, one can produce the pseudo-metric $\rho_{\mu_m}$ (as in described in Appendix~\ref{subsec:from_currents_to_pseudometrics}).

Many examples of \emph{marked} metric structures (where the relation of marked is via diffeomorphism isotopic to the identity) give rise to families $\mathcal{M}$ on $S$ satisfying Requirement 1, with the potential $\psi_m$ for $m\in \mathcal{M}$ arising as a geodesic current.
\begin{enumerate}
    \item $\Teich(S)$ or marked hyperbolic metrics~\cite[Proposition~14]{Bon88}. 
    \item Marked negatively curved Riemannian metrics~\cite[Th\'eor\`eme~1]{Otal90:SpectreMarqueNegative}. This is the space $\mathcal{M}(S)$ of negatively curved Riemannian metrics on $S$ where two metrics are equivalent if are related by a diffeomorphism isotopic to the identity.
     \item Marked negatively curved Riemannian metrics with conical singularities of angle $\geq 2\pi$~\cite[Theorem~A]{HP97:RigidityNegCurvedCone}, generalizing the previous item.
    \item Marked non-positively curved Riemannian metrics~\cite[Theorem~A]{CFF92:RigidityNonPosCurvedRiem}.
    \item Marked singular Euclidean cone metrics from quadratic differentials~\cite[Lemma~9]{DLR2010}.
   
    \item Non-positively curved Euclidean cone metrics with cone angles $\geq 2\pi$~\cite[Proposition~3.3]{BL17:RigidityFlat} (this subsumes the previous case).
\end{enumerate}

In forthcoming work~\cite{MGT25:Intersections}, the second author and D.~Thurston develop a general framework to realize length functions as geodesic currents, encompassing these and many additional examples.

\subsection*{Rigidity results}
For each of the examples above, normalized marked length spectrum rigidity is known, and thus Requirement 2 is satisfied (the numbering correlates with the list above):
\begin{enumerate}
\item For any marked hyperbolic metric $\bx \in \Teich(S)$, let $L_{\bx}$ denote the hyperbolic Liouville current associated to $\bx$. Bonahon shows that in~\cite[Theorem~12]{Bon88} that $\bx \mapsto [L_{\bx}]$ is a proper embedding.
His result boils down to an application of marked length spectrum rigidity for hyperbolic metrics, a result that dates back to work of Fricke--Klein~\cite{FK65:HyperbolicSpectrum}.
In much more recent work~\cite{ELS22:HyperbolicCone}, a rigidity result in terms of the support of the Liouville current is obtained for a generic subspace of \emph{conical} hyperbolic metrics.
\item For negatively curved Riemannian metrics of area 1. 
In this case, normalized marked length spectrum rigidity is the same as marked length spectrum rigidity, which follows from~\cite[Th\'eor\`eme~2]{Otal90:SpectreMarqueNegative}. 
Note that normalizing is necessary here since scaling a negatively curved metric by a positive number yields a different negatively curved metric but their induced pseudo-metrics are roughly similar (in terms of geodesic currents, they are multiples of each other). Another normalization could have been entropy 1.
\item The same cited paper proves rigidity.
\item The above cited paper also proves normalized rigidity for non-positively curved Riemannian metrics of area 1 (or, alternatively, entropy 1)~\cite[Lemma~1.5]{CFF92:RigidityNonPosCurvedRiem}.
\item Rigidity is likewise proven for singular Euclidean metrics of area 1 associated to quadratic differentials~\cite[Theorem~2]{DLR2010}.
\item Same cited paper above proves rigidity.
\end{enumerate}

In the next subsection we will discuss a relevant example of proper hyperbolic potentials to which our asymmetric metric extends, arising from Anosov representations.

\subsection{Extension to Anosov representations} 

\label{sec, AnosovReps}

We give a quick account of Anosov representations and discuss how they induce points in the space of metric structures $\mathcal{D}_\Gamma$ and, more generally, metric potentials $\mathcal{H}^{++}_\Gamma$ (see Appendix~\ref{sec:asymmetric} for definitions).

Anosov representations were introduced by Labourie~\cite{Lab} and extended by Guichard--Wienhard \cite{GW} to general Gromov hyperbolic groups. They are class of discrete representations with finite kernel into semisimple Lie groups sharing many features with convex cocompact representations in rank one.

\subsubsection{Structure of semisimple Lie groups}\label{subsec: preliminaries}

 Let $\g$ be a connected real semisimple algebraic group of non compact type with
 Lie algebra $\lieg$.
Fix a maximal compact subgroup $\ko<\g$ with Cartan involution $\tau$ on $\lieg$, and set $\liep:=\{v\in\lieg:\tau v=-v\}$. Fix a Cartan subspace $\liea\subset\liep$. Let $\m=\mathsf{Z}_{\ko}(\liea)$ and $\mathsf{N}_{\ko}(\liea)$ be its centralizer and normalizer, respectively. Let $\Sigma\subset\liea^*$ be the system of restricted roots; for $\alpha\in\Sigma$ the root space is
\[
\lieg_{\alpha}:=\{Y\in\lieg:\ [H,Y]=\alpha(H)Y\ \text{ for all }H\in\liea\}.
\]
Choose a positive system $\Sigma^+$ with simple roots $\Pi$. Denote by $\overline{\liea^+}$ the closed positive Weyl chamber, and by $w_0$ the longest element of the Weyl group $\w\simeq \mathsf{N}_{\ko}(\liea)/\m$; the opposition involution is $\iota:=-w_0$.

For $\Theta\subset\Pi$, set
\[
\mathfrak{p}_{\Theta}
   := \lieg_0 \oplus \bigoplus_{\alpha\in\Sigma^+} \lieg_{\alpha}
        \oplus \bigoplus_{\alpha\in\langle \Pi-\Theta\rangle} \lieg_{-\alpha}.
\]

\[
\overline{\mathfrak{p}_{\Theta}}
   := \lieg_0 \oplus \bigoplus_{\alpha\in\Sigma^+} \lieg_{-\alpha}
        \oplus \bigoplus_{\alpha\in\langle \Pi-\Theta\rangle} \lieg_{\alpha}.
\]

where $\lieg_0=\mathsf{Z}_{\lieg}(\liea)$, and $\langle \Pi - \Theta \rangle$ denotes the set of positive roots generated by roots in $\Pi-\Theta$. Let $\p_\Theta$ and $\overline{\p}_\Theta$ be the corresponding parabolics.
Every parabolic subgroup of $\g$ is conjugate to a unique $\p_\Theta$, for some $\Theta \subset \Pi$, and $\overline{\p}_\Theta$ is conjugate to $\p_{\iota(\Theta)}$, where $\iota(\Theta)\coloneqq \{ \alpha \circ \iota : \alpha \in \Theta\}$.

Let $\f_\Theta:=\g/\p_\Theta$, $\overline{\f}_\Theta:=\g/\overline{\p}_\Theta$ the \emph{flag manifolds} of $\g$.

 A pair of flags in $\overline{\f}_\Theta\times\f_\Theta$ are \emph{transverse} if they belong to $\overline{\f}_\Theta^{(2)}$, the unique open orbit of the action of $\g$ on $\overline{\f_\Theta} \times \f_\Theta$. 

\begin{example}\label{ex: roots}
When $\g=\mathsf{PSL}(V)$ for a real (resp. complex) $d$–dimensional vector space $V$, we can let $\liea$ to be the diagonal traceless matrices in an orthogonal (resp. unitary) basis. Writing $\varepsilon_i$ for the evaluation on the $i$–th diagonal entry, we have $\Sigma=\{\varepsilon_i-\varepsilon_j:i\ne j\}$, $\Sigma^+=\{\varepsilon_i-\varepsilon_j:i<j\}$, and $\Pi=\{\alpha_{i,i+1}=\varepsilon_i-\varepsilon_{i+1}:1\le i\le d-1\}$ (which will often be denoted by $\alpha_i$).
\end{example}

 For $\Theta\subset\Pi$, $p_\Theta:\liea\to\liea_\Theta:=\bigcap_{\beta\in\Pi-\Theta}\ker(\beta)$ denotes the unique projection invariant under the group $\w_\Theta \coloneqq \{ w \in \w : w|_{\liea_\Theta}=id_{\liea_\Theta}\}$.
 The dual linear space $\liea_\Theta^*$ identifies with $\{ \varphi \in \mathfrak{a}^* : \varphi \circ p_\Theta = \varphi \}$.

Let $\nsf$ be the \emph{unipotent radical} of $\p=\p_\Pi$, i.e., the connected subgroup of $\g$ associated to the Lie algebra $\sum_{\alpha \in \Sigma^+}\lieg_\liea$. If we let $\asf \coloneqq \exp \liea$, the \emph{Iwasawa decomposition} is
\[
\g = \ko \asf \nsf.
\]
In particular, $\mathcal{F} \cong \ko / \m$, and for $\xi \in \mathcal{F}$, we can find $k \in \ko$, 
so that $k \m = \xi$. Quint~\cite{Qui02:PS} defines a map $\sigma \colon \g \times \mathcal{F} \to \liea$ by \[gk = l \exp(\sigma(g,k\m))n,\]
where $n \in \nsf$, and $l \in \ko$.
In~\cite[Lemme~6.11]{Qui02:PS} it is shown that $p_\Theta \circ \sigma \colon \g \times \mathcal{F} \to \liea_\Theta$ factors through a map $\sigma_\Theta (gh,\xi)=\sigma_\Theta(g,h\cdot \xi) + \sigma_\Theta(h, \xi)$.
This map $\sigma_\Theta$ is called the \emph{$\Theta$-Busemann-Iwasawa cocycle} of $\g$.

Let $\langle \cdot, \cdot\rangle$ the inner product on $\liea^*$ dual to the Killing form of $\lieg$.
For $\varphi, \psi \in \liea^*$ set
\[
\langle \varphi,\psi \rangle \coloneqq 2\frac{(\varphi,\psi)}{(\psi,\psi)}.
\]
Given $\alpha \in \Pi$, the \emph{fundamental weight} is the functional $\omega_\alpha \in \liea^*$ defined by the formulas
$\langle \omega_\alpha, \beta\rangle = \delta_{\alpha\beta}$ for $\beta \in \Pi$. One has
\[
\omega_\alpha \circ p_\Theta= \omega_\alpha
\]
for all $\alpha \in \Theta$. In particular, $\omega_\alpha \in \liea_\Theta^*$.
Fundamental weights are related to a special set of representations of $\g$ introduced by Tits.
If $\Lambda \colon \g \to \operatorname{PGL}(V)$ is an irreducible representation, a functional $\chi \in \liea^*$ is a \emph{weight} of $\Lambda$ if the \emph{weight space} 
\[
V_\chi \coloneqq \{ \Lambda(\exp(X))\cdot v = e^{\chi(X)}v, \mbox{ for all } X \in \liea \}
\]
is non zero. 
Tits shows the existence of a unique weight $\chi_\Lambda$ which is maximal with respect to the order $\chi\geq \chi'$ if $\chi-\chi'$ is a linear combination of simple roots with non-negative coefficients. The functional $\chi_\Lambda$ is called the  \emph{highest weight} of $\Lambda$ and the representation is \emph{proximal} if the associated weight space $V_{\chi_\Lambda}$ is one dimensional.

We let $||\cdot||_\Lambda$ be a \emph{good norm}~\cite[Section~4.2]{Qui02:PS} on $V$, i.e., a norm on $V$ invariant under $\Lambda (K)$ and such 
that $\Lambda (\asf)$ consists of semihomotheties (i.e.\ diagonal matrices on an
orthonormal basis $\mathcal{E}$ of $V$).

For every $g \in G$, one has
\[
\log ||\Lambda(g)|| = \chi_\Lambda^{\omega}(a(g)), \qquad
\log \lambda_1(\Lambda (g)) = \chi_\Lambda^\omega(\lambda(g))
\]
where $\lambda_1(\Lambda(g))$ is the logarithm of the modulus of the highest eigenvalue of $\Lambda(g)$
Denote by $W_{\chi_\Lambda}$ the $\Lambda (A)$–invariant complement of $V_{\chi_\Lambda}$.
The stabilizer in $G$ of $W_{\chi_\Lambda}$ is 
$\overline{\p}_{\Theta_\Lambda}$, and thus one obtains a map of flag 
spaces
\[
(\Xi_\Lambda, \Xi_\Lambda^*) \colon 
\mathcal{F}_{\varphi_\phi}^{(2)} \longrightarrow 
\operatorname{Gr}^{(2)}_{\dim V_{\chi_\Lambda}}(V).
\]

Here
\[
\Theta_\Lambda 
= \{\, \sigma \in \Pi : \chi_\Lambda - \sigma \text{ is a weight of } \Lambda \,\}.
\]

\begin{proposition}[Tits \cite{Tits}]\label{prop: tits}
For every $\alpha\in\Pi$ there exists a finite-dimensional real representation
\[
\Lambda_\alpha:\g\to\mathsf{PGL}(V_\alpha)
\]
that is irreducible and proximal, whose highest weight is $\chi_\alpha=l_\alpha\omega_\alpha$ with $l_\alpha\in\mathbb{N}$. Moreover, $\Lambda_\alpha(\p_\Theta)$ stabilizes a line in $V_\alpha$ for every $\Theta\subset\Pi$ containing $\alpha$.
\end{proposition}

Fix such a family $\{\Lambda_\alpha\}_{\alpha\in\Pi}$. We have
\begin{equation}\label{eq: highest weight and ptheta}
\chi_\alpha\circ p_\Theta=\chi_\alpha\quad\text{for }\alpha\in\Theta,
\end{equation}
so $\chi_\alpha\in\liea_\Theta^*$ whenever $\alpha\in\Theta$.

\medskip
 The \emph{Cartan projection} $\mu:\g\to\overline{\liea^+}$ is defined by the unique element $\mu(g) \in \liea^+$ so that
\[
g\in \ko\,\exp(\mu(g))\,\ko.
\]

Let $g \in G$, and let $g=k_g z_g l_g$ denote its Cartan decomposition.
We say $g \in G$ has \emph{gap at $\Theta$} if for all $\sigma \in \Theta$ one has $\sigma(\mu_\Theta(g))>0$.
In that case, the \emph{Cartan attractor} of $g$ in $\mathcal{F}_\Theta$ is defined as
\[
U_\Theta(g) = k_g [P_\Theta].
\]

The \emph{Jordan projection} is
\begin{equation}
\lambda(g):=\lim_{n\to\infty}\frac{\mu(g^n)}{n}\in\overline{\liea^+}.
\label{eq:JordanVSCartan}
\end{equation}
For $g\in\g$ and all $\alpha\in\Pi$, we have
\begin{equation}\label{eq: spectral radious and fund weight}
\lambda_1\big(\Lambda_\alpha(g)\big)=\chi_\alpha\big(\lambda(g)\big)=l_\alpha\,\omega_\alpha\big(\lambda(g)\big).
\end{equation}

Further, we define $\mu_\Theta:=p_\Theta\circ\mu$ and $\lambda_\Theta:=p_\Theta\circ\lambda$.

Following Quint~\cite[Section~4.2]{Qui02:PS}, each $V_\alpha$ is endowed with a good norm $|| \cdot ||_\alpha$, 
so that for any $v \in \liea^+$
\[
||\Lambda_\alpha(\exp v) || = e^{\omega_\alpha(v)}.
\]

\subsubsection{Semi--similarity and proximality}

We follow~\cite[Section~2.2]{Qui02:PS}. A linear operator $T \colon V \to V$ on a finite-dimensional normed vector space is called \emph{semi-similarity} if there exist
a direct sum decomposition:
\[
V=V_T^M \oplus V_T^m
\]
into $T$-invariant subspaces so that
\begin{enumerate}
	\item $T|_{V_T^M}$ is a \emph{similarity} of ratio $|| T ||$, i.e., $||T(v)||=||T || \cdot ||v||$ for all $v \in V_T^M$, where we use the same notation $|| \cdot ||$ for the induced operator norm.
	\item all eigenvalue moduli on $V_T^M$ are strictly smaller than $||T ||$.
\end{enumerate}
We say that $T$ is a \emph{proximal} semi-similarity if, in addition, $V_T^M$ is one-dimensional.

\subsubsection{Busemann cocycle}
\label{subsec:busemann}
For $g \in G$ and a flag $\xi \in P_\Theta$, the $\Theta$-Busemann cocycle is a map $\sigma_\Theta \colon G \times P_\Theta \to E_\Theta$ so that, via the representation $\Lambda_\alpha$, satisfies
\[
\omega_\alpha(\sigma_\Theta(g,\xi)) = \log \frac{|| \Lambda_\alpha(g)v||_\alpha}{||v||_\alpha}
\]
where $v$ is in the line representing $\xi$.
There is a notion of Gromov product defined in~\cite{Sam15:Orbital}.

$\mathcal{G}_\Theta : \mathcal{F}_\Theta^{(2)} \to E_\Theta$ is defined so that
	\[
	l_\alpha \omega_\alpha(G_\Theta(p,q))
	= \log \frac{||\phi(v)||}{||\phi||_\alpha ||v||_\alpha} \]
    where $\phi \in \Xi_{\alpha}^*(x)$ and $v \in \Xi_\alpha(y)$, and $\Xi_\alpha, \Xi_\alpha^*$ are the limit maps of the representations $V_\alpha$, for each $\alpha \in \Theta$.

Recall that the \emph{Cartan basin} of $g \in \g$ is defined, for $\alpha>0$, as
\[
B_{\Theta, \alpha}(g) \coloneqq \{ x \in \mathcal{F}_\Theta : \omega_\sigma(\mathcal{G}_\Theta(U_{\iota \Theta}(g^{-1})),x) >-\alpha \mbox{ for all } \sigma \in \Theta \}.
\]
	By \cite[Lemme~6.7]{Qui02:PS}, given $\eta>0$, there exists a constant $\kappa_\eta>0$ so that if $y \in B_{\Theta, \eta}(g)$, then
	\[
	||\mu_\Theta(g) - \sigma_\Theta(g,y)|| < \kappa_\eta.
	\]

\subsubsection{Anosov representations and their length functions}\label{subsec:lengths_anosov}

 Fix $\Theta\subset\Pi$. A representation $r:\Gamma\to\g$ is \emph{$\p_\Theta$–Anosov} if there exist constants $C,c>0$ such that for every $\alpha\in\Theta$ and every $\gamma\in\Gamma$,
\begin{equation}
\alpha\big(\mu(r(\gamma))\big)\ \ge\ C\,|\gamma| - c.
\label{eq:Anosov_equation}
\end{equation}
 It follows from the definition that Anosov representations are quasi-isometric embeddings, in particular discrete with finite kernel. 
 By a result of Kapovich-Leeb-Porti~\cite[Theorem~6.15]{KLP18:Morse} it follows that if $r:\Gamma\to\g$ is $\Theta$-Anosov then $\Gamma$ is Gromov hyperbolic. We shall assume that $\Gamma$ is non-elementary and denote by $\partial\Gamma$ its Gromov boundary (see Appendix~\ref{sec:asymmetric} for a definition). We also let $\bgs$ be the space of ordered pairs of different points in $\bg$. Every infinite order element $\gamma\in\Gamma$ has a unique attracting (resp. repelling) fixed point in $\bg$, denoted by $\gamma_+$ (resp. $\gamma_-$). We let $\gh\subset\Gamma$ be the subset consisting of infinite order elements. The conjugacy class of $\gamma\in\Gamma$ is denoted by $[\gamma]$, and the set of conjugacy classes of elements of $\Gamma$ (resp. $\gh$) will be denoted by $[\Gamma]$ (resp. $[\gh]$).

A central feature of $\Theta$-Anosov representations is that they admit \textit{limit maps}. By definition, these are H\"older continuous, $r$-equivariant, dynamics preserving maps
$$\xi_r: \partial \Gamma \to \mathscr{F}_\Theta \textnormal{ and } \overline{\xi}_r: \partial \Gamma \to \overline{\mathscr{F}}_\Theta,$$
\noindent which are moreover \textit{transverse}, that is, for every $x\neq y$ in $\partial\Gamma$ one has
$$ (\overline{\xi}_r(x),\xi_r(y))\in\f^{(2)}_\Theta. $$
\noindent The limit maps exist and are unique (see \cite{BPS,GGKW, KLPanosovcharacterizations} for details).

The set of $\Theta$-Anosov representations from $\Gamma$ to $\g$ is an open subset of the space of all representations $\Gamma\to\g$~\cite{Lab,GW}.

Projective Anosov representations capture all Anosov representations.

\begin{proposition}\label{prop:proj_anosov_tits}
Let $r:\Gamma\to\g$ be  $\Theta$-Anosov. Then for every $\alpha\in\Theta$ the representation $\Lambda_\alpha\circ r:\Gamma\to\mathsf{PGL}(V_\alpha)$ is projective Anosov.
\end{proposition}

We denote by $\ha$ the space of conjugacy classes of $\p_\Theta$-Anosov representations from $\Gamma$ to $\g$.

Carvajales-Dai-Pozzetti-Wienhard in~\cite{CDPW24:Asymmetric} rely on work of Sambarino \cite{HyperconvexRepsExponentialGrowth,Quantitative} to associate to each $r\in\ha$ certain flow space which is a H\"older reparameterization of the Mineyev flow of $\Gamma$~\cite{Mineyev}. Then, they use dynamical techniques, such as  Thermodynamical Formalism to study $\ha$.

In order to define the lengths and entropies associated to Anosov representations, we introduce the Benoist  limit cone~\cite{Benoist_AsymtoticLinearGroups} for general discrete subgroups of $\g$.

\begin{definition}
The $\Theta$-\textit{limit cone} of $r\in \ha$ is the smallest closed cone $\zcone\subset \mathfrak{a}_\Theta ^+$ containing the set $\{\lambda_{\Theta}(r (\gamma)): \gamma\in \Gamma\}$. The \textit{limit cone} $\mathscr{L}_r$ of $r$ is the $\Pi$-limit cone.
\end{definition}

Under the assumption that $r$ is Zariski dense, Benoist \cite{Benoist_AsymtoticLinearGroups} shows that $\mathscr{L}_r$ and the $\Theta$-limit cone are convex cones with non empty interior.

Let $$\dzcone:=\{\varphi \in \mathfrak{a}_{\Theta}^{*}: \varphi|_{\zcone}\geq 0\}$$
\noindent be the \textit{dual cone}. We denote by $\text{int} (\dzcone)$ the interior of $\dzcone$, that is, the set of functionals in $\liea_\Theta^*$ positive on $\zcone\setminus\{0\}$.

Fix a functional $$\varphi \in \bigcap_{{r}\in\ha} \text{int} (\dzcone).$$
\noindent The above intersection is non empty: for instance, it contains $\lambda_1$ and $\omega_\alpha$ for all $\alpha\in\Pi$.

\begin{definition}
The $\varphi$-\textit{marked length spectrum} (or simply $\varphi$-\textit{length spectrum}) of $r\in\ha$ is the function $\ell^\varphi_{r}:\Gamma \to \rr_{\geq 0} $ given by
    \begin{align*}
    \ell^\varphi_{r}&(\gamma):= \varphi(\lambda_{\Theta}(r(\gamma))).
    \end{align*}
    \end{definition}
\noindent Observe that for a $\Theta$-Anosov representation $r$,  $\ell_r^\varphi(\gamma)>0$ if and only if has infinite order. Furthermore the $\varphi$-length spectrum is invariant under conjugation in $\Gamma$ and therefore descends to a function $[\Gamma]\to\rr_{\geq 0}$, which we shall denote by $\ell^\varphi_{r}$ as well.

\begin{definition}
The $\varphi$-\textit{entropy} of $r$ is defined by $$h_r^\varphi:=\displaystyle\limsup_{t\to\infty}\frac{1}{t}\log\#\{[\gamma]\in[\Gamma]:\ell_r^\varphi(\gamma)\leq t\}\in[0,\infty].$$
\end{definition}

The $\varphi$-entropy of $r$ was first introduced by Sambarino \cite{HyperconvexRepsExponentialGrowth,Quantitative}. He showed that it is defined by a true limit, is positive, finite, and coincides with the topological entropy of a certain flow associated to $r$ and $\varphi$, arising from a H\"older reparameterization of the so-called Mineyev flow on~\cite{Mineyev}, a coarse generalization of the geodesic flow for hyperbolic groups.

\subsubsection{Hyperbolic potentials associated to an Anosov representation}

There are natural pseudo-metrics associated to general $\Theta$-Anosov representations of surface groups. Some of these length functions can be
realized as geodesic currents via positive generalized cross-ratios
\cite[Proposition~2.24]{MZ19:PositivelyRatioed}, \cite[Theorem~1.6]{BIPP24}.
However, many points arise from more general pseudo-metrics that are not
geodesic currents (e.g.\ quasi-Fuchsian non-Fuchsian representations
\cite{FF22:QF,BR24:Approximating}).
More generally, for any non-elementary hyperbolic group one has similar symmetric examples coming from the Hilbert length spectrum of the representation.
However, there are also length functions that do not
arise from stable lengths of pseudo-metrics, but instead from proper hyperbolic metric
potentials, which allow asymmetry and certain negative values,  but with strictly positive stable length. See
Appendix~\ref{sec:asymmetric} for further details.

In Appendix~\ref{sec:asymmetric}, we give the precise definition of proper hyperbolic metric potential and show that metric $d$ can be extended to the space of metric potentials $\mathcal{H}^{++}_\Gamma$.

We now show that for any $\Theta$-Anosov representation $r$, and an $\varphi$ in the interior $\Theta$-limit cone of $r$, there exists a potential in $H_\Gamma^{++}$ recovering the length $\ell_r^\varphi$-length.

\begin{proposition}
\label{prop:alpha_metric}
Let $r:\Gamma\to\g$ be $\Theta$–Anosov and $\varphi \in \Int(\dzcone)$.
Then the assignment
\[
\psi_{r,\varphi}(x,y)\coloneqq \varphi(\mu_{\Theta}(r(x^{-1}y))),
\]
defines an element in $H_\Gamma^{++}$ with stable length satisfying
\[
\ell_{\psi_{r,\varphi}}(\gamma) = \ell_r^\varphi(\gamma)
\]
for every $\gamma \in \Gamma$.
\end{proposition}

\begin{proof}
	First, we claim that $
	\psi_{r,\varphi}\in H_\Gamma$. To that end, we show that for any constants $C,R>0$, there exist a constant $c_0>0$ so that for any $C$-rough geodesic segment $(\eta(i))_{i\in \{1,\cdots,k\}}$ between $x$ and $y$ (see definition in Subsection~\ref{subsec:pseudo-metrics}), and for any $z \in U_R(\eta)$, 
	\[
	\varphi(\mu_\Theta(x^{-1}y)) \geq \varphi(\mu_\Theta(x^{-1}z))+\varphi(\mu_\Theta(z^{-1}y))-c_0.
	\]
	By~\cite[Proposition~3.1]{CRS24:Joint} the claim will follow.
	By translating by $z$ on the left, we can assume that the $C$-rough geodesic passes within distance $R$ of the identity $o$.
	
	Let $c_R\coloneqq\inf \{ \min_{\alpha \in \Theta } \exp \omega_\alpha(\mathcal{G}_\Theta(\overline{\xi}_{\Lambda_\alpha}(\eta),\xi_{\Lambda_\alpha}(\zeta)) :   (\eta,\zeta)_o < R+4\delta \}$. 
    For $g \in \Gamma$, let $g_\alpha\coloneqq \Lambda_\alpha(g)$, for every $\alpha \in \Theta$.
	By the same argument as in~\cite[Lemma~7.1]{CT25:Manhattan}, we have $c_R>0$ and there exists a constant $L_R>0$ so that for a geodesic between $x$ and $y$ passing through $o$ and $|x|,|y|>L_R>0$, we have
	\[
	\min_{\alpha \in \Theta} \exp \omega_\alpha(\mathcal{G}_\Theta(U_{\iota \Theta}(x_\alpha),U_\Theta(y_\alpha))) > \frac{c_R}{2}.
	\] 
If we put $V^M_{\alpha,y} \coloneqq U_\Theta(y_\alpha)$, this says $V^M_{\alpha,y} \in B_{\Theta,\epsilon}(x_\alpha^{-1})$ with $\epsilon=c_R/2$. Taking $c_R$ smaller if necessary, so that $\epsilon \leq 1$, 
we apply Lemma~\ref{lem:trans} from which it follows there exists a constant $K>0$ so that
\[
|| \mu_\Theta(x_\alpha^{-1}y_\alpha)-\mu_\Theta(y_\alpha) - \sigma_\Theta(x_\alpha^{-1},U_\Theta(y_\alpha))) || \leq K .
\]
	By the discussion in Subsection~\ref{subsec:busemann}, we have
	\[
	|| \mu_\Theta(x^{-1}_\alpha) - \sigma_\Theta(x^{-1}_\alpha,U_\Theta(y_\alpha)) || < \kappa_{c_R/2}.
	\]
	Combining all of the above, there exist a constant $c_0 \coloneqq K + \kappa_{c_R/2}>0$ (which depends on $C$ and $R$)  so that
	\[
\varphi(\mu_\Theta(x^{-1}y)) \geq \varphi(\mu_\Theta(x^{-1})) + \varphi(\mu_\Theta(y)) - c_0.
	\]
	This proves $\psi_{r,\varphi}\in H_\Gamma$, as desired.
	To show that $\psi_{r,\varphi}\in H^{++}_\Gamma$, note that for every $\alpha \in \Theta$, the intersection $\ker \omega_\alpha \cap \zcone$ is empty, and $\ker \varphi \cap \zcone$, is also empty by assumption on $\varphi$.
	Then, $\varphi/\omega_\alpha$ is bounded below away from zero on $\zcone$, and thus by~\cite[Proposition~5.4.2]{S24:Dichotomy},
there are constants $q,Q>0$ so that
	\[
	\varphi(\mu_\Theta(r(g))) \geq q \omega_\alpha(\mu_\Theta(r(g)))-Q
	\]
    for every $g \in \Gamma$.
	Combining the fact that $\Lambda_{\alpha} \circ r$ is projective Anosov, by Proposition~\ref{prop:proj_anosov_tits}, and the inequality above, we get constants $Q',q'>0$ so that
	\[
	\varphi(\mu_\Theta(r(g))) \geq q' |g| - Q'
	\]
	for some/any word-metric $|\cdot|$, and every $g \in \Gamma$.
	By taking powers $g^n$ and dividing by $n$, and using Equation~\ref{eq:JordanVSCartan}, we get $\ell_{\psi_{r,\varphi}}=\lim_n \frac{\varphi(\mu_\Theta(r(g^n)))}{n} \geq c \ell_{|\cdot|}$.
	This shows that $\ell_{r,\varphi}(g)=\ell_{\psi_{r,\varphi}}$ and, since the word-metric $|\cdot|$ is an element in $D_\Gamma$, it follows $\psi_{r,\varphi} \in H_\Gamma^{++}$.
\end{proof}

The following is the key lemma used in the previous result.

\begin{lemma}
	Fix $\Theta\subseteq\Pi$ and a norm $|| \cdot ||$ on $E_\Theta$.
	For each $\alpha\in\Theta$, let $\Lambda_\alpha : G\to\mathrm{PGL}(V_\alpha)$ be one of the representations in Proposition~\ref{prop: tits} equipped with a good norm. Let $1 \geq \epsilon>0$. Then there exists a constant $K>0$,
	depending on $\epsilon$ such that the following holds for every $\alpha \in \Theta$.
	For every $g,h\in \Gamma$, let $g_{\alpha}\coloneqq \Lambda_\alpha(g), h_{\alpha}\coloneqq \Lambda_\alpha(h)$, and suppose \[V^M_{\alpha, h_\alpha} \in B_{\Theta, \epsilon}(g_\alpha).\]
	Then, we have
	\[
	|| \mu_\Theta(g_\alpha h)-\mu_\Theta(h_\alpha) - \sigma_\Theta(g_\alpha,U_\Theta(h_\alpha))||
	\leq K.
	\]
	\label{lem:trans}
\end{lemma}

\begin{proof}
	Fix $\alpha\in\Theta$. Then
	$h_\alpha$ is a proximal semi--similarity on $V_\alpha$ with
	\[
	V_\alpha = V_{\alpha,h}^M \oplus V_{\alpha,h}^m.
	\]
	Let $v_\alpha\in V_{\alpha,h_\alpha}^M\setminus\{0\}$.
	By \cite[Lemma~2.6]{Qui02:PS} (left inequality), we have
	\[
	\frac{||g_\alpha(v_\alpha)||}{||v_\alpha||}
	\ \le\
	\frac{||g_\alpha h_\alpha||}{||h_\alpha||}.
	\]
	
	On the other hand, apply \cite[Lemme~4.2]{Qui02:PS} to the proximal semi-similarity $g_\alpha$ to obtain
	\[
	||g_\alpha(v_\alpha)|| \ \ge\ \epsilon\,||g_\alpha||\,||v_\alpha||,
	\]
	and since
	$||g_\alpha h_\alpha||\le ||g_\alpha||\,||h_\alpha||$ we obtain
	\[
	\frac{||g_\alpha h_\alpha||}{||h_\alpha||}
	\ \le\
	\frac{1}{\epsilon}\,\frac{||g_\alpha(v_\alpha)||}{||v_\alpha||}
	\]
    for every $\alpha \in \Theta$.
	
	By \cite[Lemme~6.4]{Qui02:PS}), this can be
	rewritten as
	\[
	e^{\varpi_\alpha(\sigma_\Theta(g_\alpha,U_\Theta(h_\alpha)))}
	\ \le\
	e^{\varpi_\alpha(\mu(g_\alpha h_\alpha)-\mu(h_\alpha))}
	\ \le\
	\frac{1}{\epsilon}\,
	e^{\varpi_\alpha(\sigma_\Theta(g_\alpha,U_\Theta(h_\alpha)))}.
	\] 
 By \cite[Lemma~3.1]{Qui02:PS}, we obtain
	a uniform bound
	\[
	||\mu_\Theta(g_\alpha h_\alpha)-\mu_\Theta(h_\alpha)) - \sigma_\Theta(g_\alpha,U_\Theta(h_\alpha))\,||
	\ \le\ K,
	\]
	for some constant $K>0$ depending only on $\epsilon$,
	as claimed.
\end{proof}

\begin{remark}
The result above generalizes the case $\varphi=\omega_{\alpha}$ for every $\alpha \in \Pi$ (simple roots),
by~\cite[Lemma~7.1]{CT25:Manhattan}. 
The statement above for $\iota$-invariant $\varphi \in (\mathfrak{a}_{\Theta}^+)^*$ can be derived using~\cite[Theorem~4.5]{DK22:PSAnosov}.
\end{remark}

\begin{proof}[Proof of Theorem~\ref{extensions}]
In~Proposition~\ref{prop:asymmetric_defined} we show $d$ extends to the space of metric potentials $\mathcal{H}^{++}_\Gamma$.
By Proposition~\ref{prop:alpha_metric} for every $\Theta$-Anosov representation $r$ and every $\varphi \in \Int(\dzcone)$, there exists a proper hyperbolic potential $\psi_{r,\varphi} \in H_\Gamma$ so that
\[
\ell_{\psi_{r,\varphi}}(g) = \ell_r^{\varphi}(g)
\]
for every $g \in \Gamma$.
From this, it follows that
\[
h(\psi_{r,\varphi}) = h_r^{\varphi}.
\]
It suffices to show that restricting $d$ to these potentials induced by the functional $\varphi$ gives a non-degenerate metric on the space of $\Theta$-Anosov representations.
If $d(\psi_{r,\varphi},\psi_{r',\varphi})=0$, it follows that
\[
\ell_{\psi_{r,\varphi}}(g) h(\psi_{r,\varphi}) = \ell_{\psi_{r',\varphi}}(g) h(\psi_{r',\varphi})
\]
for every $g\in \Gamma$,
i.e.
\[
h_r^\varphi \ell_{r}^{\varphi}(g)=h_{r'}^\varphi \ell_{r'}^{\varphi}(g)
\]
for every $g\in \Gamma$.
It then follows that $r$ and $r'$ are conjugate representations by the marked length spectrum rigidity results for Anosov representations established in~\cite{BridgemanCanaryLabourieSambarino18:SimpleRoots} and more generally in~\cite{Benoist_AsymtoticLinearGroups}, \cite{Bur93:Manhattan}, \cite{DK00:Conjugacy}. For a consolidated statement, see~\cite[Theorem~6.8]{CDPW24:Asymmetric}.  
\end{proof}

In particular, unlike the setup in~\cite{CDPW24:Asymmetric}, our construction
does not rely on H\"older reparameterizations of the geodesic flow (when
$\Gamma=\pi_1(S)$) or on the Mineyev flow in the general non-elementary
hyperbolic case~\cite{Quantitative,S24:Dichotomy}.

\section{Further questions}
\label{sec:questions}
It would be interesting to understand whether analogs of \cref{horoboundary} and \cref{isometricrigidity} hold for the symmetrized Thurston metric $d_{\mathrm{sym}}$. The isometric rigidity  for $d_{\mathrm{sym}}$ does not follow directly from \cref{isometricrigidity}. Moreover, the horofunction compactification of $\mathbb{P}\mathcal{C}_{\mathrm{fill}}$ with respect to $d_{\mathrm{sym}}$ is not yet known. In fact, one would have to adapt extended map $\Psi: \mathbb{P}\mathcal{C} \rightarrow C(\mathbb{P}\mathcal{C}_{\mathrm{fill}})$ introduced in \cref{eq3} in order to prove \cref{horoboundary} for $d_{\mathrm{sym}}$.

\begin{question}
For distinct closed surfaces $S_{g}$ and $S_{g'}$, can $(\mathbb{P}\mathcal{C}_{\mathrm{fill}}(S_{g}), d_{\mathrm{sym}})$ and $(\mathbb{P}\mathcal{C}_{\mathrm{fill}}(S_{g'}), d_{\mathrm{sym}})$ be isometric?
\end{question}

Another natural, seemingly more challenging, question is to ask for quasi-isometric rigidity for any of the two metrics.

\begin{question}
For distinct closed surfaces $S_{g}$ and $S_{g'}$, is it true that $(\mathbb{P}\mathcal{C}_{\mathrm{fill}}(S_{g}), d_{\mathrm{sym}})$ (resp. $(\mathbb{P}\mathcal{C}_{\mathrm{fill}}(S_{g}), d)$ and $(\mathbb{P}\mathcal{C}_{\mathrm{fill}}(S_{g'}), d_{\mathrm{sym}})$ (resp. $(\mathbb{P}\mathcal{C}_{\mathrm{fill}}(S_{g'}), d)$ are not quasi-isometric?
\end{question}

In the Appendix~\ref{sec:asymmetric} we explain that the space of projective filling currents $(\PC_{fill},d)$ isometrically embeds into the strictly bigger space of metric structures $\mathcal{D}_{\pi_1(S)}$ originally introduced in~\cite{Fur02:Coarse} and recently studied in~\cite{Eduardo,CT25:Manhattan,StephenEduardo}.
It would be interesting to extend the horoboundary construction to this space.

\begin{question} Can one construct a horofunction boundary for the space of metric structures $\mathcal{D}_{\pi_1(S)}$, relying on an extended framework of horoboundary functions that does not assume properness of the underlying space, such as the one described in~\cite{Cor18:HilbertInfinite}? More generally, one would like to construct the horofunction boundary for $\mathcal{D}_\Gamma$ of any non-elementary hyperbolic group $\Gamma$ (see details in Appendix~\ref{sec:asymmetric}).
Does the horofunction compactification of $\mathcal{D}_\Gamma$ coincide with the space of left-invariant hyperbolic pseudo-metrics with the \emph{bounded backtracking property} (in the sense of~\cite{StephenEduardo})?
\end{question}

One of the structural properties of Thurston's asymmetric metric $d_{Th}$ on $\Teich(S)$, proven in~\cite{Th1998}, is that it is a Finsler metric. In higher rank there are analogous asymmetric Finsler metrics, although they are not known to be equivalent to the natural extensions of $d_{Th}$; see~\cite[Theorem~1.4]{CDPW24:Asymmetric}.

\begin{question}
Can $d$ be realized as a natural Finsler metric on the space of projective filling currents (or, more generally, on the space of pseudo-metrics)?
\end{question}

Addressing this question would first require defining an appropriate notion of tangent space to the space of geodesic currents. One possible approach might use the theory of continuous reparametrizations of the geodesic flow~\cite{Tholozan2019CocyclesReparametrizations}.
\appendix

\section{The asymmetric Thurston metric on the space of pseudo-metrics}
\label{sec:asymmetric}

In this appendix we discuss the most general setting in which the asymmetric
metric is defined: the space of proper hyperbolic potentials of a hyperbolic
group. This is a broad class which, in the case of surface groups, contains the
space of projective filling geodesic currents, and for a non-elementary Gromov hyperbolic group contains
the space of metric structures $\mathcal{D}_\Gamma$, consisting of classes of left-invariant Gromov
hyperbolic pseudo-metrics quasi-isometric to a word-metric.

The appendix is primarily expositional. The only genuine novelty is the
observation that the metric $d$ extends to the space of proper hyperbolic
potentials, a fact that follows from minor adaptations of the arguments
in~\cite{CT25:Manhattan}. After indicating the necessary modifications to the
proofs, we show that, when restricted to projective filling currents, this
metric recovers the asymmetric Thurston metric on geodesic currents. This is
the metric that underlies our horoboundary and rigidity results.

For detailed background, we refer the reader to
\cite{Eduardo,CT25:Manhattan,StephenEduardo,LucDid,CRS24:Joint}.

\subsection{The space of pseudo-metrics}
\label{subsec:pseudo-metrics}
Let $(X,\rho)$ be a metric space.
We say $\rho$ is $\delta$-\emph{Gromov hyperbolic} (or just \emph{hyperbolic}), if there exist $\delta>0$ so that for every $x,y,z,w \in \Gamma$, $\rho(x,z) + \rho(y,w) \leq \max \{ \rho(x,y) + \rho(z,w), \rho(x,w) + \rho(y,z) \} +2\delta$.
Basic examples are $\mathbb{H}^2$ or real trees.

Let $\Gamma$ be a Gromov hyperbolic group, i.e. a finitely generated group $\Gamma$ so that for some finite and symmetric generating set $S$ (i.e. $|S|<\infty$ and $S=S^{-1}$) the induced word-metric $\rho_S$ is Gromov hyperbolic. Fundamental examples are surface groups, i.e. $\Gamma=\pi_1(S)$ where $S$ is a closed orientable hyperbolic surface of genus at least 2, as well as free groups.

We let $\partial \Gamma$ denote its \emph{Gromov boundary}, consisting of ideal classes of divergent sequences. A sequence of group elements $(x_n)_{n=0}^{\infty}$ \emph{diverges} if $(x_n | x_m)_o$ for some $o \in \Gamma$, the identity element in $\Gamma$, as $\min \{ n, m\}$ diverges to infinity. 
Two divergent sequences $(x_n)_{n=0}^{\infty}$ and $(y_n)_{n=0}^{\infty}$ are \emph{equivalent} if $(x_n | y_m)_o$ diverges as $\min \{n,m \}$, tends to infinity. This can be thought of as the discrete version of equivalence classes of bounded-distance infinite rays.
We say $\Gamma$ is \emph{non-elementary} if the cardinality of $\partial \Gamma$, $|\partial \Gamma|$ is infinite. In the case $\Gamma=\pi_1(S)$, $\partial \Gamma$ is homeomorphic to $S^1$.

A \emph{pseudo-metric} on $\Gamma$ is a function $\rho \colon \Gamma \times \Gamma \to \mathbb{R}$ so that for every $x,y,z \in \Gamma$, we have $\rho(x,y)\geq 0$, $\rho(x,y)=\rho(y,x)$ and $\rho(x,z) \leq \rho(x,y) + \rho(y,z)$. This is, it is a distance except that the two distinct points are allowed to have distance zero. A pseudo-metric is \emph{left-invariant} if $\rho(gx, gy) = \rho(x,y)$ for every $g,x,y \in \Gamma$.

We say that two pseudo-metrics $\rho,\rho'$ are \emph{quasi-isometric} if there exist constants $c_2,c_1,a>0$ so that
\[
c_1 \rho(x,y) - a \leq \rho'(x,y) \leq c_2 \rho (x,y) + a
\]
for every $x,y \in \Gamma$.
For example, any two hyperbolic metrics on a surface $S$ induce quasi-isometric metrics on $\mathbb{H}^2$ via lifting the metric.

From now on, let us assume that $\Gamma$ is non-elementary and let $D_\Gamma$ be the set of all left-invariant pseudo-metrics on $\Gamma$ quasi-isometric to a fixed word metric on $\Gamma$.
A pseudo-metric $\rho \in D_\Gamma$ is always \emph{$C$-roughly-geodesic} for some $C$ (depending on $\rho$), meaning that for every $x,y \in \Gamma$, there exist a finite set of points $\eta(i)\in \Gamma$ with $i=0,\cdots,k$ so that $x=\eta(0),y=\eta(k)$, and $i- j - C \leq \rho(\eta(i),\eta(j)) \leq i-j +C$ for every $k \geq i>j \geq 0$.
Two pseudo-metrics $\rho,\rho'$ are said to be \emph{roughly similar} if there exist positive constants $\lambda, A$ so that $|\rho - \lambda \cdot \rho'| < A$, i.e., $\rho$ and $\rho'$ are within bounded distance after rescaling by a positive constant. If $\lambda=1$, we say they are \emph{roughly isometric}.
We denote by $\mathcal{D}_\Gamma$ the set of rough similarity classes of elements in $D_\Gamma$.
Two concrete examples of elements in $D_{\Gamma}$ are:
\begin{enumerate}
\item\emph{Word metrics:} Given any finite, symmetric, generating subset $S \subset \Gamma$, the word metric $\rho_{S}$ is an element in $D_{\Gamma}$. Indeed, word metrics are hyperbolic, and are quasi-isometric to any fixed symmetric finite generating set.
\item \emph{Geometric actions:} If $\Gamma$ acts geometrically on a roughly geodesic pseudo-metric space $(X, \rho_{X})$, then for any $w \in X$, the orbit pseudo-metric $\rho^w_X(g,h) \coloneqq \rho_{X}(gw, hw)$ belongs to $D_\Gamma$ (by Milnor-Schwarz lemma).
If we pick another basepoint $w' \in X$, we get, by triangular inequality
\[
|\rho_X^w - \rho_X^{w'}| \leq \rho_{X}(w,w').
\]
Thus $\rho_X^w$ and $\rho_X^{w'}$ are roughly isometric, and they induce points in $\mathcal{D}_\Gamma$.
\end{enumerate}

Let $o$ denote the identity element of $\Gamma$. Given $\rho \in D_\Gamma$, and $x \in \Gamma$, we define the \emph{stable length of x} in $\Gamma$ as $\ell_\rho(x) \coloneqq \lim_n \frac{1}{n} \rho(o, x^n)$ as follows 
\[
\ell_{\rho}(x) = \lim_{n \to \infty} \frac{1}{n}\rho(o, x^n).
\]
This definition only depends on the conjugacy class of $x$ in $\Gamma$, and so defines a function on the set of conjugacy classes of $\Gamma$.
It is straightforward to verify that $\ell_{\rho}(x) > 0$ if and only if $x$ is torsion free. As a consequence the subgroup $\{ x \in \Gamma \colon \rho(o,x) = 0 \}$ is torsion. Since $\Gamma$ is hyperbolic, every torsion subgroup is finite~\cite[Ch.8-Corollarie~36]{GdlH88:Notes}. Note that, for any pseudo-metric $\rho \in D_\Gamma$, and $x, y \in \Gamma$, we have 
\[
\rho (x,y) = \rho(y^{-1}x, o) \geq \ell_{\rho}(y^{-1}x).
\]
If $\rho(x,y) = 0$, then $\ell_{\rho}(y^{-1}x) = 0$ and $y^{-1}x$ is a torsion element. Hence, if $\Gamma$ is torsion-free, any pseudo-metric on $\Gamma$ is an actual metric.

There is another natural source of examples arising from \emph{Green metrics}, metrics obtained by admissible random walks on $\Gamma$, which do not naturally come from actions on a Gromov hyperbolic group~\cite{BHM2011}.

\subsection{Geodesic currents as pseudo-metrics}
\label{subsec:from_currents_to_pseudometrics}
The main example we will be concerned with in this paper are pseudo-metrics arising from geodesic currents on $S$.

If $\Gamma=\pi_1(S)$, given a geodesic current $\mu$ on $S$, and $[(X,\phi)] \in \Teich(S)$, and a point $w \in \mathcal{X}$, one can define a metric on $\Gamma$ as follows
\[
\rho_{\mu, \mathcal{X}}^w(x,y) \coloneqq \rho_{\phi_*\mu}(\phi_*(y) w, \phi_*(y) w)
\]
for every $x,y \in \Gamma$,
where $\phi_*\mu$ is the pushforward of the current $\mu$ on $S$ to a current $\wt{\phi}_*\mu$ on $X$
and $\rho_{\phi_*\mu}$ is a pseudo-metric defined on $\wt{X}$ defined as as follows~\cite{BIPP21,LucDid}: let $[p,q)$, $(p,q]$ denote half-closed half-open $X$-geodesic segments connecting $p,q \in \wt{X}$.
\[
 \rho_{\phi_*\mu}(p,q) \coloneqq \frac{1}{2} \left( \phi_*\mu(G([p,q)) +  \phi_*\mu(G((p,q]\right).
\]
Invariance of the intersection number and \cite[Lemma~4.7]{BIPP21}
\begin{equation}
\ell_{\rho_{\mu, \mathcal{X}}^w}(x)=i(\phi_*(\mu), \phi^*(x)) = i(\mu, x), \quad \mbox{ for all } x \in \Gamma,
\label{eq:ell_equal_intersection}
\end{equation}
and hence the rough isometry class of $\rho_{\mu, \mathcal{X}}^w$ only depends on $\mu$, and not on the choice of $\mathcal{X} \in \operatorname{Teich}(\Gamma)$. We will just refer to its rough similarity class as $\rho_{\mu}$.
In \cite{LucDid}, the second author studies in detail the properties of these pseudo-metrics.
They prove
\begin{enumerate}
\item $\rho_{\mu}$ is $\delta$-hyperbolic, and the $\delta$ hyperbolicity property is given in terms of $\mu$.
\item $\rho_{\mu}$ is (equivariant isometric to an) $\mathbb{R}$-tree if and only if $\mu$ is a measured lamination.
\item $\rho_{\mu} \in D_\Gamma$ if and only if $\mu$ is filling.
\end{enumerate}

The last item is also proven by Cantrell-Reyes \cite{StephenEduardo}.

\subsection{Dilation and entropy}

The following is an extension of the dilation from that of filling geodesic currents.

\begin{definition}
For two $\rho,\rho'\in D_\Gamma$, we define its \emph{dilation} by the formula
\[
\Dil(\rho,\rho') \coloneqq \sup_{[x] \in [\Gamma']} \frac{\ell_{\rho'}(x)}{\ell_\rho(x)}
\]
\end{definition}
where $[\Gamma']$ denotes the subset of conjugacy classes $[x]$, for $x \in \Gamma$, so that $\ell_\rho(x)>0$ (and, hence, $\ell_{\rho'}(x)>0$).

It is easy to see that for any triple $\rho,\rho', \rho'' \in D_\Gamma$, $\Dil(\rho,\rho'') \leq \Dil(\rho,\rho') \Dil(\rho',\rho'')$.

Another relevant quantity is that of \emph{entropy} of $\rho \in D_\Gamma$.
\[
h_\rho \coloneqq \limsup_n \log|\{ [x]\in \Gamma : \ell_\rho(x) \leq n\}|/n,
\]
which can be proven (see Cantrell-Tanaka) to coincide with the quantity
\[
v_\rho \coloneqq \lim_{r \to \infty} \frac{1}{r} \log \# B(o,r),
\]
where $B(o,r)\coloneqq \{ y \in \Gamma : d(y,o) < r \}$.
In~\cite[Proposition~3.1]{CT25:Manhattan} it is shown that, in fact,
\[
h_\rho=v_\rho.
\]
\subsection{Manhattan curves}
\label{subsec:manhattan}
Manhattan curves were first introduced in the context of isometric actions on rank-1 symmetric spaces~\cite{Bur93:Manhattan}.
We follow Cantrell-Tanaka.

\begin{definition}
Let $\rho, \rho' \in D_\Gamma$. The Manhattan curve of $\rho,\rho'$ is the boundary of the convex set
\[
C_{\rho'/\rho} \coloneqq \{ (a,b) \in \mathbb{R}^2 : \sum_{x \in \Gamma} e^{-a\rho'(o,x) - b\rho(o,x)} < \infty \}.
\]
Note that if $a=0$, the sum is the Poincar\'e sum of $\rho$ with parameter $b$.
Hence, the critical exponent of the resulting series is the infimum among $b$ for which the series converges.
Similarly, for each $t \in \mathbb{R}$, let $\theta_{\rho'/\rho}(t)$ denote the critical exponent of the series
\[
s \mapsto \sum_{x \in \Gamma} e^{-t \rho'(o,x) -\rho d(o,x)}.
\]
\end{definition}
As happens for the entropy, the Manhattan curve can be equivalently defined using stable lenghts instead of distances, and counting over conjugacy classes, as follows.
Let $C_M$ the subset of $\mathbb{R}^2$ given by
\[
\mathrm{C}_{\rho'/\rho}^\ell \coloneqq \{(a,b) \in \mathbb{R}^2 : Q(a,b) < \infty \}
\]
where 
\[
Q(a,b) \coloneqq \sum_{[x] \in [\Gamma]} \exp(-a \ell_{\rho}(x) -b\ell_{\rho'}(x)).
\]
Then $\mathrm{C}_{\rho'/\rho}^\ell=C_{\rho'/\rho}$. 
In the sequel, we will simply refer to it as $C_{\rho'/\rho}$ or $C_M$ when the specific dependence on $\rho,\rho'$ is not emphasized.
Cantrell-Tanaka prove the following.
\begin{theorem} 
For $\rho, \rho' \in D_\Gamma$.
\begin{enumerate}
\item $\theta_{\rho'/\rho}$ is convex, decreasing, and continuously differentiable.
\item $\theta_{\rho'/\rho}$ goes through the points $(0,v_\rho)$ and $(v_{\rho'},0)$, and is a straight line between these points if and only if $\rho,\rho'$ are roughly similar.
\item \[
-\theta'_{\rho'/\rho}(v_{\rho'}) \leq \frac{v_\rho}{v_{\rho'}} \leq -\theta'_{\rho'/\rho}(0),
\]
and both equalities occur if and only if $\rho$ and $\rho'$ are roughly similar.
\item $\rho,\rho'$ are roughly similar if and only if they have proportional marked length spectra.
\label{thm:cantrell-tanaka_main}
\end{enumerate}
\end{theorem}

For certain pairs of pseudo-metrics, such as word-metrics, or Green metrics, $C_{\rho'/\rho}$ is known to be real analytic.
Cantrell-Reyes also prove that for $\rho,\rho'$ not roughly similar, $\theta_{\rho'/\rho}$ is in fact strictly convex~\cite[Proposition~2.3]{CR24:Marked}.

Cantrell-Tanaka prove in \cite[Theorem 1.2]{CT25:Manhattan} the following

\begin{theorem}
For any pair $(\rho,\rho')$ with $\rho,\rho' \in D_\Gamma$, the following limit, called the \emph{mean distortion of $\rho$ and $\rho'$ } exists:
\[
\tau(\rho'/\rho) \coloneqq \lim_{r \to \infty} \frac{1}{\# B(o,r)} \sum_{x \in B(o,r)} \frac{\rho'(0,x)}{r},
\]
where the balls $B(o,r) \coloneqq \{ y \in \Gamma | \rho(y,o) \leq r\}$.
Moreover, 
\[
\tau(\rho'/\rho) \geq \frac{v}{v'}
\]
where $v$ is the exponential growth of $\rho$ and $v'$, is the exponential growth of $\rho'$.
Moreover, the following conditions are equivalent:
\begin{enumerate}
\item The equality $\tau(\rho'/\rho) = v/v'$ holds,
\item there exist a constant $c>0$ so that $\ell'([x])=c \ell([x])$ for all $[x] \in [\Gamma]$
\item $\rho$ and $\rho'$ are roughly similar.
\end{enumerate}
\label{thm:ct22}
\end{theorem}

As proven in Cantrell-Tanaka, it can be shown that the mean distortion can be equivalently defined in terms of stable lengths and conjugacy classes with non-trivial stable length, as
\[
\tau(\rho'/\rho) \coloneqq \lim_{r \to \infty} \frac{1}{|\{ [x] \in [\Gamma'] : \ell_\rho(x) \leq r \}|} \sum_{\ell_\rho(x) \leq r} \frac{\ell_{\rho'}(x)}{r}.
\]

The mean distortion is not symmetric in $\rho,\rho'$, in general. It is well-defined, finite and positive.

\subsection{Hyperbolic potentials}

In recent work~\cite{CRS24:Joint} it is observed that some of the avatars of Cantrell-Tanaka's paper~\cite{CT25:Manhattan} related to the Patterson-Sullivan construction can be promoted from $\mathcal{D}_\Gamma$ to the more general space of \emph{metric potentials} $\mathcal{H}^{++}_\Gamma$, introduced by them, which we will now survey. These are, roughly speaking, a generalization of metric structures  which are also allowed to be asymmetric and take negative values, while their stable lengths are still positive.

We will discuss how the asymmetric metric $d$ on $\mathcal{D}_\Gamma$ can be further extended to $\mathcal{H}^{++}_\Gamma$.

Following~\cite{CT25:Manhattan,CRS24:Joint}, we say $\psi$ is a \emph{hyperbolic metric potential} if it is a $\Gamma$ left-invariant function $\psi \colon \Gamma \times \Gamma \to \mathbb{R}$ and for every $\rho_0 \in D_\Gamma$ there exists a constant $C>0$ so that
\[
(x|y)_z^\psi < C(x|y)_z^{\rho_0}+C,
\]
where 
\[
(x|y)_z^\psi \coloneqq \frac{1}{2}\left( \psi(x,z) + \psi(z,y) - \psi(x,y)\right).
\]
The space of all hyperbolic metric potentials is denoted by $H_\Gamma$. 
We can similarly quotient by the equivalence relation of \emph{rough similarity}, i.e., two hyperbolic metric potentials $\psi, \psi'$ are equivalent iff there exist constants $k,A>0$ so that
\[
|k\cdot \psi - \psi'|<A.
\]
Finally, we say $\psi$ is a proper hyperbolic metric potential, and denote the space of such elements by $H^{++}_\Gamma$, iff there exist $\rho_0 \in D_\Gamma$ so that $\ell_\psi \geq \ell_{\rho_0}$. Since this is well defined up to rough similarity, we can coin the next definition.

\begin{definition}
\emph{The space of metric potentials} is defined as the space of rough similarity classes of proper hyperbolic metric potentials, and is denoted by $\mathcal{H}_\Gamma^{++}$.
\end{definition}

Following~\cite{CT25:Manhattan}, we say $\psi$ satisfies (QE) if there exists a \emph{quasi-extension} of $\psi$ to 
$\Gamma \times (\Gamma \cup \partial \Gamma)$, i.e.,
there exists a function $\psi(\cdot|\cdot)_0 \colon \Gamma \times (\Gamma \cup \partial \Gamma) \to \mathbb{R}$
and a constant $C \geq 0$, so that
\[
\limsup_{n \to \infty} \psi(x|\zeta_n) - C \leq \psi(x|\zeta)_o \leq \liminf_{n \to \infty} \psi(x|\zeta_n')_o + C
\]
for all $(x,\zeta) \in \Gamma \times (\Gamma \cup \partial \Gamma)$ and for all 
sequences $(\zeta_n),(\zeta_n') \in \zeta$.
One can then define the analogous Busemann function $(x,\zeta) \in \Gamma \times \partial \Gamma$,
\[
\beta_o^\psi(x,\zeta) \coloneqq \sup \{ \limsup_n (\psi(x,\zeta_n)-\psi(o,\zeta_n)) : (z_n) \in \zeta \}.
\]
If $\psi$ is $\Gamma$-invariant, then $\beta_0^\psi$ satisfies
\[
|\beta_0^\psi(xy,\zeta)-(\beta_0^\psi(y,x^{-1}\zeta) + \beta_o^\psi(x,\zeta))| \leq 4C.
\]
We will now extend Cantrell-Tanaka's~definition of (RG) (see~\cite[Section~2.6]{CT25:Manhattan} for pairs $(\psi,d)$, with $d \in D_\Gamma$
and $\psi \in H^{++}_\Gamma$ to pairs $(\psi',\psi)$
where $\psi \in H_\Gamma^{++}$ and $\psi' \in H_\Gamma$.
We say $\psi'$ satisfies (RG) relative to $\psi$ if there exists $C>0$ so that
\[
(x|y)_z^{\psi'} \leq C (x|y)_z^{\psi} + C
\]
for all $x,y,z \in \Gamma$.

We note that when $\psi=d \in D_\Gamma$, the latter definition is equivalent to Cantrell-Tanaka's definition for pairs $(d, \psi)$, $d \in D_\Gamma$, $\psi \in H^{++}_\Gamma$, by~\cite[Proposition~3.1]{CRS24:Joint}.

In particular, we will observe here that, in fact, Theorem~\ref{thm:ct22} can be promoted to $H^{++}_\Gamma$, as follows.

\begin{theorem}
For any pair $(\psi,\psi')$ with $\psi,\psi' \in H^{++}_\Gamma$, the following limit, called the \emph{mean distortion of $\psi$ and $\psi'$ } exists:
\[
\tau(\psi'/\psi) \coloneqq \lim_{r \to \infty} \frac{1}{\# C^\psi(r)} \sum_{x \in C^\psi(r)} \frac{\ell_{\psi'}([x])}{r},
\]
where the balls $C^\psi(r) \coloneqq \{ [y]] \in [\Gamma'] \mid \ell_\psi([y]) \leq r\}$.
Moreover, 
\[
\tau(\psi'/\psi) \geq \frac{h(\psi)}{h(\psi')}
\]
where $h(\psi)$ (resp. $h(\psi')$) is the entropy of $\psi$ (resp. $\psi'$).
Moreover, the following conditions are equivalent:
\begin{enumerate}
\item The equality $\tau(\psi'/\psi) = h(\psi)/h(\psi')$ holds,
\item there exist a constant $c>0$ so that $\ell_{\psi'}([x])=c \ell_\psi([x])$ for all $[x] \in [\Gamma]$
\item $\psi$ and $\psi'$ are roughly similar.
\end{enumerate}
\label{thm:ct22_potentials}
\end{theorem}

First, we observe the following extension of \cite[Proposition~2.7]{CT25:Manhattan}
for pairs $(\psi, \psi') \in H_\Gamma^{++} \times H_\Gamma$, with $\psi'$ satisfying (RG) with respect to $\psi$.
The proof is the same as in that proposition, replacing $\psi$ by $\psi'$ and $d$ by $\psi$, in their notation, and defining
\[
O^\psi(x,R) \coloneqq \{ \zeta \in \partial \Gamma : (\zeta|o)_x^\psi \leq R \}.
\]

\begin{proposition}
For $\psi \in H^{++}_\Gamma$, let $\psi'\in H_\Gamma$ satisfy (RG) relative to $\psi$.
Then the abscissa of convergence $\theta$ of the series in $s$ given by,
\[
\sum_{x \in \Gamma} \exp(-\psi'(o,x)-s\psi(o,x)),
\]
is finite and there exists a probability measure $\mu_\psi$ on $\partial \Gamma$
satisfying (QC) with exponent $\theta$ relative to $(\psi', \psi)$.
Moreover, every finite Borel measure $\mu$ satisfying (QC) has the property:
\[
(C')^{-1} \exp(-\psi'(o,x) - \theta \psi(o,x)) \leq \mu(O^\psi(x,R)) \leq C' 
exp(-\psi'(o,x)-\theta \psi(o,x))
\]
for all $x \in \Gamma$, where $C'>0$ is a constant depending on $R$ and the constant $C$ in the definition of (RG).
\label{prop:CT_potentials}
\end{proposition}

\begin{proof}[Proof of Theorem~\ref{thm:ct22_potentials}]
Proposition~\ref{prop:CT_potentials} allows us to extend Lemma 2.8 and 2.9 in~\cite{CT25:Manhattan} to pairs $(\psi, \psi')$ by, again, performing the same replacement.
These results can then be used to prove the analog of Corollary 2.10., Proposition 2.11, Lemma 2.12 and Corollary 2.13 in~\cite{CT25:Manhattan}.
This, in turn, is used to settle the extension of Lemma 3.5 and Lemma 3.6 in~\cite{CT25:Manhattan}.
Lemma~3.5 is used in the proof of the extension of Theorem 3.12, which we rely on.
Both Lemma 3.5 and Lemma 3.6 are used in proving the extension of~\cite[Theorem~3.7]{CT25:Manhattan}.
Corollary 2.10, Lemma 2.12, Lemma 3.5 and Theorem 3.7 are then used to prove the analog of~\cite[Theorem~3.10]{CT25:Manhattan}. 
Finally, the extensions of Theorem 3.12, Lemma 3.4, Theorem 3.7 and Theorem 3.10 settle the analog of Theorem 1.2.

Finally, we use the analogous extensions of $\tau$ in terms of conjugacy classes and stable lengths, as well as the corresponding extension of Manhattan curves which follows from the extension of~\cite[Proposition~3.1]{CT25:Manhattan}.
\end{proof}

We then have the following result, which follows from similar ideas as in~\cite[Proposition~3.1]{CT25:Manhattan}, where the dilation of $\psi'$ and $\psi$, $\Dil(\psi',\psi)$ is defined in the analogous way.
\begin{lemma}
For any pair of pseudo-metric $\psi,\psi' \in H^{++}_\Gamma$, we have
\[
\tau(\psi'/\psi) \leq \Dil(\psi',\psi).
\]
\label{lem:mdvsdil}
\end{lemma}
\begin{proof}
For large enough $L >0$, and for each $a \in \mathbb{R}$
\[
\theta(a) \coloneqq \limsup_{n \to \infty} \frac{1}{n} 
\log \sum_{\substack{x \in \Gamma \\ |\psi(o,x) - n| \leq L}} 
  e^{-a \psi'(o,x)}
\]
and
\[
\Theta(a) \coloneqq \limsup_{n \to \infty} \frac{1}{n} 
\log \sum_{\substack{[x] \in [\Gamma'] \\ |\ell_{\psi}(x) - n| \leq L}} 
  e^{-a \ell_{\psi'}(x)}.
\]
For every conjugacy class $[x] \in [\Gamma']$, we have
\[
\Dil(\psi',\psi)  \geq \frac{\ell_{\psi'}(x)}{\ell_{\psi}(x)}.
\]
In particular, for every $a > 0$ and every $n \geq 0$, we have
\[
\sum_{[x] : |\ell_{\psi}(x) - n| \leq L} e^{-a \ell_{\psi'}(x)} \geq c \cdot e^{-a \cdot \Dil \cdot n}.
\]
By the extension of \cite[Proposition 3.1]{CT25:Manhattan}, $\theta(a)=\Theta(a)$, and therefore, by 
taking logarithms, dividing by $n$, and taking $\limsup$ when $n\to \infty$, we get
\[
\theta(a) \geq -a \cdot \Dil(\psi',\psi).
\]
By the extension of \cite[Theorem 3.7.]{CT25:Manhattan}, $\theta$ is differentiable on $\mathbb{R}$ and $\theta'(a)=\tau_{a,b}$, where
$\tau_{a,b}$ is a function of the pairs of points $(a,b) \in C_M$.
Hence,
\[
-\frac{\theta(a)}{a} \leq \Dil(\psi',\psi)
\]
and taking $\lim a \to {0^+}$ and using differentiability of $\theta$,
\[
-\theta'(0)\leq \Dil(\psi',\psi).
\]
By the extension of \cite[Theorem 3.12]{CT25:Manhattan}, moreover, $\tau_{0,v} = \tau(\psi'/\psi)$.
Therefore, $-\theta'(0)=\tau(\psi'/\psi)$, and the result follows.
\end{proof}

We will use Lemma~\ref{lem:mdvsdil} to establish the following proposition. 

\begin{proposition}
The function
\[
d(\psi,\psi') \coloneqq \log \left( \Dil(\psi,\psi') \frac{h(\psi')}{h(\psi)} \right).
\]
is an asymmetric metric on $\mathcal{H}^{++}_\Gamma$, i.e., $d$ is non-negative, non-degenerate and satisfies the triangular inequality.
\label{prop:asymmetric_defined}
\end{proposition}
\begin{proof}
By Theorem \ref{thm:ct22_potentials} and Lemma \ref{lem:mdvsdil} it follows that
\[
h(\psi)/h(\psi') \leq \tau(\psi'/\psi) \leq \Dil_+(\psi',\psi).
\]
Thus, $d(\psi,\psi') \geq 0$.
Also, by this inequalities, it follows that if $d(\psi,\psi')=0$, then $h(\psi)/h(\psi') = \Dil_+(\psi',\psi)$, and hence $h(\psi)/h(\psi') = \tau(\psi'/\psi)$, so by Theorem \ref{thm:ct22_potentials}
$\psi'$ and $\psi$ are roughly similar, proving non-degeneracy of $d$.
Triangular inequality is a straightforward computation using the fact that the supremum of a sum is smaller or equal than the sum of suprema and the monotonicity of the logarithm.
\end{proof}

\subsection{Cubulations}\label{ex:cubulation}

This asymmetric distance $d$ on the space of metric potentials $\mathcal{H}^{++}_\Gamma$ generalizes the asymmetric distance on the space of metric structures $\mathcal{D}_\Gamma$. It also generalizes the asymmetric distance considered in many other contexts.
We have already seen some of those contexts in Section~\ref{sec:examples}, including many arising from geodesic currents, as well as, more generally, lengths arising from Anosov representations. Another example, less central to this paper, is that of \emph{cubulations} of $\Gamma$, i.e., geometric actions of a Gromov hyperbolic group $\Gamma$ a essential, hyperplane essential $\operatorname{CAT}(0)$-cube complex $\mathcal{X}$ of entropy 1. More generally, one can consider \emph{cuboids}~\cite{BeyrerFioravanti2021} of entropy 1, where a cuboid is a generalization of a cubulation where edge distances are allowed to be arbitrary positive real numbers as opposed to just positive integers. We refer the reader to the cited paper for the definitions of ``essential'' and ``hyperplane essential'', which are beyond the scope of this paper. This includes many examples of cubulations coming naturally from closed curves on surfaces, or invariant collections of quasi-circles associated to quasi-fuchsian 3-manifolds~\cite{BR24:Approximating}.
Marked length spectrum rigidity for geometric actions on cuboids is proven in~\cite{BeyrerFioravanti2021}.

\subsection{Geodesic currents on surface groups}

In fact, it induces an asymmetric metric on projective filling geodesic currents.

\begin{lemma}
The map $\mu \mapsto \rho_\mu$ from $\mathcal{C}_{fill}$ to $D_{\pi_1(S)}$ induces an embedding from $\PC_{\mbox{fill}}$ to $\mathcal{D}_{\pi_1(S)}$.
\label{lem:currasym}
\end{lemma}
\begin{proof}
The proof is not new (see \cite[Theorem~1.11]{StephenEduardo},~\cite[Theorem~C]{LucDid}), but we sketch it here for completeness. A geodesic current $\mu$ is filling if and only if $\rho_{\mu}$ is proper. Also, the action of $\Gamma$ on $\rho_{\mu}$ is cobounded, hence, by Milnor-Svar\c{c} lemma, $\rho_{\mu}$ is in $D_\Gamma$.
Hence, to show that the map is an embedding it suffices to show that two pseudo-metrics $\rho=\rho_\mu, \rho'=\rho_{\mu'}$ arising from filling geodesic currents $\mu, \mu'$ are roughly equivalent if and only if the currents are multiples of one another. Suppose that there are constants $k, A>0$ so that
\[
|\rho- k \cdot \rho'| < A
\]
taking $g \in \pi_1(S)$, let $\tilde{g}$ be the hyperbolic axis of the deck transformation associated to $g$
and $x \in \tilde{g} \subset \mathbb{H}^2$. We have
\[
-A < \rho_\mu(x, g^n x) - k \cdot \rho_{\mu'}(x, g^n x) < A.
\]
From the choice of $x$, it follows that 
$ \rho_\mu(x, g^n x) = i(\mu, g^n)$
and $ \rho_{\mu'}(x, g^n x)= i(\mu', g^n)$, for every $g\in \pi_1(S)$, and for every $n \in \mathbb{N}$.
Dividing the previous inequalities by $n$ and taking the limit when $n \to \infty$, we get
$i(\mu, g) = i(k \cdot \mu', g)$ for every $g \in \pi_1(S)$.
By \cite[Th\'eor\`eme 2]{Otal90:SpectreMarqueNegative}, we have $\mu= k \cdot \mu'$.
The other implication is clear.
\end{proof}

From here, we have

\begin{proposition}
The restriction of $d$ to (the  image of) $\PC_{fill}$ inside of $\mathcal{H}^{++}_{\pi_1(S)}$ induces an asymmetric distance $d_{\PC_{fill}}$ on $\PC_{fill}$ given by
\[
d_{\PC_{fill}}(x,y) = \log \left( \sup_{\bc \in \cc} \frac{i(x,\bc)}{i(y,\bc)} \cdot \frac{h(x)}{h(y)} \right).
\]
\label{prop:asymmetric_geodesics}
\end{proposition}
\begin{proof}
We have
\[
d_{\PC_{fill}}(x,y) = \log \left( \sup_{g \in \Gamma} \frac{\ell_{\rho_y}(g)}{\ell_{\rho_x}(g)} \cdot \frac{h(\rho_y)}{h(\rho_x)} \right)=\log \left( \sup_{\bc \in \cc} \frac{i(y,\bc)}{i(x,\bc)} \cdot \frac{h(y)}{h(x)} \right)
\]
where the first equality is by definition of the metric $d$, and the second by Equation~\ref{eq:ell_equal_intersection} we have $i(x,\bc)=\ell_{\rho_x}(\bc)$, and hence also $h(x) = h(\rho_x)$.
\end{proof}

\subsection{The asymmetric metric on $\mathcal{D}_\Gamma$ is geodesic}

Given two pseudo-metrics $\rho,\rho'$, let $\theta_{\rho'/\rho}(t)$ be their Manhattan curve.

The next result is contained in~\cite[Theorem~1.2]{StephenEduardo}, and says, in particular, that there are unparameterized geodesic segments between any two points in $\mathcal{D}_\Gamma$.

\begin{proposition}
    Given two pseudo-metrics $\rho,\rho' \in D(\Gamma)$. For $t \in \mathbb{R}$ let $\rho_t \coloneqq t \rho + \theta(t) \rho'$.  Then, for every $r < s < t$, we have 
\[d_{\mathrm{sym}}([\rho_r],[\rho_t]) = d_{\mathrm{sym}}([\rho_r],[\rho_s]) + d_{\mathrm{sym}}([\rho_s],[\rho_t]).\]
\label{prop:geod_sym}
\end{proposition}

The following is a general lemma stating that geodesics for the symmetrized metric are geodesics for the asymmetric metric.

\begin{proposition}
Let $(X,d)$ be a metric space equipped with an asymmetric metric $d$, and let $d_{\mathrm{sym}}(x,y) \coloneqq d(x,y)+d(y,x)$ for every $x,y \in X$. If $\gamma$ is a $d_{\mathrm{sym}}$ geodesic between $x, y \in X$, then $\gamma$ is a $d$ geodesic between $x$ and $y$, and a $d$ geodesic between $y$ and $x$.
\label{prop:dsym_implies_d_geodesic}
\end{proposition}

\begin{proof}
Suppose $\gamma \colon [0,1] \to X$ is a geodesic wrt $d_{\mathrm{sym}}$, i.e., for every $0 \leq s < t < u \leq 1$ we have
\[
d_{\mathrm{sym}}(\gamma(s),\gamma(u)) = d_{\mathrm{sym}}(\gamma(s),\gamma(t)) + d_{\mathrm{sym}}(\gamma(t),\gamma(u)).
\]
Expanding in the definition of $d_{\mathrm{sym}}$ and rewriting terms on one side, we get
\[
\left[d(\gamma(s),\gamma(u)) - d(\gamma(s),\gamma(t)) - d(\gamma(t),\gamma(u)) \right] + \left[d(\gamma(u),\gamma(s)) - d(\gamma(u),\gamma(t)) - d(\gamma(t),\gamma(s)) \right]=0.
\]
Since both bracketed terms are $\leq 0$ by triangular inequality of $d$, it follows that both of them have to be zero. Hence, the claim follows.
\end{proof}

\bibliographystyle{alpha} 
\bibliography{horoboundary}

\begin{thebibliography}{GGKW17}

\bibitem[AB06]{AB06:InfiniteDimensional}
Charalambos~D. Aliprantis and Kim~C. Border.
\newblock {\em Infinite Dimensional Analysis}.
\newblock Springer-Verlag Berlin Heidelberg, 2006.

\bibitem[Aze24]{A24:HoroboundaryQualitative}
Aitor Azemar.
\newblock A qualitative description of the horoboundary of the {T}eichm\"uller metric.
\newblock {\em Algebr. Geom. Topol.}, 24(7):3919--3984, 2024.

\bibitem[Bal95]{Nonpositive}
Werner Ballmann.
\newblock {\em Lectures on spaces of nonpositive curvature}, volume~25 of {\em DMV Seminar}.
\newblock Birkh\"{a}user Verlag, Basel, 1995.
\newblock With an appendix by Misha Brin.

\bibitem[BCLS18]{BridgemanCanaryLabourieSambarino18:SimpleRoots}
Martin Bridgeman, Richard Canary, Fran\c{c}ois Labourie, and Andres Sambarino.
\newblock Simple root flows for {H}itchin representations.
\newblock {\em Geom. Dedicata}, 192:57--86, 2018.

\bibitem[Ben97]{Benoist_AsymtoticLinearGroups}
Yves Benoist.
\newblock Propri\'{e}t\'{e}s asymptotiques des groupes lin\'{e}aires.
\newblock {\em Geom. Funct. Anal.}, 7(1):1--47, 1997.

\bibitem[BF21]{BeyrerFioravanti2021}
Jonas Beyrer and Elia Fioravanti.
\newblock Cross-ratios on cat(0) cube complexes and marked length-spectrum rigidity.
\newblock {\em Journal of the London Mathematical Society}, 104(5):1973--2015, 2021.

\bibitem[BHM11]{BHM2011}
Sébastien Blachère, Peter Haïssinsky, and Pierre Mathieu.
\newblock Harmonic measures versus quasiconformal measures for hyperbolic groups.
\newblock {\em Annales scientifiques de l’École Normale Supérieure, Série 4}, 44(4):683--721, 2011.

\bibitem[BIPP21]{BIPP21}
M.~Burger, A.~Iozzi, A.~Parreau, and M.~B. Pozzetti.
\newblock Currents, systoles, and compactifications of character varieties.
\newblock {\em Proc. Lond. Math. Soc. (3)}, 123(6):565--596, 2021.

\bibitem[BIPP24]{BIPP24}
Marc Burger, Alessandra Iozzi, Anne Parreau, and Maria~Beatrice Pozzetti.
\newblock Positive crossratios, barycenters, trees and applications to maximal representations.
\newblock {\em Groups Geom. Dyn.}, 18(3):799--847, 2024.

\bibitem[BL17]{BL17:RigidityFlat}
Anja Bankovic and Christopher~J. Leininger.
\newblock Marked-length-spectral rigidity for flat metrics.
\newblock {\em Transactions of the American Mathematical Society}, 370(3):1867–1884, 2017.

\bibitem[Bon86]{Bon86}
Francis Bonahon.
\newblock Bouts des vari\'{e}t\'{e}s hyperboliques de dimension {$3$}.
\newblock {\em Ann. of Math. (2)}, 124(1):71--158, 1986.

\bibitem[Bon88]{Bon88}
Francis Bonahon.
\newblock The geometry of {T}eichm\"{u}ller space via geodesic currents.
\newblock {\em Invent. Math.}, 92(1):139--162, 1988.

\bibitem[BPS19]{BPS}
Jairo Bochi, Rafael Potrie, and Andr\'es Sambarino.
\newblock Anosov representations and dominated splittings.
\newblock {\em Jour. Europ. Math. Soc.}, 11:3343--3414, 2019.

\bibitem[BR24]{BR24:Approximating}
Nic Brody and Eduardo Reyes.
\newblock Approximating hyperbolic lattices by cubulations, 2024.

\bibitem[Bur93]{Bur93:Manhattan}
Marc Burger.
\newblock Intersection, the manhattan curve, and patterson--sullivan theory in rank 2.
\newblock {\em International Mathematics Research Notices}, 1993(7):217--225, 1993.

\bibitem[CDPW24]{CDPW24:Asymmetric}
Le{\'o}n Carvajales, Xian Dai, Beatrice Pozzetti, and Anna Wienhard.
\newblock Thurston’s asymmetric metrics for anosov representations.
\newblock {\em Groups, Geometry, and Dynamics}, 2024.
\newblock Published online first on July 24, 2024.

\bibitem[CFF92]{CFF92:RigidityNonPosCurvedRiem}
C.~Croke, A.~Fathi, and J.~Feldman.
\newblock The marked length-spectrum of a surface of nonpositive curvature.
\newblock {\em Topology}, 31(4):847--855, 1992.

\bibitem[CMGR]{CMGR26:GreenMetrics}
Stephen Cantrell, D{\'\i}dac Mart{\'\i}nez-Granado, and Eduardo Reyes.
\newblock Density of green metric for hyperbolic groups.
\newblock in preparation.

\bibitem[CR24]{CR24:Marked}
Stephen Cantrell and Eduardo Reyes.
\newblock Marked length spectrum rigidity from rigidity on subsets, 2024.

\bibitem[CR25]{StephenEduardo}
Stephen Cantrell and Eduardo Reyes.
\newblock Manhattan geodesics and the boundary of the space of metric structures on hyperbolic groups.
\newblock {\em Comment. Math. Helv.}, 100(1):11--59, 2025.

\bibitem[CRS24]{CRS24:Joint}
Stephen Cantrell, Eduardo Reyes, and Cagri Sert.
\newblock The joint translation spectrum and manhattan manifolds, 2024.

\bibitem[CT24]{CT24:Invariant}
Stephen Cantrell and Ryokichi Tanaka.
\newblock Invariant measures of the topological flow and measures at infinity on hyperbolic groups.
\newblock {\em J. Mod. Dyn.}, 20:215--274, 2024.

\bibitem[CT25]{CT25:Manhattan}
Stephen Cantrell and Ryokichi Tanaka.
\newblock The {M}anhattan curve, ergodic theory of topological flows and rigidity.
\newblock {\em Geom. Topol.}, 29(4):1851--1907, 2025.

\bibitem[DK00]{DK00:Conjugacy}
Fran\c{c}oise Dal'Bo and Inkang Kim.
\newblock A criterion of conjugacy for zariski dense subgroups.
\newblock {\em Comptes Rendus de l'Acad{\'e}mie des Sciences - Series I - Mathematics}, 330(8):647--650, 2000.

\bibitem[DK22]{DK22:PSAnosov}
Subhadip Dey and Michael Kapovich.
\newblock Patterson-{S}ullivan theory for {A}nosov subgroups.
\newblock {\em Trans. Amer. Math. Soc.}, 375(12):8687--8737, 2022.

\bibitem[DLR10]{DLR2010}
Moon Duchin, Christopher~J. Leininger, and Kasra Rafi.
\newblock Length spectra and degeneration of flat metrics.
\newblock {\em Invent. Math.}, 182(2):231--277, 2010.

\bibitem[Dow91]{Dow91:Choquet}
Tomasz Downarowicz.
\newblock The {C}hoquet simplex of invariant measures for minimal flows.
\newblock {\em Israel J. Math.}, 74(2-3):241--256, 1991.

\bibitem[DPS12]{DPS12:HoroboundaryNegative}
Fran\c~coise Dal'bo, Marc Peign\'e, and Andrea Sambusetti.
\newblock On the horoboundary and the geometry of rays of negatively curved manifolds.
\newblock {\em Pacific J. Math.}, 259(1):55--100, 2012.

\bibitem[DRMG25]{LucDid}
Luca De~Rosa and Dídac Martínez-Granado.
\newblock Dual spaces of geodesic currents.
\newblock {\em Journal of Topology}, 2025.
\newblock to appear.

\bibitem[ELS22]{ELS22:HyperbolicCone}
Viveka Erlandsson, Christopher~J. Leininger, and Chandrika Sadanand.
\newblock Hyperbolic cone metrics and billiards.
\newblock {\em Adv. Math.}, 409:Paper No. 108662, 74, 2022.

\bibitem[ES22]{ES22:GeodesicCount}
Viveka Erlandsson and Juan Souto.
\newblock {\em Geodesic currents and {M}irzakhani's curve counting}, volume xii-151 of {\em Progress in Mathematics}.
\newblock Birkhäuser Cham, 2022.

\bibitem[FF22]{FF22:QF}
Ethan Fricker and Alex Furman.
\newblock Quasi-fuchsian vs negative curvature metrics on surface groups.
\newblock {\em Israel Journal of Mathematics}, 251:365--378, 2022.

\bibitem[FK65]{FK65:HyperbolicSpectrum}
Robert Fricke and Felix Klein.
\newblock {\em Vorlesungen \"uber die {T}heorie der automorphen {F}unktionen. {B}and 1: {D}ie gruppentheoretischen {G}rundlagen. {B}and {II}: {D}ie funktionentheoretischen {A}usf\"uhrungen und die {A}ndwendungen}, volume Bande 3, 4 of {\em Bibliotheca Mathematica Teubneriana}.
\newblock Johnson Reprint Corp., New York; B. G. Teubner Verlagsgesellschaft, Stuttgart, 1965.

\bibitem[FM11]{FrancavigliaMartino2011}
Stefano Francaviglia and Armando Martino.
\newblock Metric properties of outer space.
\newblock {\em Publicacions Matemàtiques}, 55(2):433--473, 2011.

\bibitem[Fur02]{Fur02:Coarse}
Alex Furman.
\newblock {\em Coarse-Geometric Perspective on Negatively Curved Manifolds and Groups}, pages 149--166.
\newblock Springer Berlin Heidelberg, Berlin, Heidelberg, 2002.

\bibitem[GdlH90]{GdlH88:Notes}
\'E. Ghys and P.~de~la Harpe, editors.
\newblock {\em Sur les groupes hyperboliques d'apr\`es {M}ikhael {G}romov}, volume~83 of {\em Progress in Mathematics}.
\newblock Birkh\"auser Boston, Inc., Boston, MA, 1990.
\newblock Papers from the Swiss Seminar on Hyperbolic Groups held in Bern, 1988.

\bibitem[GGKW17]{GGKW}
Olivier Guichard, Fran\c{c}ois Guéritaud, Fanny Kassel, and Anna Wienhard.
\newblock Anosov representations and proper actions.
\newblock {\em Geom. Topol.}, 21:485--584, 2017.

\bibitem[GK17]{GK11:GeometricallyFinite}
Fran\c~cois Gu\'eritaud and Fanny Kassel.
\newblock Maximally stretched laminations on geometrically finite hyperbolic manifolds.
\newblock {\em Geom. Topol.}, 21(2):693--840, 2017.

\bibitem[GM25]{GM25:MargulisSpacetimes}
Krishnendu Gongopadhyay and Neelanjan Mondal.
\newblock Thurston's asymmetric metric on margulis spacetimes, 2025.

\bibitem[Gro81]{Gro}
M.~Gromov.
\newblock Hyperbolic manifolds, groups and actions.
\newblock In {\em Riemann surfaces and related topics: {P}roceedings of the 1978 {S}tony {B}rook {C}onference ({S}tate {U}niv. {N}ew {Y}ork, {S}tony {B}rook, {N}.{Y}., 1978)}, Ann. of Math. Stud., No. 97, pages 183--213. Princeton Univ. Press, Princeton, NJ, 1981.

\bibitem[GW12]{GW}
Olivier Guichard and Anna Wienhard.
\newblock Anosov representations: domains of discontinuity and applications.
\newblock {\em Invent. Math.}, 190:357--438, 2012.

\bibitem[HP97]{HP97:RigidityNegCurvedCone}
Sa'ar Hersonsky and Fr\'{e}d\'{e}ric Paulin.
\newblock On the rigidity of discrete isometry groups of negatively curved spaces.
\newblock {\em Comment. Math. Helv.}, 72(3):349--388, 1997.

\bibitem[Hub61]{entropy1}
H.~Huber.
\newblock Zur analytischen theorie hyperbolischer raumformen und bewegungsgruppen. ii. nachtrag zu math. annalen 142, 385-398 (1961).
\newblock {\em Mathematische Annalen}, 143:463--464, 1961.

\bibitem[KL18]{KL18:Finsler}
Michael Kapovich and Bernhard Leeb.
\newblock Finsler bordifications of symmetric and certain locally symmetric spaces.
\newblock {\em Geom. Topol.}, 22(5):2533--2646, 2018.

\bibitem[KLP17]{KLPanosovcharacterizations}
Michael Kapovich, Bernhard Leeb, and Joan Porti.
\newblock Anosov subgroups: Dynamical and geometric characterizations.
\newblock {\em European Journal of Mathematics}, 3:808--898, 2017.

\bibitem[KLP18]{KLP18:Morse}
Michael Kapovich, Bernhard Leeb, and Joan Porti.
\newblock A {M}orse lemma for quasigeodesics in symmetric spaces and {E}uclidean buildings.
\newblock {\em Geom. Topol.}, 22(7):3827--3923, 2018.

\bibitem[KN10]{KN10:HeisenbergGroup}
Tom Klein and Andrew Nicas.
\newblock The horofunction boundary of the {H}eisenberg group: the {C}arnot-{C}arath\'eodory metric.
\newblock {\em Conform. Geom. Dyn.}, 14:269--295, 2010.

\bibitem[Lab06]{Lab}
Fran\c{c}ois Labourie.
\newblock Anosov flows, surface groups and curves in projective space.
\newblock {\em Invent. Math.}, 165:51--114, 2006.

\bibitem[LOS78]{LOS78:Poulsen}
J.~Lindenstrauss, G.~Olsen, and Y.~Sternfeld.
\newblock The {P}oulsen simplex.
\newblock {\em Ann. Inst. Fourier (Grenoble)}, 28(1):vi, 91--114, 1978.

\bibitem[LS14]{LS14:Horoboundary_Teich}
Lixin Liu and Weixu Su.
\newblock The horofunction compactification of the {T}eichm\"uller metric.
\newblock In {\em Handbook of {T}eichm\"uller theory. {V}ol. {IV}}, volume~19 of {\em IRMA Lect. Math. Theor. Phys.}, pages 355--374. Eur. Math. Soc., Z\"urich, 2014.

\bibitem[LSZ15]{LiuSuZhong2015}
Lixin Liu, Weixu Su, and Youliang Zhong.
\newblock On metrics defined by length spectra on teichmüller spaces of surfaces with boundary.
\newblock {\em Annales Academiæ Scientiarum Fennicæ. Mathematica}, 40:617--644, 2015.

\bibitem[LW11]{LW11:Hilbert}
Bas Lemmens and Cormac Walsh.
\newblock Isometries of polyhedral {H}ilbert geometries.
\newblock {\em J. Topol. Anal.}, 3(2):213--241, 2011.

\bibitem[MGT26]{MGT25:Intersections}
D{\'\i}dac Mart{\'\i}nez-Granado and Dylan~P. Thurston.
\newblock The intersection dual of geodesic currents.
\newblock in preparation, 2026.

\bibitem[Min05]{Mineyev}
Igor Mineyev.
\newblock Flows and joins of metric spaces.
\newblock {\em Geom. Topol.}, 9:403--482, 2005.

\bibitem[MOT24]{MOT24:Ball}
Giuseppe Martone, Charles Ouyang, and Andrea Tamburelli.
\newblock A closed ball compactification of a maximal component via cores of trees.
\newblock {\em Algebr. Geom. Topol.}, 24(7):3693--3717, 2024.

\bibitem[MT18]{MT18:RandomWalks}
Joseph Maher and Giulio Tiozzo.
\newblock Random walks on weakly hyperbolic groups.
\newblock {\em J. Reine Angew. Math.}, 742:187--239, 2018.

\bibitem[MZ19]{MZ19:PositivelyRatioed}
Giuseppe Martone and Tengren Zhang.
\newblock Positively ratioed representations.
\newblock {\em Comment. Math. Helv.}, 94(2):273--345, 2019.

\bibitem[OR23]{Eduardo}
Eduardo Oreg\'{o}n-Reyes.
\newblock The space of metric structures on hyperbolic groups.
\newblock {\em J. Lond. Math. Soc. (2)}, 107(3):914--942, 2023.

\bibitem[OT21]{OT24:Blaschke}
Charles Ouyang and Andrea Tamburelli.
\newblock Limits of {B}laschke metrics.
\newblock {\em Duke Math. J.}, 170(8):1683--1722, 2021.

\bibitem[OT23]{OT23:SO23}
Charles Ouyang and Andrea Tamburelli.
\newblock Length spectrum compactification of the {$\rm SO_0(2,3)$}-{H}itchin component.
\newblock {\em Adv. Math.}, 420:Paper No. 108997, 37, 2023.

\bibitem[Ota90]{Otal90:SpectreMarqueNegative}
Jean-Pierre Otal.
\newblock Le spectre marqué des longueurs des surfaces à courbure négative.
\newblock {\em Ann. of Math. (2)}, 131(1):151, 1990.

\bibitem[Phe01]{phelps2001choquet}
Robert~R. Phelps.
\newblock {\em Lectures on Choquet's Theorem}, volume 1757 of {\em Lecture Notes in Mathematics}.
\newblock Springer, Berlin, 2nd edition, 2001.

\bibitem[Pou61]{poulsen1961simplex}
E.~B. Poulsen.
\newblock A simplex with dense extreme points.
\newblock {\em Annales de l'institut Fourier}, 11:83--87, 1961.

\bibitem[PS16]{PapadopoulosSu2016}
Athanase Papadopoulos and Weixu Su.
\newblock Thurston’s metric on teichmüller space and the translation distances of mapping classes.
\newblock {\em Annales Academiæ Scientiarum Fennicæ Mathematica}, 41(2):867--879, 2016.

\bibitem[PS24]{PanSu2024GeometryThurstonMetric}
Huiping Pan and Weixu Su.
\newblock The geometry of the thurston metric: A survey.
\newblock In Ken’ichi Ohshika and Athanase Papadopoulos, editors, {\em In the Tradition of Thurston III: Geometry and Dynamics}, pages 7--43. Springer, 2024.
\newblock First Online: 19 March 2024.

\bibitem[PT07]{GuPa}
Athanase Papadopoulos and Guillaume Th\'{e}ret.
\newblock On the topology defined by {T}hurston's asymmetric metric.
\newblock {\em Math. Proc. Cambridge Philos. Soc.}, 142(3):487--496, 2007.

\bibitem[Qui02]{Qui02:PS}
J.-F. Quint.
\newblock Mesures de {P}atterson-{S}ullivan en rang sup\'erieur.
\newblock {\em Geom. Funct. Anal.}, 12(4):776--809, 2002.

\bibitem[RS19]{RS19}
Kasra Rafi and Juan Souto.
\newblock Geodesic currents and counting problems.
\newblock {\em Geom. Funct. Anal.}, 29(3):871--889, 2019.

\bibitem[Sam14a]{HyperconvexRepsExponentialGrowth}
Andr\'{e}s Sambarino.
\newblock Hyperconvex representations and exponential growth.
\newblock {\em Ergodic Theory Dynam. Systems}, 34(3):986--1010, 2014.

\bibitem[Sam14b]{Quantitative}
Andr\'{e}s Sambarino.
\newblock Quantitative properties of convex representations.
\newblock {\em Comment. Math. Helv.}, 89(2):443--488, 2014.

\bibitem[Sam15]{Sam15:Orbital}
Andr\'es Sambarino.
\newblock The orbital counting problem for hyperconvex representations.
\newblock {\em Ann. Inst. Fourier (Grenoble)}, 65(4):1755--1797, 2015.

\bibitem[Sam24]{S24:Dichotomy}
Andr\'es Sambarino.
\newblock A report on an ergodic dichotomy.
\newblock {\em Ergodic Theory Dynam. Systems}, 44(1):236--289, 2024.

\bibitem[Sap24]{Jenya}
Jenya Sapir.
\newblock An extension of the {T}hurston metric to projective filling currents.
\newblock {\em Geom. Dedicata}, 218(3):Paper No. 78, 18, 2024.

\bibitem[Sas22]{Sas22:CuspedCurrents}
Dounnu Sasaki.
\newblock Currents on cusped hyperbolic surfaces and denseness property.
\newblock {\em Groups, Geometry, and Dynamics}, 16(3):1077--1117, 2022.

\bibitem[Shi25]{Shi25:Thurston}
Jiajun Shi.
\newblock Thurston's asymmetric metric for the space of flat metrics, 2025.

\bibitem[Tho19]{Tholozan2019CocyclesReparametrizations}
Nicolas Tholozan.
\newblock Cocycles, cohomology and reparametrizations of the geodesic flow.
\newblock Unpublished PDF; available at ENS, \url{https://www.math.ens.psl.eu/~tholozan/Annexes/CocyclesReparametrizations2.pdf}., 2019.

\bibitem[Thu98]{Th1998}
William~P. Thurston.
\newblock Minimal stretch maps between hyperbolic surfaces.
\newblock {\em arXiv:math GT/9801039}, 1998.

\bibitem[Tit71]{Tits}
Jacques Tits.
\newblock Représentations linéaires irréductibles d'un groupe réductif sur un corps quelconque.
\newblock {\em Journal für die reine und angewandte Mathematik (Crelle)}, 1971(247):196--220, 1971.

\bibitem[Tri22]{T22:ThurstonCompactification}
Marie Trin.
\newblock Thurston's compactification via geodesic currents: The case of non-compact finite area surfaces, 2022.

\bibitem[Wal14a]{CW2015}
Cormac Walsh.
\newblock The horoboundary and isometry group of {T}hurston's {L}ipschitz metric.
\newblock In {\em Handbook of {T}eichm\"{u}ller theory. {V}ol. {IV}}, volume~19 of {\em IRMA Lect. Math. Theor. Phys.}, pages 327--353. Eur. Math. Soc., Z\"{u}rich, 2014.

\bibitem[Wal14b]{Co14:FiniteHilbert}
Cormac Walsh.
\newblock The horofunction boundary and isometry group of the {H}ilbert geometry.
\newblock In {\em Handbook of {H}ilbert geometry}, volume~22 of {\em IRMA Lect. Math. Theor. Phys.}, pages 127--146. Eur. Math. Soc., Z\"urich, 2014.

\bibitem[Wal18]{Cor18:HilbertInfinite}
Cormac Walsh.
\newblock Hilbert and {T}hompson geometries isometric to infinite-dimensional {B}anach spaces.
\newblock {\em Ann. Inst. Fourier (Grenoble)}, 68(5):1831--1877, 2018.

\end{thebibliography}

\end{document}